\input amstex.tex
\input amsppt.sty

%\topmatter
%\title{Extensions of valuation rings}\endtitle
%\author{Shai Sarussi}\endauthor
%\address{Department of Mathematics, Bar-Ilan University, Ramat-Gan 52900,
%Israel}\endaddress

%\endtopmatter

%\footline{\ifnum\pageno=0\hfill\else\hss\tenrm\folio\hss\fi}
%\topinsert\vskip1.8truecm\endinsert
 \centerline{\bf Quasi-Valuations Extending a Valuation } \vskip8pt
$${\vbox{\halign{\hfil\hbox{#}\hfil\qquad&\hfil\hbox{#}\hfil\cr
$$Shai Sarussi$^*$\cr
Department of Mathematics\cr Bar Ilan University\cr Ramat-Gan
52900, Israel\cr}}}$$

\topmatter \abstract {Suppose $F$ is a field with valuation $v$
and valuation ring $O_{v}$, $E$ is a finite field extension and
$w$ is a quasi-valuation on $E$ extending $v$. We study
quasi-valuations on $E$ that extend $v$; in particular, their
corresponding rings and their prime spectrums. We prove that these
ring extensions satisfy INC (incomparability), LO (lying over),
and GD (going down) over $O_{v}$; in particular, they have the
same Krull Dimension. We also prove that every such
quasi-valuation is dominated by some valuation extending $v$.

Under the assumption that the value monoid of the quasi-valuation
is a group we prove that these ring extensions satisfy GU (going
up) over $O_{v}$, and a bound on the size of the prime spectrum is
given. In addition, a 1:1 correspondence is obtained between
exponential quasi-valuations and integrally closed quasi-valuation
rings.

Given $R$, an algebra over $O_{v}$, we construct a quasi-valuation
on $R$; we also construct a quasi-valuation on $R \otimes _{O_{v}}
F$ which helps us prove our main Theorem. The main Theorem states
that if $R \subseteq E$ satisfies $R \cap F=O_{v}$ and $E$ is the
field of fractions of $R$, then $R$ and $v$ induce a
quasi-valuation $w$ on $E$ such that $R=O_{w}$ and $w$ extends
$v$; thus $R$ satisfies the properties of a quasi-valuation ring.

%Suppose $F$ is a field with valuation $v$ and valuation ring
%$O_{v}$, $E$ is a finite field extension, and $w$ is a
%quasi-valuation on $E$ extending $v$. We study quasi-valuations on
%$E$ that extend $v$. We prove that every such quasi-valuation is
%dominated by some valuation extending $v$. Introducing a notion
%called PIM, we prove that the quasi-valuation
%rings satisfy GD (going down) over $O_{v}$. %and conclude that these
%rings satisfy the height formula.

}\endabstract \thanks {$^*$This paper is part of the author's
forthcoming doctoral dissertation.} \endthanks
\endtopmatter

\heading \S 0 Introduction
\endheading

Recall that a valuation on a field $F$ is a function $v : F
\rightarrow \Gamma \cup \{  \infty \}$, where $\Gamma$ is a
totally ordered abelian group and $v$ satisfies the following
conditions:

(A1) $v(0) = \infty$;

 (A1') $v(x) \neq \infty$ for every $0 \neq x\in F$;

(A2) $v(xy) = v(x)+v(y)$ for all $x,y \in F$;

(A3) $v(x+y) \geq \min \{ v(x),v(y) \}$ for all $x,y \in F$.

   Valuation theory has long been a key tool in commutative algebra,
with applications in number theory and algebraic geometry. It has
become a useful tool in the study of division algebras, and used
in the construction of various counterexamples such as Amitsur's
construction of noncrossed products division algebras. See [Wad]
for a comprehensive survey.

Generalizations of the notion of valuation have been made
throughout the last few decades.
%Morandi's value functions (cf. [Mor]) and Tignol and Wadsworth's
%gauges [TW] are examples of such generalizations. Also see [KZ]
%for Manis valuations and PM-valuations. Let us review some of
%these theories.
Knebusch and Zhang (cf. [KZ]) have studied valuations in the sense
of Bourbaki [Bo, section 3]. Thus they were able, omitting (A1'),
to study valuations on any commutative ring rather than just on an
integral domain. They define the notion of Manis valuation (the
valuation is onto the value group) and show that an $R$-Pr\"ufer
ring is related to Manis valuations in much the same way that a
Pr\"ufer domain is related to valuations of its quotient field.

%Pseudo-valuations (whose axioms are the closest to the axioms of
% quasi-\break valuations) were also studied as generalizations of
%valuations. A Pseudo-valuation satisfy the same conditions as
%quasi-valuation but its values lie inside a totally ordered
%abelian group $\Gamma$; actually, much of the work on
%pseudo-valuations was by taking $\Gamma=\Bbb R$. As examples one
%can mention Cohn (cf. [Co]) who has given necessary and sufficient
%conditions for a non-discrete topological field to have its
%topology induced by a pseudo-valuation; Huckaba (cf. [Hu]) has
%given, using filtration, necessary and sufficient conditions for a
%pseudo-valuation to be extended to an overring and
%Mahdavi-Hezavehi has obtained matrix valuations from matrix
%pseudo-valuations [cf. MH].

A monoid $M$ is called a {\it totally ordered} monoid if it has a
total ordering $\leq$, for which $a \leq b$ implies $a+c \leq b+c$
for every $c \in M$. When we write $b > a $ we mean $b \geq a$ and
$b \neq a$.

In this paper we study quasi-valuations, which are generalizations
of valuations. A {\it quasi-valuation} on a ring $R$ is a function
$w : R \rightarrow ~M \cup \{  \infty \}$, where $M$ is a totally
ordered abelian monoid, to which we adjoin an element $\infty$
greater than all elements of $M$, and $w$ satisfies the following
properties:

(B1) $w(0) = \infty$;

(B2) $w(xy) \geq w(x) + w(y)$ for all $x,y \in R$;

(B3) $w(x+y) \geq \min \{ w(x), w(y)\}$ for all $x,y \in R$.

In the literature this is called a pseudo-valuation when the
target $M$ is a totally ordered abelian group, usually taken to be
$(\Bbb R,+)$. As examples one can mention Cohn (cf. [Co]) who gave
necessary and sufficient conditions for a non-discrete topological
field to have its topology induced by a pseudo-valuation; Huckaba
(cf. [Hu]) has given necessary and sufficient conditions for a
pseudo-valuation to be extended to an overring, and
Mahdavi-Hezavehi has obtained "matrix valuations" from matrix
pseudo-valuations ([cf. MH]). We use the terminology
quasi-valuation to stress that the target monoid need not be a
group. Moreover, our study concentrates on quasi-valuations
extending a given valuation on a field.

The minimum of a finite number of valuations with the same value
group is a quasi-valuation. For example, the $n$-adic
quasi-valuation on $\Bbb Q$ (for any positive $n \in \Bbb Z$)
already has been studied in [Ste]. (Stein calls it the $n$-adic
valuation.) It is defined as follows: for any $0 \neq \frac{c}{d}
\in \Bbb Q$ there exists a unique $e \in \Bbb Z$ and integers $a,b
\in \Bbb Z$, with $b$ positive, such that
$\frac{c}{d}=n^{e}\frac{a}{b}$ with $n \nmid a$, $(n,b)=1$ and
$(a,b)=1$. Define $w_{n}(\frac{c}{d})=e$ and $w_{n}(0)=\infty$.

%Let $R$ be a UFD, let $F$ be its field of fractions, and fix a
%positive $n \in R$ (not necessarily prime). Every $0 \neq
%\frac{c}{d} \in F$ can be written uniquely
%$\frac{c}{d}=n^{e}\frac{a}{b}$, where $a,b \in R$, $b>0$, $n \nmid
%a$, $(n,b)=1$ and $(a,b)=1$. Define $w_{n}(\frac{c}{d})=e$ and
%$w_{n}(0)=\infty$. $w_{n}$ is called the n-adic quasi-valuation on
%$F$. A detailed analysis of the n-adic quasi-valuation on $\Bbb Q$
%is given in [Ste, ch. 16].

A quasi-valuation is a much more flexible tool than a valuation;
for example, quasi-valuations exist on rings on which valuations
cannot exist.

Three main classes of rings were suggested throughout the years as
the noncommutative version of a valuation ring. These three types
are invariant valuation rings, total valuation rings, and Dubrovin
valuation rings. They are interconnected by the following diagram:
$$\{\text{invariant valuation rings}\}\!\subset\!\{\text {total valuation rings} \}
\!\subset\! \{\text {Dubrovin valuation rings}\}.$$

Morandi (cf. [Mor]) has studied Dubrovin valuation rings and their
ideals. It turned out that unlike a valuation on a field, the
value group of a Dubrovin valuation ring $B$ does not classify the
ideals in general but does so when $B$ is integral over its
center.

At the outset and in section 1 one sees that there are an enormous
amount of quasi-valuations, even on $\Bbb Z$, so in order to
obtain a workable theory one needs further assumptions. Morandi
(cf. [Mor]) defines a value function which is a quasi-valuation
satisfying a few more conditions. Given an integral Dubrovin
valuation ring $B$ of a central simple algebra $S$, Morandi shows
that there is a value function $w$ on $S$ with $B$ as its value
ring (the value ring of $w$ is defined as the set of all $x \in S$
such that $w(x) \geq 0$). Morandi also proves the converse, that
if $w$ is a value function on $S$, then the value ring is an
integral Dubrovin valuation ring.

The use of value functions allow a number of results about
invariant valuation rings to be extended to Dubrovin valuation
rings. Also, their use has led to simpler and more natural proofs
of a number of results on invariant valuation rings.

Tignol and Wadsworth (cf. [TW]) have developed a powerful theory
which utilizes filtrations. They consider the notion of gauges,
which are the surmultiplicative value functions for which the
associated graded algebra is semisimple, and which also satisfy a
defectlessness condition. The gauges are defined on finite
dimensional semisimple algebras over valued fields with arbitrary
value groups. Tignol and Wadsworth also show the relation between
their value functions and Morandi's.

Quasi-valuations generalize both value functions and gauges (the
axioms of quasi-valuations are contained in the axioms of value
functions and gauges). Although we do not obtain Tignol and
Wadsworth's decisive results concerning the Brauer group, we do
get a workable theory, in which we are able to answer questions
regarding the structure of rings using quasi-valuation theory. One
difference with value functions, for example, is that on a field,
Morandi's construction of a value function is automatically a
valuation, while that is not the case for quasi-valuations. Gauges
on a field reduce to exponential quasi-valuations, a special case
of quasi-valuations. We believe that, even over a field,
quasi-valuations are a natural generalization of valuations, and
enrich valuation theory even for commutative algebras.

Whereas for a valuation on a field $M$ is automatically a group
since \break $v(x^{-1})=-v(x)$, the situation is different with
quasi-valuations (since \break $w(x^{-1}) + w(x) \leq w(1)$). Thus
it is more natural to map $w$ to a monoid. In case $w : R
\rightarrow M \cup \{ \infty \}$ is a quasi-valuation for $M$ a
group, we say that $w$ is {\it group-valued} (rather than calling
it a pseudo-valuation); we do not require $w$ to be surjective.

 Note that it does not follow from the axioms of a
quasi-valuation that $w(1)$ is necessarily 0. We shall see
examples in which $w(1) \neq 0$.

We list here some of the common symbols we use for $v$ a valuation
on a field $F$ and $w$ a quasi-valuation on a ring $R$ (usually we
write $E$ instead of $R$):

$O_{v}=\{x\in F \ | \ v(x)\geq0\}$; the valuation ring.

$I_{v}=\{x\in F \ | \ v(x)>0\}$; the valuation ideal.

$O_{w}=\{x\in R \ | \ w(x)\geq0\}$; the quasi-valuation ring.

$I_{w}=\{x\in R \ | \ w(x)>0\}$; the quasi-valuation ideal.

$J_{w}$; the Jacobson radical of $O_{w}$.

$\Gamma_{v}$; the value group of the valuation $v$.

$M_{w}$; the value monoid of the quasi-valuation $w$, i.e., the
submonoid of $M$ generated by $  w(R\setminus \{ 0 \})$.

Note that $w$ is group-valued iff $M_{w}$ is cancellative, since
any ordered abelian cancellative monoid has a group of fractions.

Here is a brief overview of this paper. %This paper is organized as
In section 1 we present some general properties of
quasi-valuations on rings as well as some basic examples.
Thereafter we work under the assumption that $F$ denotes a field
with a valuation $v$, $E/F$ usually is a finite field extension,
and $w$ is a quasi-valuation on $E$ such that $w|_{F}=v$. In
section 2 we discuss some of the basic results regarding
quasi-valuations and their corresponding rings; most of the
results in this section are valid in the more general case where
$E$ is a finite dimensional $F-$algebra. In section 3 we prove
that a quasi-valuation ring satisfies INC and LO over $O_{v}$; in
fact, LO is valid in the case where $E$ is a finite dimensional
$F-$algebra. In section 4 we introduce a notion called PIM
(positive isolated monoid); the PIMs enable us to prove that a
quasi-valuation ring satisfies GD over $O_{v}$. We also describe a
connection between the closure of a PIM and the prime ideals of
$O_{w}$. In sections 5 and 6 we assume that $w(E \setminus \{ 0
\})$ is torsion over $\Gamma_{v}$, the value group of the
valuation. In section 5 we generalize some properties of valuation
rings and we discuss the set of all expansions of a
quasi-valuation ring (and one of his maximal ideals). We show that
a quasi-valuation ring satisfies GU over $O_{v}$. We also
construct a quasi-valuation that arises naturally from $w$; this
quasi-valuation is one of the key steps to give a bound on the
size of the prime spectrum of the quasi-valuation ring. In section
6 we prove that any valuation whose valuation ring contains
$O_{w}$ dominates $w$. We also show that $O_{w}$ satisfies the
height formula (since we have the property GU). In section 7 we
discuss exponential quasi-valuations; we do not assume that
$w(E\setminus \{ 0 \})$ is torsion over $\Gamma_{v}$. Instead, we
assume the weaker hypothesis that the value monoid is weakly
cancellative. We prove that an exponential quasi-valuation must be
of the form $w= \min \{u_1,...,u_k\}$ for valuations $u_i$ on $E$
extending $v$. We obtain a 1:1 correspondence between $ \{ \text
{exponential quasi-valuations extending\ } v\}$ and $ \{
\text{integrally closed quasi-valuation rings}\}$. We also deduce
that the number of exponential quasi-valuations is bounded by
$2^{[E:F]}-1$. In section 8 we construct a total ordering on a
suitable amalgamation of $M_{w}$ and $\Gamma_{\text{div}}$ that
allows us to compare elements of $\Gamma_{\text{div}}$ with
elements of $M_{w}$. Then we show that there exists a valuation
$u$ whose valuation ring contains $O_{w}$ and $u$ dominates $w$.
In section 9 we review some of the notions of cuts and we present
the construction of the cut monoid (of a totally ordered abelian
group) which is an $\Bbb N$-strictly ordered abelian monoid. Then
we introduce the filter quasi-valuations, which are
quasi-valuations that can be defined on any $O_{v}-$algebra. The
filter quasi-valuation is induced by an $O_{v}-$algebra and the
valuation $v$; its values lie inside the cut monoid of
$\Gamma_{v}$. This gives us our main theorem, which enables us to
apply the methods developed in the previous sections, as indicated
below in Theorem 9.37. In section 10 we show that filter
quasi-valuations and their cut monoids satisfy some properties
which general quasi-valuations and their monoids do not
necessarily share. This enables us to prove a stronger version of
Theorem 8.15 when dealing with filter quasi-valuations. We also
show that the filter quasi-valuation construction respects
localization at prime ideals of $O_{v}$. Finally, we present the
minimality of the filter quasi-valuation with respect to a natural
partial order.   
%we define an equivalence relation for quasi-valuations extending
%$v$ on $E$ having the same value ring; we also define a partial
%order on $\Cal W_{R} / \sim$. Then, we prove that the equivalence
%class of the filter quasi-valuation is the minimal one.
%Finally, we present a minimality property of the filter
%quasi-valuation.

%We also show that the equivalence class of the
%filter quasi-valuation is the minimal one with respect to a
%natural partial order.

%This paper is organized as follows: in section 1 we give some
%basic definitions and results following [Sa]. In section 2 we show
%that any quasi-valuation is dominated by some valuation. In
%section 3 we introduce a notion called PIM (positive isolated
%monoid) and prove that the quasi-valuation ring $O_{w}$ satisfies
%GD over the valuation ring $O_{v}$. In section 4 we construct an
%$\Bbb N$-strictly ordered abelian monoid; called the cuts
%monoid and introduce the filter quasi-valuations. This gives us
%our main theorem,
%our main theorem is that for any field $F$ with valuation $v$, any
%finite field extension $E$ of $F$, and any ring $R \subseteq E$
%whose field of fractions is $E$ and $R \cap F=O_{v}$, there is a
%quasi-valuation $w$ of $E$ extending $v$ such that $R=O_{w}$.
%which enables us to apply the methods of [Sa] as indicated below
%in Theorem 4.31. In section 5 we show that filter quasi-valuations
%and their cut monoids satisfy some properties which
%general quasi-valuations and their monoids do not necessarily
%share.

\heading \S 1 General Quasi-Valuations and Examples
\endheading

In this section we present some of the basic definitions and
properties regarding quasi-valuation on rings. We also present
some examples of quasi-valuations. In particular, we give examples
of quasi-valuations on integral domains which cannot be extended
to their fields of fractions.

\demo {Remark 1.1} If the monoid $M$ is cancellative then $w(1)
\leq 0$.
\enddemo

\demo {Proof} $w(1)=w(1^{2})\geq w(1)+w(1)$. \enddemo

Here is a trivial example of a quasi-valuation on $\Bbb Z$.

\demo {Example 1.2} $w(0)=\infty$ and $w(z)=-1$ for all $z \neq
0$. \enddemo

 \proclaim {Lemma 1.3} Let $w$ be a quasi-valuation
on a ring $R$ and suppose that $w(-1)=0$. Then $w(a)=w(-a)$ for
any $a \in R$. \endproclaim

\demo {Proof} Let $a \in R$; then $$w(-a)=w(-1 \cdot a)\geq
w(-1)+w(a)=w(a).$$ By symmetry $w(a)=w(-(-a))\geq w(-a)$.\enddemo

The following lemma generalizes a well known lemma from valuation
theory, with the same proof. We prove it here for the reader's
convenience.

 \proclaim {Lemma 1.4} Let $w$ be a quasi-valuation
on a ring $R$ and suppose that $w(-1)=0$. Let $a,b \in R$;

$(i)$ If $w(a) \neq w(b)$ then $w(a+b)=\min \{ w(a), w(b)\}.$

$(ii)$ If $w(a+b)>w(a)$ then $w(b)=w(a)$.

\endproclaim

\demo {Proof} We prove the first statement; the second follows
easily. Assume $w(a)<w(b)$; then $w(a+b) \geq w(a)$. On the other
hand,
$$w(a)=w(a+b-b) \geq \min \{ w(a+b), w(b)\}=w(a+b).$$

Q.E.D.
\enddemo

%\proclaim {Lemma 0.5} If $w$ is a quasi-valuation on $R$ with
%$w(-1)=0$ and $a,b \in ~R$ such that $w(a+b)>w(a)$, then
%$w(b)=w(a)$.\endproclaim

%\demo {Proof} Follows at once from Lemma 0.4.
%\enddemo

In order to be able to use techniques from valuation theory, we
define the following special elements with respect to a
quasi-valuation.

\demo {Definition 1.5} Let $w$ be a quasi-valuation on a
commutative ring $R$. An element $c \in R$ is called {\it stable}
with respect to $w$ if $$w(cx)=w(c)+w(x)$$ for every $x \in R$.

\enddemo

\proclaim {Lemma 1.6} Let $w$ be a quasi-valuation on a
commutative ring $R$ such that %$M_{w}$ is weakly cancellative and
$w(1)=0$. Let $c$ be an invertible element of $R$. Then $c$ is
stable iff $w(c)=-w(c^{-1})$. \endproclaim

\demo {Proof} ($\Rightarrow$) $c$ is stable and invertible; thus
$$0=w(1)=w(cc^{-1})=w(c)+w(c^{-1}),$$ i.e., $w(c)=-w(c^{-1})$.

\ \ \ \ \ \ ($\Leftarrow$) Let $x \in R$; we have
$$w(x)=w(c^{-1}cx)\geq w(c^{-1})+w(cx)$$ $$\geq
w(c^{-1})+w(c)+w(x)=w(x).$$ Hence equality holds, and since
$w(c^{-1})=-w(c)$ is invertible, we have $$w(cx)=w(c)+w(x).$$
%First note that since $M_{w}$ is weakly cancellative, $w(1)=0$
%(indeed w(1)=w(1)+w(1))

Q.E.D.
\enddemo

\demo {Definition 1.7}
%A {\it quadratic quasi-valuation} $w$ on
%a ring $R$ is a quasi-valuation such that for every $x \in R$,
%$w(x^{2})=2w(x)$.
An {\it exponential quasi-valuation} $w$ on a ring $R$ is a
quasi-valuation such that for every $x \in R$,
$$(1) \ \ \ \   w(x^{n})=nw(x), \ \ \ \forall n \in \Bbb N.$$
\enddemo

We note that any exponential quasi-valuation $w$ with $M_{w}$
cancellative satisfies $w(1)~=w(-1)=0$.

Usually we need only check the following special case of (1):

\proclaim {Proposition 1.8} Let $R$ be a ring and let $M$ be a
totally ordered cancellative abelian monoid. Let $w : R
\rightarrow M \cup \{ \infty \}$ be a quasi-valuation satisfying
the condition: for every $x \in R$, $w(x^{2})=2w(x)$. Then $w$ is
an exponential quasi-valuation.
\endproclaim

\demo {Proof} It is obvious that $w(x^{n})=nw(x)$ for $n=1,2$.
Assume by induction that the statement is true for every $k<n$.

If $n$ is even,

$$w(x^{n})=w((x^{n/2})^{2})=2w(x^{n/2})=2(n/2)w(x)=nw(x).$$

If $n$ is odd,

$$w(x^{n})\geq nw(x)$$ and

$$(n+1)w(x)=w(x^{n+1})=w(x^{n}x)$$$$\geq w(x^{n})+w(x)\geq nw(x)+w(x)=(n+1)w(x).$$

Therefore equality holds and cancelling $w(x)$, we get
$w(x^{n})=nw(x)$.

Q.E.D.
\enddemo

We present here some examples of quasi-valuations:

%for any $0 \neq \frac{c}{d} \in \Bbb Q$ there exists a unique $e
%\in \Bbb Z$ and integers $a,b \in \Bbb Z$, with $b$ positive, such
%that $\frac{c}{d}=n^{e}\frac{a}{b}$ with $n \nmid a$, $(n,b)=1$
%and $(a,b)=1$. Define $w_{n}(\frac{c}{d})=e$ and
%$w_{n}(0)=\infty$.

\demo{Example 1.9} Let $R$ be a UFD, let $F$ be its field of
fractions, and fix a non-zero $n \in R$ (not necessarily prime).
For any $0 \neq \frac{c}{d} \in F$ there exists a unique $e \in R$
such that $\frac{c}{d}=n^{e}\frac{a}{b}$, where $a,b \in R$, $n
\nmid a$, $(n,b)=1$ and $(a,b)=1$. Define $w_{n}(\frac{c}{d})=e$
and $w_{n}(0)=\infty$. $w_{n}$ is called the n-adic
quasi-valuation on $F$. A detailed analysis of the n-adic
quasi-valuation on $\Bbb Q$ is given in [Ste, ch. 16].

%\footnote{*This paper is part of the author's forthcoming doctoral
%dissertation.}

%\vfill\eject

Note that whenever $$n=\prod_{i=1}^r p_{i}$$ where the $p_{i}$'s
are distinct non-zero primes, the n-adic quasi-valuation $w_{n}$
is equal to $$\min_{1 \leq i \leq r } \{w_{p_{i}}\};$$ where each
$w_{p_{i}}$ is a valuation.\enddemo

\demo{Example 1.10} Let $R$ be a ring and let $\{w_{i}\}_{i \in I
}$ be a set of quasi-valuations having the same value monoid
(i.e., for every $i,j \in I$, $M_{w_{i}}=M_{w_{j}}=:M$ ). Assume
that for any $r \in R$, $\{ w_{i}(r) \mid \ i \in I \}$ is finite.
This is automatic when $I$ is finite but also holds for an
arbitrary set of valuations over a number field (i.e., the
$w_{i}$'s are valuations) since for every $r$, $w_{i}(r)=0$ for
almost all $i$. Define $$u=\min_{i \in I}\{w_{i} \}.$$ Then $u$ is
a quasi-valuation with values inside the monoid $ M$. \enddemo

Example 1.10 will be generalized and become a tool for computing
Krull dimension in section 5.

  Let $R$ be an integral domain and assume that $v$ is a valuation
on $R$. Then $v$ is automatically a valuation on $F$, the field of
fractions of $R$. However, the situation is different in the
quasi-valuation case. A quasi-valuation on an integral domain
cannot always be extended to its field of fractions. For example,
let $v$ be a valuation on a field $F$ with value group $\Gamma$
and let $\alpha \in \Gamma$ be a positive element. Define
$w_{\alpha}:F \rightarrow \Gamma \cup \{ \infty \}$ by
$$ w_{\alpha}(x)=\cases v(x)  \text{ \ \  if     } v(x) < \alpha \\  \infty \text{\ \ \ \ \ otherwise.} \endcases$$
$w_{\alpha}$ is a quasi-valuation on $O_{v}$ but obviously cannot
be extended to $F$ (since $w^{-1}(\{ \infty \})$ is an ideal of
$F$).

We shall now present an example of a quasi-valuation on $\Bbb Z$
with value group $\Gamma=\Bbb Z$ such that $w(x) \neq \infty$ for
$x \neq 0$ and for which the quasi-valuation cannot be extended to
a quasi-valuation on its field of fractions $\Bbb Q$, such that
the value monoid is torsion over $\Gamma$.

\demo{Example 1.11} Let $v$ be any p-adic valuation on $\Bbb Z$;
we define a function $w: \Bbb Z \rightarrow \Bbb Z \cup \{ \infty
\}$ by $w(x)=[v(x)]^{2}$ (with the convention
$\infty^{2}=\infty$). Note that im$(w) \subseteq \Bbb N \cup \{ 0,
\infty \}$; it is not difficult to see that $w$ is a
quasi-valuation on $\Bbb Z$. Indeed, $w(0)=\infty$; if $x,y \in
\Bbb Z$, then
$$w(xy)=[v(xy)]^{2}=[v(x)+v(y)]^{2}$$$$=[v(x)]^{2}+2v(x)v(y)+[v(y)]^{2}
\geq w(x)+w(y).$$ Finally, if $x,y \in \Bbb Z$, then
$$w(x+y)=[v(x+y)]^{2} \geq [\min \{v(x),v(y)\}]^{2}$$ $$=\min \{[v(x)]^{2},[v(y)]^{2}\}= \min \{w(x),w(y) \}.$$ Now, let $x,y$ be non-zero integers
such that $w(x), w(y)>0$ and let $n \in \Bbb N$; we have

$$ w(y ^{n})= w(y ^{n}  x  \frac{1}{x}) \geq w(y ^{n}  x)+w(\frac{1}{x}) $$

$$= [v(y ^{n}  x)]^{2}+ w(\frac{1}{x})= [v(y ^{n}) + v(x)]^{2}+ w(\frac{1}{x})$$$$=
[v(y ^{n})]^{2}+2v(y ^{n})v(x) + [v(x)]^{2}+
w(\frac{1}{x})$$$$=w(y ^{n})+2v(y ^{n})v(x) + w(x)+
w(\frac{1}{x}).$$

Thus, cancelling $w(y ^{n})$, we get $w(\frac{1}{x}) \leq -2v(y
^{n})v(x)-w(x)$ for every $n \in \Bbb N$. Hence $w(\frac{1}{x})$
is less than every $z \in \Bbb Z$. Thus there cannot exist an $n
\in \Bbb N$ such that $nw(\frac{1}{x}) \in \Bbb Z$. Therefore $w$
cannot be extended to $\Bbb Q$ with value monoid torsion over
$\Gamma$.\enddemo

Let us consider an example of such a quasi-valuation extending a
valuation in an infinite dimensional extension.

\demo{Example 1.12} Let $v$ denote the trivial valuation on $\Bbb
Q$ and let $w$ denote the $\lambda ^{2}$-adic quasi-valuation (see
Example 1.9) on $\Bbb Q[\lambda]$, where $\lambda$ is a
commutative indeterminate. $w$ is a quasi-valuation which is not a
valuation since $w(\lambda ^{2})=1>0=2w(\lambda)$ and obviously
$w$ extends $ v$.\enddemo

In order to avoid pathological cases, we shall study
quasi-valuations extending a given valuation on a field.

{\it In this paper, $F$ denotes a field with a valuation $v$,
$E/F$ usually is a finite field extension with $n=[E:F]$, and $w$
is a quasi-valuation on $E$ such that $w|_{F}=v$.} We shall see
that the theory is surprisingly rich in these cases, especially
for exponential quasi-valuations.

Some of the results are valid in the more general case where $F$
is a field and $E$ is a finite dimensional $F$-algebra. When
possible, we shall discuss this more general scope.

%We present here two examples of quasi-valuations: 1. Let $R$ be a
%UFD and let $F$ be its field of fractions. Every $0 \neq
%\frac{c}{d} \in F$ can be written uniquely
%$\frac{c}{d}=n^{e}\frac{a}{b}$ where $n \nmid a$, $(n,b)=1$ and
%$(a,b)=1$. Define $w(\frac{c}{d})=e$, $w(0)=\infty$. $w$ is called
%the n-adic quasi-valuation on $F$. A detailed analysis of the
%n-adic quasi-valuation on $\Bbb Q$ is given in [Ste, ch. 16]. 2.
%If $v_{1}$ and $v_{2}$ are two distinct discrete valuations on a
%field, then $w= \min \{ v_{1},v_{2} \}$ is a quasi-valuation which
%is not a valuation. This example will be generalized and become a
%tool for computing Krull dimension in section 2.

%\footnote{*This paper is part of the author's forthcoming doctoral
%dissertation.}

%\vfill\eject

\heading \S 2 Basic Properties and Examples
\endheading

In this section, until Example 2.7, we consider the more general
case for which $E$ is a (not necessarily commutative) finite
dimensional $F$-algebra.

Note that $\Gamma_{v}$ embeds in $M_{w}$ since $w$ extends $v$.
Also note that since $w$ extends $v$ and $v$ is a valuation, then
by Lemma 1.6 every nonzero $x\in F$ is stable with respect to $w$;
$0$ is obviously stable. So, every $x\in F$ is stable with respect
to $w$. Furthermore, $w(-1)=v(-1)=0$. Hence, applying induction to
Lemma 1.4, for any $x_1,x_2,...x_n$ such that there exists a
single $k$ satisfying $w(x_k)=\min_{1 \leq i \leq n}\{
w(x_{i})\}$, we have
$$ w(\sum_{i=1}^ {n} x_{i})=w(x_{k}) .$$

\proclaim {Lemma 2.1} The elements of $w(E\setminus \{ 0 \})$ lie
in $\leq n$ cosets over $\Gamma_{v}$, where we recall
$[E:F]=n$.\endproclaim

\demo {Proof} Let $b_{1},...,b_{m}\in E\setminus \{ 0 \}$ such
that $w(b_{i})+\Gamma_{v}\neq w(b_{j})+\Gamma_{v}$ for every
$i\neq j$. We show that the set $\{b_{k}\}^{m}_{k=1}$ is linearly
independent. Assume to the contrary, that there exist $\alpha_{k}
\in F$ not all $0$ such that $\sum \alpha_{i} b_{i}=0$. Note that
$w(\alpha_{i} b_{i})=w(\alpha_{i})+w( b_{i})$ and therefore the
value of each $\alpha_{i} b_{i}$ lies in a different coset; in
particular, $w(\alpha_{i} b_{i}) \neq w(\alpha_{j} b_{j})$ for
every $i\neq j$. Hence $$w(\sum \alpha_{i} b_{i})=\min_{1 \leq i
\leq n} \{w(\alpha_{i} b_{i})\}\neq w(0),$$ a contradiction.

Q.E.D.\enddemo

\demo {Remark 2.2} Let $O_{v}$ be a valuation ring of $F$ and let
$E$ be an $F$-algebra with $[E:F]=n$. Then every set
$\{b_{1},b_{2},...,b_{m}\} \subseteq ~E$ with $m>n$ satisfies the
following: there exist $\alpha_{i}\in O_{v}$ such that
$\sum_{i=1}^{m} \alpha_{i}b_{i}=0$
 and $\alpha_{i_{0}}=1$ for some $1 \leq i_{0}
\leq m$, i.e., $$b_{i_{0}}=-\sum_{i\neq i_{0}} \alpha_{i}b_{i}.$$
% $b_{i_{0}}=-\sum_{i\neq i_{0}} \alpha_{i}b_{i}$ for some $1 \leq
% i_{0} \leq m$ and $\alpha_{i}\in F$.
\enddemo

\demo {Proof} The $b_{i}$ are linearly dependent over $F$, so
there exist $\beta_{i}\in F$ not all $0$ such that $\sum_{i=1}^{m}
\beta_{i}b_{i}=0$. We choose a $\beta_{i_{0}}$ such that
$v(\beta_{i_{0}})\leq v(\beta_{i})$ for all $1\leq i \leq m$ and
multiply by $\beta_{i_{0}}^{-1}$. Writing $\alpha_{i}=
\beta_{i_{0}}^{-1}\beta_{i}$, we have $\alpha_{i}\in O_{v}$ and
$\alpha_{i_{0}}=1$.

Q.E.D.\enddemo

\proclaim {Lemma 2.3} Let $R \supseteq O_{v}$ be a ring contained
in $E$. Then every f.g. $O_{v}$-submodule of $R$ is spanned by
$\leq n$ generators.\endproclaim

\demo {Proof} By Remark 2.2 every set
$\{b_{1},b_{2},...,b_{n+1}\}$ of $n+1$ elements can be reduced to
a set of $n$ elements which span the same submodule.\enddemo

\proclaim{Corollary 2.4} Every f.g. ideal of $O_{w}$ is generated
by $\leq n$ generators.\endproclaim

Note that if an ideal $I$ is not f.g., one cannot argue that $I$
is generated by $\leq n$ generators. Indeed, let $v$ be a
non-discrete valuation; thus $O_{v}$ itself is not Noetherian.
Now, take $E=F$ and $w=v$.

%The following lemma is valid in the more general case where $E$ is
%a finite dimensional $F$-algebra.

\proclaim {Lemma 2.5} $[O_{w}/I_{v}O_{w}:O_{v}/I_{v}]\leq
n$.\endproclaim

\demo {Proof} Let $\overline{a_{1}},...,\overline{a_{m}}\in
O_{w}/I_{v}O_{w}$ be linearly independent over $O_{v}/I_{v}$. Let
$a_{1},...,a_{m}$ be their representatives. We prove that
$a_{1},...,a_{m}$ are linearly independent over $F$. Assume to the
contrary and apply Remark 2.2 to get $\sum_{i=1}^m
\alpha_{i}a_{i}=0$ where $\alpha_{i} \in O_{v}$ and
$\alpha_{i_{0}}=1$ for some $1 \leq i_{0} \leq m$. So, we have
$\sum_{i=1}^m \overline{\alpha_{i}} \overline{a_{i}}=0$ where
$\overline{\alpha_{i_{0}}}=\overline{1}$. This contradicts the
linear independence of the $\overline{a_{i}}$'s.

Q.E.D.\enddemo

\proclaim{Corollary 2.6} $[O_{w}/I_{w}:O_{v}/I_{v}]\leq
n$\endproclaim

%$[O_{w}/J_{w}:O_{v}/I_{v}]\leq n$ .
\demo {Proof} Note that $I_{w} \supseteq I_{v}O_{w}$ and thus we
have the natural epimorphism $$O_{w}/I_{v}O_{w} \twoheadrightarrow
O_{w}/I_{w}.$$ The result now follows from Lemma 2.5.

%Regarding $J_{w}$, we have proved in Lemma 3.3 that $I_{w}
%\subseteq K$ for every maximal ideal $K$ of $O_{w}$. Therefore
%$I_{w} \subseteq J_{w}$, and again we have the natural epimorphism
%$O_{w}/I_{w} \twoheadrightarrow O_{w}/J_{w}$.

Q.E.D.\enddemo

\demo {Example 2.7} Let $F$ be a field with valuation $v$, and let
$E$ be a finite dimensional $F$-algebra with basis $B=\{ b_{1}=1,
b_{2},....,b_{n} \}$. Assume that $B$ satisfies the following: for
every $1\leq i,j \leq n$, $b_{i}b_{j}=\delta b_{k}$ for some $1
\leq k \leq n$ and $\delta \in O_{v}$ (in other words the basis is
projectively multiplicative with respect to $O_{v}$). Define
$$w(\sum_{i=1}^{n} \alpha_{i}b_{i})=\min_{1 \leq i \leq n} \{
v(\alpha_{i}) \}.
$$ It is
obvious that $w$ extends $v$; we show that $w$ is a
quasi-valuation. Indeed, assuming
$$v(\alpha_{i_{0}})=w(\sum_{i=1}^{n} \alpha_{i}b_{i}) \leq
w(\sum_{i=1}^{n} \beta_{i}b_{i})=v(\beta_{i_{1}}),$$ we have
                          $$w(\sum_{i=1}^{n} \alpha_{i}b_{i}+\sum_{i=1}^{n}
\beta_{i}b_{i})=w(\sum_{i=1}^{n} (\alpha_{i}+\beta_{i})b_{i})$$
 $$=\min_{1 \leq i \leq n} \{ v(\alpha_{i}+\beta_{i}) \} \geq \min_{1
\leq i \leq n} \{
 \min \{v(\alpha_{i}),v(\beta_{i})\} \}$$ $$=v(\alpha_{i_{0}})= w(\sum_{i=1}^{n} \alpha_{i}b_{i}).$$

 Now, $$w(\sum_{i=1}^{n} \alpha_{i}b_{i}\sum_{j=1}^{n}
\beta_{j}b_{j})= w(\sum_{i,j} \alpha_{i}\beta_{j}b_{i}b_{j})=
w(\sum_{i,j} \alpha_{i}\beta_{j}\delta_{i,j}b_{i,j}),$$ for
$\delta_{i,j} \in O_{v}$ and $b_{i,j} \in B$. By the previous
calculation, we have $$w(\sum_{i,j}
\alpha_{i}\beta_{j}\delta_{i,j}b_{i,j}) \geq \min_{i,j} \{
w(\alpha_{i}\beta_{j}\delta_{i,j}b_{i,j}) \}= \min_{i,j} \{
v(\alpha_{i}\beta_{j}\delta_{i,j}) \}$$ $$ \geq \min_{i,j} \{
v(\alpha_{i}\beta_{j}) \}= \min_{i,j} \{
v(\alpha_{i})+v(\beta_{j}) \}= \min_{i} \{ v(\alpha_{i})\}+
\min_{j} \{ v(\beta_{j})\}=$$ $$ w(\sum_{i=1}^{n}
\alpha_{i}b_{i})+w(\sum_{j=1}^{n} \beta_{j}b_{j}).$$

\enddemo

\proclaim {Lemma 2.8} Let $E/F$ be a finite field extension and
let $w$ be a quasi-valuation on $E$ extending a valuation $v$ on
$F$. Then $w(x) \neq \infty$ for all $0\neq x \in E$.\endproclaim

\demo {Proof} If $w(x)>0$, writing $\sum_{i=0}^{n}
\alpha_{i}x^{i}=0$ for $\alpha_{i}\in O_{v}$ and $\alpha_{0} \neq
0$, we have
$$w(x) \leq w(\alpha_{0})=v(\alpha_{0}) < \infty.$$

Q.E.D.\enddemo

The following example will be used to provide negative answers to
several natural questions.

\demo {Example 2.9} Let $v$ be a valuation on a field $F$ and let
$(M,+)$ be a totally ordered abelian monoid containing
$\Gamma_{v}=v(F \setminus \{0\})$. Let $E=F [e]$ be a Kummer field
extension with $e \notin F$ and $e^{2} \in O_{v}$. Define $w$ on
$E$ in the following way:
$$w (a + be)= \min \{ v(a), v(b)- \gamma \}$$ for some invertible
element $ \gamma \in M$, $ \gamma \geq 0$. We shall prove that $w$
is indeed a quasi-valuation extending $v$: assume $w (a + be) \leq
w (c + de)$; then

$$ w (a + be+c + de)=\min \{ v(a+c), v(b+d)-\gamma  \} $$
$$\geq  \min \{ \min \{v(a),v(c)\}, \min \{v(b),v(d)\}-\gamma  \}$$
$$= \min \{ \min \{v(a),v(c)\}, \min \{v(b)-\gamma,v(d)-\gamma\}  \}$$
$$= \min \{ v(a),v(b)-\gamma\}=w (a + be).$$

Next,

$$ w ((a + be)(c + de))=\min \{ v(ac+bde^{2}), v(ad+bc)-\gamma  \} $$
$$\geq  \min \{ \min \{v(ac),v(bde^{2})\}, \min \{v(ad),v(bc)\}-\gamma  \}$$
$$= \min \{ \min \{v(a)+v(c),v(b)+v(d)+v(e^{2})\}, \min \{v(a)+v(d)-\gamma,v(b)+v(c)-\gamma\}  \}$$
$$= \min \{ v(a)+v(c),v(b)+v(d)+v(e^{2}),v(a)+v(d)-\gamma,v(b)+v(c)-\gamma  \}$$
$$\geq \min \{ v(a)+v(c),v(b)+v(d),v(a)+v(d)-\gamma,v(b)+v(c)-\gamma  \}.$$

Note that  $$(2) \ \ w (a + be)+w (c + de)= \min \{ v(a),
v(b)-\gamma \}+ \min \{ v(c), v(d)-\gamma \}.$$ So, if the sum on
the right side of (2) is $v(a)+v(c)$, then $$v(a)+v(c) \leq \min
\{ v(a)+v(d)-\gamma, v(b)+v(c)-\gamma, v(b)-\gamma+v(d)-\gamma
\}$$ and thus $\leq v(b)+v(d)$. We note that if $\gamma > 0$ then
in this case we have an equality $w ((a + be)(c + de))=w (a +
be)+w (c + de)$. If the sum on the right side of (2) is
$v(a)+v(d)-\gamma$ or $v(b)-\gamma+v(c)$, we have a similar
situation and again an equality (if $\gamma >0$). If the sum on
the right side of (2) is $v(b)-\gamma+v(d)-\gamma$, then
$$v(b)-\gamma+v(d)-\gamma \leq \min \{ v(a)+v(c), v(b)+v(d), v(a)+v(d)-\gamma,
v(b)+v(c)-\gamma \}.$$ So, we have $w ((a + be)(c + de)) \geq w (a
+ be)+w (c + de).$

Note that if $\gamma \in \Gamma_{v}$ then taking $c \in F$ such
that $v(c)= \gamma$, we have $O_{w}=O_{v}[ce]$. If $\gamma \notin
\Gamma_{v}$ then $$O_{w}= \bigcup_{c \in F, \  v(c)>\gamma}
O_{v}[ce].$$

For example if one takes $F= \Bbb Q$, $v=v_{p}$, $M=\Bbb Z$ and
$\gamma=n \in \Bbb N \cup \{ 0\}$; then $O_{w}=\Bbb
Z_{p}[p^{n}e]$.
%(where $\Bbb Z_{p}$ is the localization of $\Bbb Z$ in $p$)).
\enddemo

\heading \S 3 INC, LO and K-dim \endheading

For any ring $R$, we denote by K-dim$R$ the classical Krull
dimension of $R$, by which we mean the maximal length of the
chains of prime ideals of $R$. Our next goal is to prove that
$O_{w}$ satisties INC (incomparability) and LO (lying over) over
$O_{v}$; cf. [Row, p. 184-192].
%K-dim$O_{w}=$K-dim$O_{v}$, by means of the standard comparisons of
%the prime spectrum: INC (incomparability), GU (going up), and LO
%(lying over); cf. [Row, p. 184-192].

For the reader's convenience we review these notions here. Let $C
\subseteq R$ be rings. We say that $R$ satisfies INC over $C$ if
whenever $Q_{0} \subset Q_{1}$ in $\text{Spec}(R)$ we have $Q_{0}
\cap C \subset Q_{1} \cap C$. We say that $R$ satisfies LO over
$C$ if for any $P \in \text{Spec}(C)$ there is $Q \in
\text{Spec}(R)$ lying over $P$, i.e., $Q \cap C=P$.

\demo {Remark 3.1} Let $v$ be a valuation with valuation ring
$O_{v}$ and let $I $ be an ideal of $O_{v}$. If $x \in I$ and $y
\in O_{v}$ such that $v(y) \geq v(x)$, then $y \in I$. Indeed,
\linebreak $O_{v}y \subseteq O_{v}x \subseteq ~I$.\enddemo

\demo {Remark 3.2} Let $C \subseteq R$ be rings. Let
$I_{1}\subseteq I_{2}$ be ideals of $R$ such that $I_{1}\cap
C=I_{2}\cap ~C$. If $x\in I_{2}$ and $\sum_{i=0}^m
\alpha_{i}x^{i}\in I_{2}$ for $\alpha_{i}\in C$, then
$\alpha_{0}\in I_{1}$. Indeed, $\sum_{i=1}^m \alpha_{i}x^{i}\in
I_{2}$ (since $x\in I_{2}$), so $\alpha_{0}\in I_{2}$. Therefore
$\alpha_{0}\in I_{1}$ (since $\alpha_{0}\in I_{2}\cap C=I_{1} \cap
C$).\enddemo

\proclaim {Lemma 3.3} Let $C \subseteq R$ be commutative rings.
Let $Q_{1}\subseteq Q_{2}$ be prime ideals of $R$ with $Q_{1}\cap
C=Q_{2}\cap C$. If $\sum_{i=0}^m \alpha_{i}x^{i}\in Q_{1}$ for
$\alpha_{i}\in C$ and $x \in Q_{2} \setminus Q_{1}$; then each
$\alpha_{i} \in Q_{1}$.\endproclaim

\demo {Proof} By Remark 3.2, $\alpha_{0} \in Q_{1}$. Note that $x
\sum_{i=1}^m \alpha_{i}x^{i-1} \in Q_{1}$ and since $x \notin
Q_{1}$ we have $\sum_{i=1}^m \alpha_{i}x^{i-1} \in Q_{1}$. Again,
we apply Remark 3.2 and get $\alpha_{1} \in Q_{1}$. We continue by
induction and get each $\alpha_{i} \in Q_{1}$.

Q.E.D.\enddemo

\demo{Definition 3.4} Let $C \subseteq R$ be rings. We call $R$
{\it quasi-integral} over $C$ if for every $r \in R$ there exists
a polynomial $\sum c_{i} r^{i}=0$ with $c_{i} \in C$ and
$c_{i_{0}}=1$ for some $i_{0}$ (depending on $r$).\enddemo

\demo{Remark 3.5} Recall that $E$ is a finite field extension of a
field $F$ with valuation $v$ and valuation ring $O_{v}$. Note that
for any $E \supseteq R \supseteq O_{v}$, $R$ is quasi-integral
over $O_{v}$, by Remark 2.2 applied to $\{ 1,r,...,r^{n} \}$. In
particular, the quasi-valuation ring $O_{w}$ is quasi-integral
over $O_{v}$.\enddemo

\demo {Remark 3.6} If $R$ is quasi-integral over $C$ and $Q
\vartriangleleft R$, then $R/Q$ is quasi-integral over $C /C \cap
Q$.\enddemo

\demo {Proof} For any $r+Q \in R/Q$, there exists a polynomial
$\sum_{i=0}^n c_{i}r^{i}=0$ for $c_{i} \in C$ and $c_{i_{0}}=1$
for some $i_{0}$. We pass to the image in $R/Q$ over $(C+Q)/Q
\approx C/C \cap Q$.

Q.E.D.\enddemo

\proclaim{Theorem 3.7} Any commutative quasi-integral ring
extension $R$ of $O_{v}$ satisfies INC over $O_{v}$.\endproclaim

\demo{Proof} Let $Q_{1} \subset Q_{2}$ be prime ideals of $R$.
Assume to the contrary, that \linebreak $Q_{1} \cap O_{v}=Q_{2}
\cap O_{v}$. Let $x \in Q_{2} \setminus Q_{1}$ and by assumption
write $\sum_{i=0}^n \alpha_{i}x^{i}=0$ for $\alpha_{i} \in O_{v}$,
$\alpha_{i_{0}}=1$ for some $i_{0}$ . By Lemma 3.3,
$\alpha_{i_{0}} \in Q_{1}$, a contradiction.

Q.E.D.\enddemo

\proclaim{Corollary 3.8} If $R$ is commutative and quasi-integral
over $O_{v}$, then K-dim$R \leq$K-dim$O_{v}$.\endproclaim

\demo {Proof} By Theorem 3.7 $R$ satisfies INC over $O_{v}$; thus
K-dim$R \leq$K-dim$O_{v}$.

Q.E.D.\enddemo

\proclaim{Corollary 3.9} K-dim$O_{w} \leq$K-dim$O_{v}$, for any
quasi-valuation $w$ extending $v$ on a finite dimensional field
extension $E$ of $F$.\endproclaim

\demo {Proof} $O_{w}$ is quasi-integral over $O_{v}$ by Remark
3.5; now, use Corollary 3.8.

Q.E.D.\enddemo

%Note that if $Q \in Spec (O_{w})$, then $$O_{v}/(Q \cap
%O_{v})\approx (O_{v}+Q)/Q$$ is a subring of $O_{w}/Q$. We write
%$P=Q \cap O_{v}$. Although one can easily see that $O_{v}/P$ is a
%valuation ring, we write its corresponding valuation explicitly:
%$$\overline{v}(x+P)=v(x) \text{\ if\ } x \notin P \text{\ and\ }
%\overline{v}(x+P)=\infty \text{\ otherwise}.$$ Note that
%$\overline{v}$ is well defined; indeed if $x+P=x'+P$ then $x-x'
%\in P$. Assume on the contrary that $v(x) \neq v(x')$, say $v(x)<
%v(x')$; then
%$$v(x')>v(x)=v(x-x').$$ Therefore, by Remark 3.1, $x' \in P$ and
%we have $$\overline{v}(x+P)=\overline{v}(x'+P)=\infty.$$ Note that
%$\Gamma_{\overline{v}} \subseteq \Gamma_{v}$, and
%$$\overline{v}(x+P) \leq \overline{v}(y+P) \text{ \ iff\ } v(x) \leq v(y).$$

The following three lemmas are valid in the more general case
where $E$ is a finite dimensional $F$-algebra.

%\proclaim {Remark 2.10} Let $A$ be an integral domain with field
%of fractions $F$, $E$ a finite dimensional $F$-algebra and $B$ a
%subring of $E$ with $B \cap F=A$. Let $I$ be a proper ideal of
%$A$; then every $x \in IB$ can be written as $x=ar$ for $a \in I$,
%$r \in B$.
%\endproclaim

\proclaim {Lemma 3.10} Let $I$ be a proper ideal of $O_{v}$; then
every $x \in IO_{w}$ can be written as $x=ar$ for $a \in I$, $r
\in O_{w}$.
\endproclaim

\demo {Proof} Every $x \in IO_{w}$ is of the form $\sum_{i=1}^{n}
a_{i}r_{i}$ for $a_{i}\in I$, $r_{i}\in O_{w}$; so take
$a=a_{i_{0}}$ with minimal $v$-value and write $x=ar$ for
appropriate $r \in O_{w}$.

Q.E.D.\enddemo

\proclaim {Lemma 3.11} Let $I$ be a proper ideal of $O_{v}$; then
$IO_{w}$ is a proper ideal of $O_{w}$ and $IO_{w} \cap O_{v}=I$.
\endproclaim

\demo {Proof} Let $x \in IO_{w}$; by Lemma 3.10 we may write
$x=ar$ for $a \in I$ and $r \in O_{w}$. Thus, since $w(ar) \geq
w(a)=v(a)>0$, $1 \notin IO_{w}$ and $IO_{w}$ is a proper ideal of
$O_{w}$. Now, let $x  \in IO_{w} \cap O_{v}$ and write $x=ar$ for
$a \in I$ and $r \in O_{w}$. Then $r \in F \cap O_{w}$ and thus $r
\in O_{v}$; hence $IO_{w} \cap O_{v}=I$.

Q.E.D.\enddemo

 \proclaim {Lemma 3.12} $O_{w}$ satisfies LO over $O_{v}$.
\endproclaim

\demo {Proof} Let $P \in \text{Spec} (O_{v})$ and denote $T=\{ A
\vartriangleleft O_{w} \mid A \cap O_{v}=P\}$. By Lemma 3.11,
$PO_{w}$ is a proper ideal of $O_{w}$ lying over $P$. Thus $T \neq
\emptyset $. Now, $T$ with the partial order of containment
satisfies the conditions of Zorn's Lemma and thus there exists $Q
\vartriangleleft O_{w}$ lying over $P$, maximal with respect to
containment. We shall prove that $Q \in \text{Spec} (O_{w})$.
Indeed, let $A,B$ be ideals in $O_{w}$ with $AB \subseteq Q$ and
assume that $A,B \nsubseteq Q$. Then
$$(Q+A) \cap O_{v} \neq P \text { \ and \ } (Q+B) \cap O_{v}
\neq P;$$ so there exist $a \in A, \ b \in B, \ q,q' \in Q$ such
that $q+a,q'+b \in O_{v} \setminus P$. However
$$(q+a)(q'+b)=qq'+qb+aq'+ab \in Q \cap O_{v}=P,$$ a contradiction.

Q.E.D\enddemo

 \heading \S 4 PIMs and GD\endheading

%( Recall that an isolated subgroup of a totally ordered abelian
%group $\Gamma$ is a subgroup $H$ of $\Gamma$ such that $\{ \gamma
%\in \Gamma | 0 \leq \gamma \leq h \} \subseteq H$ for any $h \in
%H$. Note that the quotient group $\overline{\Gamma} = \Gamma / H$
%is totally ordered, by setting $\overline{\gamma} \geq
%\overline{0}$ iff $\gamma \geq h$ for some $h \in H$. There is a
%one to one correspondence between the set of all prime ideals of a
%valuation ring and $G(\Gamma)$ the set of all isolated subgroups
%of $\Gamma$. The rank of $\Gamma$ is the order type of $G(\Gamma)
%\setminus \{ \Gamma \}$ (see [End, chap. 13]). Also recall that if
%$\Gamma$ and $\Delta$ are totally ordered abelian groups such that
%$\Gamma \subseteq \Delta$  and $ \Delta /\Gamma$ is a torsion
%group, then $\text{rank} \Gamma = \text{rank} \Delta$. See [End,
%chap. 13] or [Bo, section 4] for further discussion.

In this section we introduce the notion of PIM. We show that the
convex hull of $H^{\geq0}$ in $M$, where $H$ is an isolated
subgroup of $\Gamma$, is a PIM. We use the closures of these PIMs
to prove that $O_{w}$ satisfies GD over $O_{v}$. We also obtain a
connection between the height of a prime ideal $P$ of $O_{v}$ and
the height of a prime ideal of $O_{w}$ lying over $P$. Finally, we
present a connection between the closure of a PIM and the prime
ideals of $O_{w}$.

 From now through Corollary 4.10 we consider abstract
properties of a totally ordered abelian monoid $M$ containing a
group ~$\Gamma$.

Let $N \subseteq M$, we denote
$$N^{\geq 0}= \{ m \in N \mid m \geq 0 \} $$ and $$N^{<0}= \{
\alpha \in N \mid \alpha < 0 \}$$ In particular, $M^{\geq 0}= \{ m
\in M \mid m \geq 0 \}. $ Note that $M^{\geq 0}$ is a submonoid of
$M$.

\demo{Definition 4.1} A subset $N \subseteq M^{\geq 0}$ is called
a {\it positive isolated monoid} (PIM) if $N$ is a submonoid of
$M^{\geq 0}$ and for every $\gamma \in N$ we have $$\{ \delta \in
M^{\geq 0} \mid \delta \leq \gamma \} \subseteq N.$$\enddemo

Recall that an isolated subgroup of a totally ordered abelian
group $\Gamma$ is a subgroup $H$ of $\Gamma$ such that $\{ \gamma
\in \Gamma | 0 \leq \gamma \leq h \} \subseteq H$ for any $h \in
H$ (some texts call such subgroups "convex" or "distinguished");
see [End, p. 47].
%Note that the quotient group $\overline{\Gamma} = \Gamma / H$ is
%totally ordered, by setting $\overline{\gamma} \geq \overline{0}$
%iff $\gamma \geq h$ for some $h \in H$. There is a one to one
%correspondence between the set of all prime ideals of a valuation
%ring and the set $G(\Gamma)$ of all isolated subgroups of $\Gamma$
%given by $P \mapsto \{ \alpha \in \Gamma \mid \alpha \neq v(p)$
%and $\alpha \neq -v(p)$ for all $p \in P\}$. The rank of $\Gamma$
%is the order type of $G(\Gamma) \setminus \{ \Gamma \}$ ([End, p.
%47]).

For any totally ordered set $A$ and a subset $B \subseteq A$ the
convex hull of $B$ in $A$ is defined by
$$hull_{A}(B)=\{a \in A \mid \exists b_{1},b_{2} \in B \text{ such
that } b_{1} \leq a \leq b_{2} \}.$$
%In particular, for totally
%ordered abelian groups $\Gamma \subseteq \Delta$ and an isolated
%subgroup $H$ of $\Gamma$, the convex hull of $H$ in $\Delta$ is an
%isolated subgroup of $\Delta$.

Let $H$ be an isolated subgroup of $\Gamma$; it is not difficult
to see that $$hull_{M}(H^{\geq0})=[hull_{M}(H)]^{\geq0}.$$

%\demo{Definition 3.2} Let $T$ be any subset of $\Gamma$. We denote
%%$$T^{\geq0}= \{ \alpha \in T \mid \alpha \geq 0 \},$$
%$$T^{<0}= \{ \alpha \in T \mid \alpha < 0 \}$$ and

%$$c_{M}({T}) = \{ m \in M \mid -n\alpha \leq m \leq n\alpha \text{\ for
%some\ }  \alpha \in T,\ \alpha \geq 0\}.$$

%*****************************************************************************

%$$c_{M}({H}) = \{ m \in M \mid 0 \leq m \leq \alpha \text{\ for
%some\ } \alpha \in H^{\geq0}) \}. $$ $c_{M}(H)$ is called the convex
%hull of $H$ inside $M$; when there is no danger of confusion we
%shall write $c(H)$ instead of $c_{M}(H)$\enddemo

%It will be understood from the context whether we consider a set
%of elements grater or equal to 0 in $M$ or in $\Gamma$; thus, we
%shall omit the subscripts $M$ and $\Gamma$.

%\demo{Remark 3.3} $c({H^{\geq0}}) = \{ \sum m_{i} \mid 0 \leq
%m_{i} \leq \alpha_{i}$ for  $m_{i}\in M^{\geq 0}$, $\alpha_{i} \in
%H^{\geq0}) \}$.
%\enddemo
%Indeed If $0 \leq m_{i} \leq
%\alpha_{i}$ for  $m_{i}\in M^{\geq 0}$ and $\alpha_{i} \in H^{\geq0})$, then
%$\sum m_{i} \leq \sum \alpha_{i} \in H^{\geq0})$ since $H^{\geq0})$ is
%closed under addition. The other direction is trivial.\enddemo

%Assume on the contrary that $\gamma_{i_{0}} > \alpha$ for every
%$\alpha \in H^{\geq0})$; then $\sum \gamma_{i} \geq \gamma_{i_{0}} >
%\alpha$ for every $\alpha \in H^{\geq0})$, a contradiction.
%\enddemo

\demo{Remark 4.2} Let $H \leq \Gamma$ be an isolated subgroup.
$hull_{M}(H^{\geq0})$ is a PIM and $hull_{M}(H^{\geq0}) \cap
\Gamma=H^{\geq0}$. We prove that $hull_{M}(H^{\geq0})$ is indeed a
PIM. Let $m \in hull_{M}(H^{\geq0})$ and $0 \leq m' \leq m $. Then
there exists $\alpha \in H^{\geq0}$ such that $0 \leq m' \leq m
\leq \alpha$; thus $m' \in hull_{M}(H^{\geq0})$.
$hull_{M}(H^{\geq0})$ is a monoid since $H^{\geq0}$ is closed
under addition.
\enddemo

\demo{Remark 4.3} The PIMs of $M^{\geq 0}$ are linearly ordered by
inclusion.\enddemo

\demo{Proof} Assume $N_{1} \neq N_{2}$ and assume there exists
$\gamma \in N_{1} \setminus N_{2}$. Then for any $\delta \in N_{2}
$ we have $\gamma > \delta$ (since otherwise $\gamma \in N_{2}$),
implying $\delta \in N_{1}$ by assumption that $N_{1}$ is
isolated. Thus $N_{2} \subseteq N_{1}$.\enddemo

\demo{Remark 4.4} The correspondence $H^{\geq0} \rightarrow
hull_{M}(H^{\geq0})$ from $$\{ \text{positive part of isolated
subgroups of\ } \Gamma \} \rightarrow \{\text{PIMs} \text{\ of\ }
M^{\geq 0}\}$$ is an injection.\enddemo

\demo{Proof} Recall that the isolated subgroups are linearly
ordered, and thus their positive parts are also totally ordered.
Assume that $(H_{1})^{\geq0} \supset (H_{2})^{\geq0} $ and let
\linebreak $\alpha \in (H_{1})^{\geq0} \setminus (H_{2})^{\geq0}$.
Then $$\alpha \in hull_{M}((H_{1})^{\geq0}) \setminus
hull_{M}((H_{2})^{\geq0}).$$

Q.E.D.\enddemo

\demo{Definition 4.5} A PIM $N \subseteq M^{\geq 0}$ is said to be
{\it lie over} the positive part $H^{\geq0}$ of an isolated
subgroup $H$ if $N \cap \Gamma = H^{\geq0}$.
\enddemo

\demo{Example 4.6} Let $N$ be a PIM of $M^{\geq 0}$. Then $$H=(N
\cap \Gamma) \cup \{ -\alpha \mid \alpha \in  N \cap \Gamma\}$$ is
an isolated subgroup of $\Gamma$, and $N$ lies over $H^{\geq0}$.
\enddemo

\demo{Remark 4.7} $hull_{M}(H^{\geq0})$ is the minimal PIM lying
over $H^{\geq0}$; i.e., for any PIM $N \subset
hull_{M}(H^{\geq0})$ we have $N \cap \Gamma \subset H^{\geq0} $.
\enddemo

\demo{Proof} Take $\gamma \in hull_{M}(H^{\geq0}) \setminus N$;
then $0 \leq \gamma \leq \alpha$ for some $\alpha \in H^{\geq0}$.
Assuming $N$ lies over $H^{\geq0}$ we get $\alpha \in N$ and thus
$\gamma \in N$, a contradiction.

Q.E.D.\enddemo

\demo{Definition 4.8} Let $H \leq \Gamma$ be an isolated subgroup
and let $hull_{M}(H^{\geq0})$ be its corresponding PIM. The {\it
closure} of $hull_{M}(H^{\geq0})$, denoted
$\overline{hull_{M}(H^{\geq0})}$, is the set
$$\bigcup \{ N_{i} \mid N_{i} \text{\ is a PIM lying over\ } H^{\geq0}
\}.$$
\enddemo

\proclaim{Lemma 4.9} Notation as in Definition 4.8.
$\overline{hull_{M}(H^{\geq0})}$ is a PIM of $M^{\geq 0}$ lying
over $H^{\geq0}$.\endproclaim

\demo{Proof} It is clear that $\overline{hull_{M}(H^{\geq0})}$ is
a PIM and that $H^{\geq0} \subseteq
\overline{hull_{M}(H^{\geq0})}$. On the other hand, if $\alpha \in
\overline{hull_{M}(H^{\geq0})} \cap \Gamma$ then $\alpha \in
N_{i}\cap \Gamma$ for some $N_{i}$ lying over $H^{\geq0}$. Thus
$\alpha \in H^{\geq0}$.

Q.E.D.\enddemo

\proclaim{Corollary 4.10} $\overline{hull_{M}(H^{\geq0})}$ is the
maximal PIM lying over $H^{\geq0}$.\endproclaim

We note at this point that for an isolated subgroup $H$ of
$\Gamma$ and $m \in M^{\geq 0}$, $m \notin hull_{M}(H^{\geq0})$
iff $m> \alpha$ for every $\alpha \in H^{\geq0}$. However, one can
have $m> \alpha$ for every $\alpha \in H^{\geq0}$ and still $m \in
\overline{hull_{M}(H^{\geq0})}$. We shall see in Lemma 10.6 an
example of such an element (we denote it by $\Cal C$), where the
monoid is the cut monoid. Examples 10.7 and 10.8 also provide
examples of such elements. We shall discuss more about this fact
after these examples.

\heading Going Down \endheading

We shall now use the notion of PIM to prove Going Down (GD) for
quasi-valuation rings over valuation rings. Recall (see [Row, p.
191]) that in case $C \subseteq R$ are rings, we say that $R$
satisfies GD (going down) over $C$ if for any $P_{0} \subseteq
P_{1} \in \text{Spec} (C)$ and every $Q_{1} \in \text{Spec} (R)$
lying over $P_{1}$ there is $Q_{0} \subseteq Q_{1}$ in
$\text{Spec} (R)$ lying over $P_{0}$.

As usual, $v $ is a valuation on a field ~$F$, $O_{v}$ its
valuation ring, $E$ is a finite field extension of $F$ and $w$ is
a quasi-valuation on $E$ extending $v$ with quasi-valuation ring
$O_{w}$ and value monoid $M_{w}$.

Recall that there is a one to one correspondence between the set
of all prime ideals of a valuation ring and the set
$G(\Gamma_{v})$ of all isolated subgroups of $\Gamma_{v}$ given by
$P \mapsto \{ \alpha \in \Gamma_{v} \mid \alpha \neq v(p)$ and
$\alpha \neq -v(p)$ for all $p \in P\}$. (See [End. p. 47]).

\proclaim{Lemma 4.11} Let $P \in \text{Spec}(O_{v})$, $Q \in
\text{Spec}(O_{w})$ such that $Q \cap O_{v} = P$ and $x \in O_{w}
\setminus Q$. Let $H$ denote the isolated subgroup corresponding
to $P$ and let $hull_{M_{w}}(H^{\geq0})$ denote the corresponding
PIM. Then $$w(x) \in \overline{hull_{M_{w}}(H^{\geq0})}.$$
\endproclaim

\demo{Proof} Assume to the contrary that $w(x) \notin
\overline{hull_{M_{w}}(H^{\geq0})}$. Let $$N= \{ m \in
(M_{w})^{\geq 0} \mid m \leq nw(x) \text{\ for some\ }  n \in \Bbb
N \}.$$ $N$ is clearly a PIM strictly containing
$\overline{hull_{M_{w}}(H^{\geq0})}$; thus $N \cap \Gamma_{v}
\supset H^{\geq0}$. Let $\alpha \in (N \cap \Gamma_{v}) \setminus
H^{\geq0}$ and let $a \in F$ with $v(a)= \alpha$. Note that $a \in
P$ since $v(a) \notin H^{\geq0}$. By the definition of $N$, there
exists
 $n \in \Bbb N$ such that
$$w(x^{n}) \geq nw(x) \geq v(a).$$ Therefore $$x^{n} \in aO_{w}
\subseteq PO_{w} \subseteq Q$$ and thus $x \in Q$, a
contradiction.

Q.E.D.\enddemo

For a monoid $M$ and subsets $N_{1},N_{2} \subseteq M$, we denote
$$N_{1}+N_{2}=\{ m_{1}+m_{2} \mid m_{1} \in N_{1}, m_{2} \in
N_{2}\}.$$

 \proclaim{Lemma 4.12} $O_{w}$ satisfies GD over
$O_{v}$.\endproclaim

\demo{Proof} Let $P_{1} \subset P_{2} \in \text{Spec}(O_{v})$ and
let $Q_{2} $ be a prime ideal of $O_{w}$ such that $Q_{2} \cap
O_{v} = P_{2}$. We shall prove that there exists a $Q_{1} \in
\text{Spec}(O_{w})$ such that $Q_{1} \cap O_{v} = P_{1}$ and
$Q_{1} \subset Q_{2}$. We denote $S_{1} = O_{v} \setminus P_{1}$,
$S_{2} = O_{w} \setminus Q_{2}$ and $$S = \{ s_{1} s_{2} | s_{1}
\in S_{1}, s_{2} \in S_{2} \}.$$ Note that $S$ is a multiplicative
monoid. We shall prove $S \cap P_{1}O_{w} = \emptyset$ and then
every ideal $Q_{1}$ which contains $P_{1}O_{w}$, maximal with
respect to $S \cap Q_{1} = \emptyset$, is prime. Note that such
$Q_{1}$ satisfies the required properties, $Q_{1} \cap O_{v} =
P_{1}$ and $Q_{1} \subset Q_{2}$.

Let $H_{i} \leq \Gamma_{v}$ $(i=1,2)$ be the isolated subgroups
corresponding to $P_{i}$ and let $hull_{M_{w}}((H_{i})^{\geq0})$
be their corresponding PIMs in $(M_{w})^{\geq 0}$.

Let $x \in P_{1}O_{w}$ and write, by Lemma 3.10, $x = pa$ where $p
\in P_{1}$, $a \in O_{w}$. Then
$$w(x) = w(pa)  \geq
w(p) = v(p)  \notin H_{1}.$$ %for an appropriate $\alpha \in
%\Gamma_{v} \setminus H_{1}$.
%Let $x \in P_{1}O_{w}$, and write $x = \sum_{i=1}^{n} p_{i}a_{i}$
%where $p_{i} \in P_{1}$, $a_{i} \in O_{w}$. Then $$w(x) \geq
%\min_{1 \leq i \leq n} \{ w(p_{i}a_{i}) \} \geq w(p_{i_{0}}) =
%v(p_{i_{0}}) = \alpha \notin H_{1}$$ for an appropriate $p_{i_{0}}
%\in P_{1}$, $\alpha \in \Gamma_{v} \setminus H_{1}$.
We claim $w(x) \notin hull_{M_{w}}((H_{1})^{\geq0})$. Indeed,
otherwise, $w(x) \leq \beta$ for some $\beta \in (H_{1})^{\geq0}$,
a contradiction.

Now, let $y \in S$, write $y=ab$ for $a \in S_{1}$, $b \in S_{2}$.
By Lemma 4.11, $$w(b) \in
\overline{hull_{M_{w}}((H_{2})^{\geq0})}.$$ As for $a$, $a \in
S_{1}$ and thus $w(a)=v(a) \in H_{1}$. In particular $w(a) \in
hull_{M_{w}}((H_{1})^{\geq0})$. Note that $P_{1} \subset P_{2}$
and thus $H_{2} \subset H_{1}$ and $hull_{M_{w}}((H_{2})^{\geq0})
\subset hull_{M_{w}}((H_{1})^{\geq0})$. Thus,
$$\overline{hull_{M_{w}}((H_{2})^{\geq0})} \subset
hull_{M_{w}}((H_{1})^{\geq0}).$$ (Indeed, since $H_{2} \subset
H_{1}$, $hull_{M_{w}}((H_{1})^{\geq0})$ contains an $\alpha \in
(H_{1})^{\geq0} \setminus (H_{2})^{\geq0}$ and any PIM $N_{2}$
lying over $(H_{2})^{\geq0}$ does not contain $\alpha$). Also note
that $a$ is stable with respect to $w$. So we have
$$w(ab) = w(a) + w(b) \in hull_{M_{w}}((H_{1})^{\geq0}) + \overline{hull_{M_{w}}((H_{2})^{\geq0})} =
hull_{M_{w}}((H_{1})^{\geq0}).$$ %where

%$$hull_{M_{w}}((H_{1})^{\geq0}) + \overline{hull_{M_{w}}((H_{2})^{\geq0})}=\{m_{1}+m_{2} \mid
%m_{1} \in hull_{M_{w}}((H_{1})^{\geq0}), m_{2} \in
%\overline{hull_{M_{w}}((H_{2})^{\geq0})}\}$$

So, we have proved that $$ w(S) \subseteq
hull_{M_{w}}((H_{1})^{\geq0}) $$ and $$w(P_{1}O_{w}) \subseteq
(M_{w})^{\geq 0} \setminus hull_{M_{w}}((H_{1})^{\geq0}).$$
Therefore $ S \cap P_{1}O_{w} = \emptyset.$

Q.E.D

\enddemo

Note: it is possible to prove GD in the more general case in which
$E$ is a finite dimensional $F$-algebra, but the proof is more
complicated and uses the notion of filter quasi-valuations (to be
intruduced in section 9).

\proclaim{Corollary 4.13} K-dim$O_{w} =$K-dim$O_{v}$.
\endproclaim

 \demo {Proof} By Remark 3.5 and Theorem 3.7, $O_{w}$
satisfies INC over $O_{v}$; thus, we get K-dim$O_{w}  \leq
$K-dim$O_{v}$. By Lemma 3.12, $O_{w}$ satisfies LO over $O_{v}$
and by Lemma 4.12 , $O_{w}$ satisfies GD over $O_{v}$; thus
K-dim$O_{w} \geq $K-dim$O_{v}$.

Q.E.D.\enddemo

If $R$ is a ring and $P \in \text{Spec} (R)$, we write $h_{R}(P)$
for the height of $P$ in $R$. Note that if $R$ is a valuation ring
and K-dim$R<\infty$ then two prime ideals have the same height iff
they are equal, since the ideals are linearly ordered by
inclusion.

To avoid complicated notation, we shall write $h_{v}(P)$ instead
of $h_{O_{v}}(P)$ for $ P \in \text{Spec} (O_{v})$.% and $h_{w}(Q)$
%instead of $h_{O_{w}}(Q)$ for $ Q \in \text{Spec} (O_{w})$.

%Recall that we denote $h_{v}(P)$ for the height of $P \in
%\text{Spec}(O_{v})$ and $h_{w}(Q)$ for the height of $Q \in
%\text{Spec}(O_{w})$.

%If $R$ is a ring and $P \in Spec(R)$, we write $h_{R}(P)$ for the
%height of $P$ in $R$. Note that if $R$ is a valuation ring then
%two prime ideals have the same height iff they are equal. (since
%the ideals are linearly ordered by inclusion). To avoid
%complicated notation, we shall write $h_{v}(P)$ instead of
%$h_{O_{v}}(P)$ for $ P \in Spec(O_{v})$ and $h_{w}(Q)$ instead of
%$h_{O_{w}}(Q)$ for $ Q \in Spec(O_{w})$.

\proclaim{Lemma 4.14} Let $(F,v)$ be a valued field with valuation
ring $O_{v}$ and K-dim$O_{v}<\infty$. Let $E/F$ be a finite field
extension and let $R$ be a subring of $E$ lying over $O_{v}$.
Assume that $R$ satisfies GD over $O_{v}$. Let $Q \in \text{Spec}
(R)$ and $P \in \text{Spec} (O_{v})$; then $$P=Q \cap O_{v}
\text{\ iff\ } h_{R}(Q)=h_{v}(P).$$
\endproclaim

\demo{Proof} $(\Rightarrow)$ By Remark 3.5 and Theorem 3.7, $R$
satisfies INC over $O_{v}$; thus, we get $h_{R}(Q) \leq h_{v}(P)$.
By assumption, $R$ satisfies GD over $O_{v}$; thus $h_{v}(P) \leq
h_{R}(Q)$.

$(\Leftarrow)$ Write $Q \cap O_{v}=P'$. Then, by $(\Rightarrow)$,
$h_{R}(Q)=h_{v}(P')$ and thus $h_{v}(P) = h_{v}(P')$. Therefore
$P'=P$.

Q.E.D.\enddemo

\proclaim{Lemma 4.15} Let $Q \in \text{Spec}(O_{w})$, $P \in
\text{Spec}(O_{v})$ and assume K-dim$O_{v}< \infty$; then $$P=Q
\cap O_{v} \text{\ iff\ } h_{w}(Q)=h_{v}(P).$$
\endproclaim

\demo{Proof} $O_{w} \subseteq E$ lies over $O_{v}$, and by Lemma
4.12 $O_{w}$ satisfies GD over $O_{v}$. Now use Lemma 4.14.

Q.E.D.\enddemo

%By Remark 3.5 and Theorem 3.7, $O_{w}$ satisfies INC

Let $R$ be a ring; $R$ is said to satisfy the height formula if
for any ideal $P \in \text{Spec} (R)$, we have
$$ \text{K-dim}R=h_{R}(P)+\text{K-dim}(R/P).$$

The following lemma is a well known result, seen by matching
chains of prime ideals:

\proclaim{Lemma 4.16} If $C \subseteq R$ are commutative rings
such that $C$ satisfies the height formula and $R$ satisfies GU,
GD and INC over $C$, then $R$ satisfies the height formula.

\endproclaim

In Theorem 5.16 we shall prove that if $w(E \setminus \{0\})$ is
torsion over $\Gamma_{v}$ (and thus $M_{w}$ is actually a group)
then $O_{w}$ satisfies GU over $O_{v}$.

Since $O_{w}$ satisfies INC and GD over $O_{v}$, we deduce by
Lemma 4.16 the following theorem,

\proclaim{Theorem 4.17} If $O_{w}$ satisfies GU over $O_{v}$ then
$O_{w}$ satisfies the height formula.
\endproclaim

\proclaim{Lemma 4.18} Let $(F,v)$ be a valued field with valuation
ring $O_{v}$. Let $E$ be a finite field extension and let $R$ be a
subring of $E$ lying over $O_{v}$. Assume that $R$ satisfies GD
over $O_{v}$. Let $I $ be an ideal of $R$ containing $P \in
\text{Spec} (O_{v})$ and assume that $I \subseteq Q$ for all $Q
\in \Cal Q$ where $ \Cal Q = \{ Q \in \text{Spec} (R)| Q \cap
O_{v} = P \}$. Then $\sqrt{I} = \bigcap_{Q \in \Cal Q} Q$.

\endproclaim

\demo{Proof} $\sqrt{I}$ is a radical ideal and thus $\sqrt{I} =
\bigcap_{j \in J} S_{j}$ where the $S_{j}$ are the ideals in
$\text{Spec} (R)$ containing $I$. Now, we must have $S_{j} \cap
O_{v} \supseteq P$ for every such $S_{j}$ since $I \supseteq P$.
However, assuming some $S_{j_{0}}$ satisfies $S_{j_{0}} \cap O_{v}
\supset P$; we have, by GD, $Q_{i_{0}} \subset S_{j_{0}}$ for some
$Q_{i_{0}} \in \Cal Q$. Thus, since $I \subseteq Q$ for all $Q \in
\Cal Q$, we can take each $S_{j}$ to be in $\Cal Q$; finally,
since $\sqrt{I} \subseteq \bigcap_{Q \in \Cal Q} Q$, we have
$\sqrt{I} = \bigcap_{Q \in \Cal Q} Q.$

Q.E.D.\enddemo

\proclaim{Proposition 4.19} Let $P \in \text{Spec} (O_{v})$, $H$
its corresponding isolated subgroup in ~$\Gamma_{v}$ and
$hull_{M_{w}}(H^{\geq0})$ the corresponding PIM in $(M_{w})^{\geq
0}$. Let
$$I = \{ x \in O_{w}| w(x) \notin \overline{hull_{M_{w}}(H^{\geq0})}
\}$$ and let $ \Cal Q = \{ Q \in \text{Spec}(O_{w})| Q \cap O_{v}
= P \}$. Then
$$\sqrt{I} = \bigcap_{Q \in \Cal Q} Q.$$
\endproclaim

\demo{Proof} %First note that $\sqrt{I} \cap O_{v} = P$. Indeed, $P
%\subseteq I$ and thus $P \subseteq \sqrt{I}$. Conversely, if
%$w(x^{n}) \notin \overline{hull_{M_{w}}(H^{\geq0})}$ and $x \in O_{v}$, then
%$x^{n} \in P$ and therefore $x \in P$.
First note that $P \subseteq I$ since $a \in P$ implies $v(a)
\notin H$ and thus $w(a)=v(a) \notin
\overline{hull_{M_{w}}(H^{\geq0})} $. Let $x \in I$; then $$w(x)
\notin \overline{hull_{M_{w}}(H^{\geq0})}$$ and
thus by Lemma 4.11, $x \in Q$ for every $Q \in \Cal Q$. %Therefore
%$$ \sqrt{I} \subseteq \bigcap_{Q \in \Cal Q} Q.$$ Now,
%by [Sa, Remark 3.5 and Theorem 3.7] $O_{w}$ satisfies INC over
%$O_{v}$, and
 Now, by Lemma 4.12, $O_{w}$ satisfies GD over $O_{v}$.
Finally, use Lemma 4.18.

%*******For the rest of the proof, we need GD and Lemma 3.16.
%Thus, the rest of the proof is the same as in Proposition
%2.4.4.*********

Q.E.D.\enddemo

We shall use Proposition 4.19 to prove a similar result in section
6 (Proposition 6.5) where we get that the value monoid is a group.

%We note that [Sa, Proposition 4.7] can be seen to be a consequence
%of Proposition 4.19 because under the hypothesis that $M_{w}_{w}$ is
%torsion over the group $\Gamma_{v}$, the positive part of
%$\sqrt{H}$ is the unique PIM lying over $H^{\geq0}$ and thus
%equals $\overline{hull_{M}(H^{\geq0})}$.

\heading \S 5 The Prime Spectrum and GU \endheading

In this section we study expansions of $(O_{w},K)$ where $K$ is a
maximal ideal of $O_{w}$. We construct a quasi-valuation on a
localization of $O_{w}$, which enables us to prove that $O_{w}$
satisfies GU over $O_{v}$. We also obtain a bound on the size of
the prime spectrum of $O_{w}$.

 We start with the following easy lemma:

\proclaim{Lemma 5.1} Let $(G,+)$ be an abelian group and let $M
\supseteq G$ be a monoid such that $M$ is torsion over $G$; i.e.,
for every $m \in M$ there exists an $n \in \Bbb N$ such that $n m
\in G$. Then $M$ is a group.
\endproclaim

\demo{Proof} Let $m \in M$ and take $n \in \Bbb N$ such that $nm=g
\in G$. Thus \linebreak $-g \in G \subseteq M$ and $-m=-g+(n-1)m
\in M$.

Q.E.D.\enddemo

\demo {Example 5.2} Take $F= \Bbb Q$, $v=v_{p}$, $M=\Bbb R$ and
$\gamma=\pi$ in Example 2.9; then $w(E\setminus \{ 0 \})=\Bbb Z
\cup \{ z-\pi \mid z\in \Bbb Z \}$ is not torsion over
$\Gamma_{v}=\Bbb Z$.\enddemo

Note that Example 5.2 indicates that even in the case of a finite
field extension, $M_{w} $ (the monoid generated by $w(E\setminus
\{ 0 \})$) is not necessarily torsion over $\Gamma_{v}$. Also note
that (since $w$ in Example 5.2 is group-valued) one can easily
take the subgroup $\Gamma_w$ of $\Bbb R$ generated by
$w(E\setminus \{ 0 \})$ and get a quasi-valuation with value group
that is not torsion over $\Gamma_{v}$ (note for example that $\pi
\in \Gamma_w \setminus M_{w}$ in Example 5.2).

{\it We assume throughout sections 5 and 6 that $E$ is a finite
dimensional field extension with quasi-valuation $w$ extending $v$
on $F$ and $w(E\setminus \{ 0 \})$ is torsion over $\Gamma_{v}$
(which is the case in many of the examples). Note this implies
that $M_{w}$ is torsion over $\Gamma_{v}$. Therefore, by Lemma
5.1, $M_{w}$ is a group; we denote it as $\Gamma_{w}$.}

We note that some of the results presented in sections 5, 6 and 7
stand only in the case that $M_{w}$ is a group. Some of the
results remain true in the more general case that $M_{w}$ is a
monoid, but the proofs are much more complicated. Moreover, there
are some new interesting results when $M_{w}$ is required to be a
monoid; we shall present them in sections 8, 9 and 10.

%The following results (Lemma 3.5 to Theorem 3.6) are valid in the
%more general case where $E$ is a finite dimensional $F$-algebra
%(of course, under the assumption that $w(E\setminus \{ 0 \})$ is
%torsion over $\Gamma_{v}$).

 \proclaim{Lemma 5.3} Let $K$ be
a maximal ideal of $O_{w}$. Then $I_{w} \subseteq K$. \endproclaim

\demo{Proof} Assume to the contrary that there exists a nonzero
element $x \in O_{w}$, $x \notin K$, with $w(x)>0$. Then
$K+<x>=O_{w}$; i.e., there exist $m \in K$ and $y \in O_{w}$ such
that
 $$  (1) \ \  m+xy=1.$$
Since $w(xy)\geq w(x)+w(y)\geq w(x)>0$, one has $w(m)=0$.
Furthermore, $m^{-1} \notin O_{w}$ implies $w(m^{-1})<0$. Now,
multiplying Equation (1) by $m^{-1}$ we get
$$1+xym^{-1}=m^{-1}.$$ So, $w(xym^{-1})=w(m^{-1})<0$ (since
$w(m^{-1})<0$ and w(1)=0). Then,
$$w(m^{-1})=w(xym^{-1})\geq w(x)+w(y)+w(m^{-1})\geq w(m^{-1}).$$
Therefore, we have equality. Now, cancel $w(m^{-1})$ from both
sides and get $w(x)+w(y)=0$. Thus $w(x)=0$, a contradiction.

Q.E.D.\enddemo

%\proclaim {Lemma 3.4} $[O_{w}/I_{v}O_{w}:O_{v}/I_{v}]\leq
%n$.\endproclaim

%\demo {Proof} Let $\overline{a_{1}},...,\overline{a_{m}}\in
%O_{w}/I_{v}O_{w}$ be linearly independent over $O_{v}/I_{v}$. Let
%$a_{1},...,a_{m}$ be their representatives. We prove that
%$a_{1},...,a_{m}$ are linearly independent over $F$. Assume the
%contrary and apply Remark 2.2 to get $\sum_{i=1}^m
%\alpha_{i}a_{i}=0$ where $\alpha_{i} \in O_{v}$ and
%$\alpha_{i_{0}}=1$ for some $1 \leq i_{0} \leq m$. So, we have
%$\sum_{i=1}^m \overline{\alpha_{i}} \overline{a_{i}}=0$ where
%$\overline{\alpha_{i_{0}}}=\overline{1}$. This contradicts the
%linear independence of the $\overline{a_{i}}$'s.

%Q.E.D.\enddemo

\proclaim{Corollary 5.4} $[O_{w}/J_{w}:O_{v}/I_{v}]\leq
n$.\endproclaim %$[O_{w}/I_{w}:O_{v}/I_{v}]\leq n$;

\demo {Proof} By Lemma 5.3, $I_{w} \subseteq K$ for every maximal
ideal $K$ of $O_{w}$. Thus $I_{w} \subseteq J_{w}$, and we have
the natural epimorphism $$O_{w}/I_{w} \twoheadrightarrow
O_{w}/J_{w}.$$ The result now follows from Corollary 2.6.

%Note that $I_{w} \supseteq I_{v}O_{w}$ and thus we have the
%natural epimorphism $$O_{w}/I_{v}O_{w} \twoheadrightarrow
%O_{w}/I_{w}.$$ Regarding $J_{w}$, we have proved in Lemma 5.3 that
%$I_{w} \subseteq K$ for every maximal ideal $K$ of $O_{w}$.
%Therefore $I_{w} \subseteq J_{w}$, and again we have the natural
%epimorphism $O_{w}/I_{w} \twoheadrightarrow O_{w}/J_{w}$.

Q.E.D.\enddemo

\proclaim{Theorem 5.5} The quasi-valuation ring $O_{w}$ has $\leq
n$ maximal ideals.\endproclaim

 \demo {Proof} $O_{w}/J_{w}$ is  a semisimple ring. Take any set $\{K_{i}\}^{t}_{i=1}$ of
maximal ideals of $O_{w}$. Then $J_{w} \subseteq \cap_{i=1}^{t}
K_{i}$; so we have an epimorphism
$$O_{w}/J_{w} \twoheadrightarrow \bigoplus^{t}_{i=1}
O_{w}/K_{i}.$$ Thus by Corollary 5.4,
$$n \geq [O_{w}/J_{w}:O_{v}/I_{v}] \geq \sum_{i=1}^{t} [O_{w}/K_{i}:O_{v}/I_{v}] \geq \sum_{i=1}^{t} 1=t,$$
proving $t \leq n$.

 Q.E.D.\enddemo

Note that if $O_{w}$ is a valuation ring then $O_{w}$ has one
maximal ideal.

The proofs of the next two results are standard.

\proclaim{Lemma 5.6} Let $E/F$ be a field extension, $O_{v}$ a
valuation ring of $F$ with valuation ideal $I_{v}$ and $O_{u}$ a
valuation ring of $E$ containing $O_{v}$ with valuation ideal
$I_{u}$ containing $I_{v}$. Then $O_{u} \cap F=O_{v}$.

\endproclaim

\demo {Proof} Assume to the contrary, that there exists an $a \in
(O_{u} \cap F) \setminus O_{v}$. Then
$$a^{-1} \in I_{v}  \subseteq I_{u},$$ so $1=aa^{-1} \in
I_{u}$, a contradiction.

Q.E.D.\enddemo

\proclaim{Proposition 5.7} Suppose $E/F$ is a finite dimensional
field extension and let $R \subseteq E$ be a ring such that $R
\cap F=O_{v}$. Let $K$ be a maximal ideal of $R$ such that $K
\supseteq I_{v}$; then there exists a valuation $u$ of $E$ that
extends~ $v$ such that $O_{u} \supseteq R$ and $I_{u} \supseteq
K$.\endproclaim

\demo {Proof} We look at the pair $(R,K)$ and we take the
collection of all pairs $(R_{\alpha},I_{\alpha})$ where $R
\subseteq R_{\alpha} \subseteq E$ and $K  \subseteq I_{\alpha}
\vartriangleleft R_{\alpha} $. We order these pairs with the
partial order of containment, i.e.,
$$(R_{\alpha},I_{\alpha}) \leq (R_{\beta},I_{\beta}) \text{\ iff\
} R_{\alpha} \subseteq R_{\beta} \text{\ and\ } I_{\alpha}
\subseteq I_{\beta}.$$ Zorn's Lemma is applicable to yield a
maximal pair $(O_{u},I_{u})$. We claim that $O_{u}$ is a valuation
ring. Indeed, let $x \in E$ and assume to the contrary that $x,
x^{-1} \notin O_{u}$; then, by [Kap, p. 35, Th. 55], $I_{u}$ is a
proper ideal of $O_{u}[x]$ or $O_{u}[x^{-1}]$. Either way, we have
contradicted the maximality of $(O_{u},I_{u})$. Now, by Lemma 5.6,
$O_{u} \cap F=O_{v}$; i.e., the valuation $u$ which corresponds to
$O_{u}$ extends $v$.

Q.E.D.\enddemo

\demo{Definition 5.8} We call such a maximal pair $(O_{u},I_{u})$
obtained in Proposition 5.7 an {\it expansion} of $(R,K)$ to $E$.
In other words, $(O_{u},I_{u})$ is an expansion of $(R,K)$ to $E$
if $R \subseteq O_{u} \subseteq E$, $K \subseteq I_{u}$ ($I_{u}$
is a proper ideal of $O_{u}$) and for any pair $(R',K')$
satisfying $O_{u} \subseteq R' \subseteq E$ and $I_{u} \subseteq
K'$ ($K'$ is a proper ideal of $R'$), we have $O_{u}=R'$ and
$I_{u}=K'$. We suppress $I_{u}$ when it is not relevant.\enddemo

\proclaim{Corollary 5.9} There exists a valuation $u$ on the field
$E$ that extends $v$, such that $$O_{u} \supseteq
O_{w}.$$\endproclaim

\demo {Proof} Take $R=O_{w}$ in Proposition 5.7; take any maximal
ideal $K$ of $O_{w}$, and note that $K \supseteq I_{w}$ by Lemma
5.3.

Q.E.D.\enddemo

Note that one can expand the pair $(O_{w},K_{i})$ for every
maximal ideal $K_{i}$ of $O_{w}$. Also note that if $K_{i} \neq
K_{j}$ then $(O_{u_{i}},I_{u_{i}}) \neq (O_{u_{j}},I_{u_{j}})$
where $(O_{u_{i}},I_{u_{i}})$ is the expansion of $(O_{w},K_{i})$
and $(O_{u_{j}},I_{u_{j}})$ is the expansion of $(O_{w},K_{j})$.
There is a well known theorem from valuation theory (cf. [End, p.
97, Th. 13.7]) which says: let
$$\Cal U=\{O_{u} \subseteq E \mid O_{u} \text{\ is a valuation
ring of\ } E \text{\ and\ } O_{u} \cap F=O_{v} \};$$ then $|\Cal
U| \leq [E:F]_{sep}$.

We denote $\Cal K=\{$maximal ideals of $O_{w} \}$ and $$\Cal
U_{w}=\{O_{u} \subseteq E \mid O_{u} \text{\ is an expansion of\ }
(O_{w},K) \text{\ for some\ } K\in \Cal K \}.$$ We have:

 \demo{Remark 5.10} $\Cal U_{w} \subseteq \Cal U$.
\enddemo

\demo{Proof} Let $O_{u} \in \Cal U_{w}$; then $O_{u}$ is a
valuation ring of $E$ with valuation ideal $I_{u}$. Now, by the
definition of $\Cal U_{w}$, $O_{u} \supseteq O_{w} \supseteq
O_{v}$ and $I_{u} \supseteq K$ for some maximal ideal of $O_{w}$.
By Lemma 5.3, every maximal ideal $K \vartriangleleft O_{w}$
contains $I_{w}$ which contains $I_{v}$; thus $I_{u} \supseteq
I_{v}$. Therefore, by Lemma 5.6, $O_{u} \in \Cal U$.
\enddemo

\proclaim{Corollary 5.11} $\mid \Cal K \mid \leq \mid \Cal U_{w}
\mid \leq \mid \Cal U \mid \leq [E:F]_{sep}$.
\endproclaim

For a ring $R \subseteq E$ we denote by $I_{E}(R)$ the integral
closure of $R$ in $E$.

The following remark is known from valuation theory. (See, for
example, [End, p.69, Th. 10.8]):

\demo {Remark 5.12} Let $R$ be a subring of $E$, then $I_{E}(R)=
\bigcap\{O_{u}  \mid O_{u} $ is a valuation ring of $E$ containing
$R \} = \bigcap\{O_{u}  \mid O_{u} $ is a valuation ring of $E$
containing $R $, and $I_{u} \cap R$ is a maximal ideal of  $R\}$.
In particular, $I_{E}(O_{w})= \bigcap\{O_{u} \mid O_{u} $ is a
valuation ring of $E$ containing $O_{w}$, and $I_{u} \cap O_{w}$
is a maximal ideal of  $O_{w}\}$.\enddemo

\demo {Remark 5.13} $\Cal U_{w}=\{O_{u}  \mid O_{u} $ is a
valuation ring of $E$ containing $O_{w} $, and $I_{u} \cap O_{w}$
is a maximal ideal of  $O_{w}\}$. Indeed, $\subseteq$ is obvious.
Conversely, if $O_{u} \supseteq O_{w}$ is a valuation ring of $E$
and $I_{u} \cap O_{w}$ is a maximal ideal of $O_{w}$, then
$(O_{u},I_{u})$ is an element in the collection described in
Proposition 5.7. Assuming $(O_{u},I_{u}) \leq (O_{u'},I_{u'})$ for
some ring $O_{u'} \subseteq E$ and a maximal ideal $I_{u'}
\vartriangleleft O_{u'}$, we conclude that $O_{u'} $ is a
valuation ring of $E$. Thus, $O_{u} \subset O_{u'}$ implies
$I_{u'} \subset I_{u}$; so we must have $O_{u} = O_{u'}$ and
$I_{u} = I_{u'}$. I.e., $O_{u} \in \Cal U_{w}$. We conclude, by
Remark 5.12, that $$I_{E}(O_{w})= \bigcap_{O_{u} \in \Cal U_{w}}
O_{u}.$$

 %If $O_{u} \notin \Cal U_{w}$, then there exists
%a valuation ring $O_{u'} \supset O_{u}$; i.e., $(O_{u},I_{u})<
%(O_{u'},I_{u'})$ but then $I_{u'} \subset I_{u}$. We conclude, by
%Remark 5.12, that $I_{E}(O_{w})= \bigcap_{O_{u} \in \Cal U_{w}}
%O_{u}$.

\enddemo

\proclaim{Lemma 5.14} Let $O_{u} \in \Cal U_{w}$ and let $I_{u}$
be its maximal ideal; then \linebreak $I_{u} \cap I_{E}(O_{w})$ is
a maximal ideal of $I_{E}(O_{w})$.\endproclaim

\demo {Proof} $I_{u} \cap I_{E}(O_{w})$ is a prime ideal of
$I_{E}(O_{w})$; thus $(I_{u} \cap I_{E}(O_{w})) \cap O_{w}$ is a
prime ideal of $O_{w}$. However, $I_{u} \cap O_{w}$ is a maximal
ideal in $O_{w}$. (Indeed, $I_{u}$ contains a maximal ideal of
$O_{w}$). Now, by [Row, p.187, Corollary 6.33] (which says for $R$
integral over $C$ and $Q \in \text{Spec} (R)$ lying over $P \in
\text{Spec} (C)$, that $P$ is maximal in $C$ iff $Q$ is maximal in
$R$), $I_{u} \cap I_{E}(O_{w})$ is a maximal ideal of
$I_{E}(O_{w})$.

Q.E.D.\enddemo

Let $H$ be an isolated subgroup of a totally ordered abelian group
$\Gamma$. Recall that the quotient group $\overline{\Gamma} =
\Gamma / H$ is totally ordered, by setting $\overline{\gamma} \geq
\overline{0}$ iff $\gamma \geq h$ for some $h \in H$. There is a
one to one correspondence between the set of all prime ideals of a
valuation ring and the set $G(\Gamma)$ of all isolated subgroups
of $\Gamma$. % given by $P \mapsto \{ \alpha \in \Gamma \mid \alpha
%\neq v(p)$ and $\alpha \neq -v(p)$ for all $p \in P\}$.
 The rank of $\Gamma$ is the order type of $G(\Gamma) \setminus \{
\Gamma \}$ (see [End, p. 47]).

%For any totally ordered set $A$ and a subset $B \subseteq A$ the
%convex hull of $B$ in $A$ is defined by
%$$hull_{A}(B)=\{a \in A \mid \exists b_{1},b_{2} \in B \text{ such
%that } b_{1} \leq a \leq b_{2} \}.$$ In particular, for totally
We note that for totally ordered abelian groups $\Gamma \subseteq
\Delta$ and an isolated subgroup $H$ of $\Gamma$, the convex hull
of $H$ in $\Delta$, $hull_{\Delta}(H)$, is an isolated subgroup of
$\Delta$ (we call it the corresponding isolated subgroup of $H$ in
$\Delta$).

%For any subset $\Omega \subseteq \Gamma$, the convex hull of
%$\Omega$ in $\Gamma$ is defined as follows
%$${CH}_{\Gamma}(\Omega)=\{ \delta \in \Gamma \mid -n\gamma \leq
%\delta \leq n\gamma \text{ for some } \gamma \in \Omega \text{\
%with\ } \gamma \geq 0 \text{\ and\ } n \in \Bbb N\}.$$

Also recall that if $\Gamma$ and $\Delta$ are totally ordered
abelian groups such that $ \Gamma \subseteq \Delta  $ and $ \Delta
/ \Gamma$ is a torsion group, then $\text{rank} \Gamma =
\text{rank} \Delta$; i.e., there is a one to one correspondence
between the set of isolated subgroups of $\Gamma$ and the set of
isolated subgroups of $\Delta$. It is easy to see (when $ \Delta /
\Gamma$ is a torsion group) that
$$hull_{\Delta}(H)=\{ \delta \in \Delta \mid n\delta \in H \text{
for some }n \in \Bbb N\}.$$ See [End, Chapter 13] or [Bo, Section
4] for further discussion.

%$\Omega$, denoted $ {CH}_{\Gamma}(\Omega)$, is the isolated
%subgroup containing all $\delta \in \Gamma$ for which there exists
%$\gamma \in \Omega$ with $\gamma \geq 0$ and $n \in \Bbb N$ such
%that $-n\gamma \leq \delta \leq n\gamma$.

\proclaim{Theorem 5.15} Let $P$ be a prime ideal of $O_{v}$, $H
\leq \Gamma_{v}$ its corresponding isolated subgroup and
$hull_{\Gamma_{w}}(H) \leq \Gamma_{w}$ the corresponding isolated
subgroup. Let
 $f: \Gamma_{w} \rightarrow \Gamma_{w} / hull_{\Gamma_{w}}(H)$ be
the natural epimorphism and let $$\widetilde{w}:E \rightarrow
(\Gamma_{w} / hull_{\Gamma_{w}}(H)) \cup \{\infty \}$$ be a map
such that $\widetilde{w}(0) = \infty$ and $\widetilde{w}(x) =
f(w(x))=w(x) + hull_{\Gamma_{w}}(H)$ for $x \neq 0$. Then
$\widetilde{w}$ is a quasi-valuation on $E$ satisfying

(1) $O_{\widetilde{w}} \supseteq O_{w}$.

(2) $\widetilde{w}$ extends $\widetilde{v}$, the corresponding
valuation of $(O_{v})_{P}$, and $O_{\widetilde{w}}$ lies

\ \ \ \ over $(O_{v})_{P}$, i.e., $O_{\widetilde{w}} \cap F =
(O_{v})_{P}$.

(3) $O_{\widetilde{w}} = O_{w}S^{-1}$ where $S = O_{v} \setminus
P$.
\endproclaim

\demo{Proof} Note that $f$ is an epimorphism of ordered groups, so
$\alpha \leq \beta \in \Gamma_{w}$ implies \break $f(\alpha) \leq
f(\beta)$. Now, let $x,y $ be two nonzero elements of $E$. We
obviously have $w(xy) \geq w(x) + w(y)$ and thus
$$\widetilde{w}(xy)=f(w(xy)) \geq f(w(x) + w(y))
$$ $$= f(w(x)) + f(w(y)) = \widetilde{w}(x) + \widetilde{w}(y).$$
Next, assume $\widetilde{w}(x) \leq \widetilde{w}(y)$; we have
$$\widetilde{w}(x+y) = f(w(x+y)) \geq f(\min \{ w(x), w(y) \}).$$
Now, if $\widetilde{w}(x) = \widetilde{w}(y)$ then $f(w(x)) =
f(w(y))$ and thus $$f(\min \{ w(x), w(y) \})=f( w(x)) =
\widetilde{w}(x).$$ If $\widetilde{w}(x) < \widetilde{w}(y)$ then
$w(x)<w(y)$ and again $$f(\min \{ w(x), w(y) \})=f( w(x)) =
\widetilde{w}(x).$$ So, we have
$$\widetilde{w}(x+y)  \geq
\min \{ \widetilde{w}(x), \widetilde{w}(y) \}.$$ Note that if $x
\in O_{w}$ then $w(x) \geq 0$ and thus $f(w(x)) \geq
\overline{0}$, so $x \in O_{\widetilde{w}}$ and (1) is proved.

To prove (2) we denote by $g:\Gamma_{v} \rightarrow \Gamma_{v}/H$
the natural epimorphism and note that $\forall x \neq 0$,
$\widetilde{v}(x) = g(v(x))=v(x) + H$. Also note that one can view
$\Gamma_{v}/H$ inside $\Gamma_{w} / hull_{\Gamma_{w}}(H)$ via the
natural
monomorphism. %Thus $\forall x \neq 0$  in $F$,
%$$\widetilde{v}(x) = v(x) + H = v(x) + hull_{\Gamma_{w}}(H) $$ $$= w(x) +
%hull_{\Gamma_{w}}(H) = \widetilde{w}(x).$$
In this sense, $\widetilde{w}$ is an extension of $\widetilde{v}$.
Now, the fact that $O_{\widetilde{w}} \cap F = (O_{v})_{P}$
follows immediately.

We now prove (3). Let $rs^{-1} \in O_{w}S^{-1}$ then $w(s) = v(s)
\in H$. Write $w(s) = h$. We have $$w(rs^{-1}) = w(r) - w(s) \geq
-h \in H.$$ Hence $\widetilde{w}(rs^{-1}) \geq \overline{0}$ and
$rs^{-1} \in O_{\widetilde{w}}$. On the other hand, let $x \in
O_{\widetilde{w}}$; then $\widetilde{w}(x) \geq \overline{0}$
i.e., $w(x) \geq j$ for some $j \in hull_{\Gamma_{w}}(H)$. Note
that if $x \in O_{w}$ it is obvious that $x \in O_{w}S^{-1}$. So,
we may assume $x \in O_{\widetilde{w}} \setminus O_{w}$. We write
$w(x) = -h'$ for some $h' \in (hull_{\Gamma_{w}}(H))^{\geq 0}$.
(Indeed if $w(x) < -h'$ for all $h' \in
(hull_{\Gamma_{w}}(H))^{\geq 0}$ then $x \notin
O_{\widetilde{w}}$). We take $t \in \Bbb N$ such that $th' \in H$
and pick $s \in S$ such that $v(s) = th'$, then
$$w(xs) = w(x) + w(s) = (t-1)h' \geq 0.$$ Thus $xs \in
O_{w}$ and $x = (xs)s^{-1}$ for $xs \in O_{w}$ and $s \in O_{v}
\setminus P$.

Q.E.D.\enddemo

 We say that $R$ satisfies GU over $C$ if for any $P_{0}
\subseteq P_{1}$ in $\text{Spec} (C)$ and every $Q_{0} \in
\text{Spec} (R)$ lying over $P_{0}$ there is $Q_{0} \subseteq
Q_{1}$ in $\text{Spec} (R)$ lying over $P_{1}$; cf. [Row, p. 185]

\proclaim{Theorem 5.16} $O_{w}$ satisfies GU over
$O_{v}$.\endproclaim

\demo {Proof} Let $P_{0} \subseteq P_{1} \in \text{Spec} (O_{v})$
and assume $Q_{0} \in \text{Spec} (O_{w})$ such that $Q_{0} \cap
O_{v}=P_{0}$. We prove that there exists a $Q_{0} \subseteq Q \in
\text{Spec} (O_{w})$ such that $Q \cap O_{v}=P_{1}$. By [Row,
p.186, Lemma 6.30] there exists an ideal $Q \supseteq Q_{0}$,
maximal with respect to $Q \cap O_{v} \subseteq P_{1}$ and any
such $Q$ is in $\text{Spec} (O_{w})$. We shall prove that $Q \cap
O_{v}=P_{1}$. Denote $S_{1}=O_{v} \setminus P_{1}$; then, by
Theorem 5.15, $O_{w}S_{1}^{-1}$ has a quasi-valuation
$\widetilde{w}$ with value group such that $\widetilde{w}$ extends
the valuation corresponding to $O_{v}S_{1}^{-1}$. Now, from the
maximality of $Q$ and the injective order preserving mapping
between ideals of $O_{w}S_{1}^{-1}$ and ideals of $O_{w}$ disjoint
from $S_{1}$,
%one to one correspondence between prime ideal of $O_{v}$ disjoint
%from $S_{1}$ and prime ideals of $O_{w}S_{1}^{-1}$,
we deduce that $QS_{1}^{-1}$ is a maximal ideal of
$O_{w}S_{1}^{-1}$. However, by Lemma 5.3, every maximal ideal of
$O_{w}S_{1}^{-1}$ must contain $P_{1}S_{1}^{-1}$. Therefore $Q
\cap O_{v}=P_{1}$.

Q.E.D.\enddemo

We are now able to deduce Corollary 4.13 (in the special case
where the value monoid of the quasi-valuation is a group) without
using the going down property.

\proclaim{Corollary 5.17} K-dim$O_{w} =$K-dim$O_{v}$.
\endproclaim

\demo {Proof} By Remark 3.5 and Theorem 3.7, $O_{w}$ satisfies INC
over $O_{v}$; by Theorem 5.16, $O_{w}$ satisfies GU over $O_{v}$.
Thus by [Row, p.185, Remark ~6.26], K-dim$O_{w} =$K-dim$O_{v}$.

Q.E.D.\enddemo

It is known from valuation theory (see [End, p.102, Th. 13.14])
that for every valuation ring $O_{u} \subseteq E$ such that $O_{u}
\cap F = O_{v}$, one has $$\text{K-dim}O_{v}=\text{K-dim}O_{u}.$$
Also, it is well known that if $R$ is integral over $C$, then
K-dim$R=$K-dim$C$. So, we have the following Theorem:

\proclaim{Theorem 5.18}
K-dim$O_{v}=$K-dim$O_{w}=$K-dim$I_{E}(O_{w})=$K-dim$O_{u}$.\endproclaim

At this point, in preparation for computing Krull dimension in the
next theorem, we pause to build another quasi-valuation $w'$ which
arises naturally from $w$. Let $O_{u}$ be a valuation ring of $E$
such that $O_{u} \cap F=O_{v}$. We first note that since $E$ is
algebraic over $F$, $\Gamma_{u}$ is torsion over $\Gamma_{v}$ (see
[End, p.99-100, Thm.13.9 and 13.11]). We recall our assumption
that $w(E\setminus \{ 0 \})$ is torsion over $\Gamma_{v}$.
Therefore $\Gamma_{u}$ and $\Gamma_{w}$ can be embedded in
$\Gamma_{\text{div}}$, where $\Gamma_{\text{div}}$ is the
divisible hull of
$\Gamma_{v}$. % i.e., the quotient of $\Gamma_{v} \times (\Bbb N
%\setminus \{ 0 \})$ by the equivalence relation $\sim$ defined by
%$(\gamma , n) \sim (\delta, m)$ iff $n \delta = m \gamma$, endowed
%with the usual addition and total ordering:
%$$(\gamma , n)+ (\delta, m)=(m \gamma+ n \delta,mn)$$ and
%$$(\gamma , n) \leq (\delta, m) \text{\ iff     }  m \gamma \leq n
%\delta.$$ So, $\Gamma_{\text{div}}$ is an ordered group.
 Note that $\Gamma_{\text{div}} / \Gamma_{v}$ is a torsion group, while on
the other hand, for any totally ordered abelian group $\Delta$
containing $\Gamma_{v}$ and such that $\Delta / \Gamma_{v}$ is a
torsion group, there is an embedding of $\Delta$ in
$\Gamma_{\text{div}}$. For more details on $\Gamma_{\text{div}}$
see [End, p.100-102].

Now, we define $w': E \rightarrow \Gamma_{\text{div}} \cup \{
\infty \}$ by
$$w'(x)= \min \{ u(x) \mid u \text{\ is a valuation whose valuation ring
is in\ } \Cal U_{w} \} .$$ Note that $\mid \Cal U_{w} \mid <
\infty$, so $w'$ is well defined. Also note that $w'$ is an
exponential quasi-valuation on $E$ that extends $v$ and $O_{w'}$
is the intersection of all $O_{u} \in \Cal U_{w}$. Thus, by Remark
5.13, $O_{w'}=I_{E}(O_{w})$.

\proclaim{Theorem 5.19} Let $O_{v} \subseteq R \subseteq E$ be a
ring whose Jacobson radical satisfies $J \supseteq I_{v}$. Then,
K-dim$R=$K-dim$O_{v}$, and $R$ has a finite number of maximal
ideals.\endproclaim

\demo {Proof} Let $\Cal U_{R}= \{ O_{u} \mid O_{u} $ is a
valuation ring of $E$ containing $R $ and $I_{u} \cap R$ is a
maximal ideal of $R\}$. By Remark 5.12 we have $$I_{E}(R)=
\bigcap_{O_{u} \in \Cal U_{R}} O_{u}.$$

%By Remark 5.12 we have $I_{E}(R)= \bigcap\{O_{u} \mid O_{u} $ is a
%valuation ring of $E$ containing $R $ and $I_{u} \cap R$ is a
%maximal ideal of  $R\}$.

Now, Let $O_{u} \in \Cal U_{R}$ and let $I_{u}$ denote its maximal
ideal; since $I_{u} \cap R$ is a maximal ideal of $R$, we get
$I_{u} \cap O_{v}=I_{v}$ (since $I_{u} \supseteq J \supseteq
I_{v}$) and thus by Lemma 5.6, $O_{u} \cap F=O_{v}$. Therefore,
$\Cal U_{R} \subseteq \Cal U$.

As above, we define $w'': E \rightarrow \Gamma_{\text{div}} \cup
\{ \infty \}$ by
$$w''(x)= \min \{ u(x) \mid u \text{\ is a valuation on $E$ whose valuation ring
is in\ } \Cal U_{R} \} .$$ $w''$ is well defined since $\mid \Cal
U_{R} \mid < \infty$. Also, $w''$ is an exponential
quasi-valuation on $E$ that extends $v$ and since $O_{w''}$ is the
intersection of all $O_{u} \in \Cal U_{R}$, we have
$O_{w''}=I_{E}(R)$.

 %Therefore, $I_{E}(R)=\bigcap \{ O_{u} \mid O_{u} $
%is an expansion of $(R,K)$ for some maximal ideal $K$ of $R \}$.
%Thus, $I_{E}(R)$ is the ring corresponding to the exponential
%quasi-valuation $w'$ extending $v$ (namely, $I_{E}(R)=O_{w'}$
%where $w'$ is the quasi-valuation that was defined before Theorem
%5.19).

By Corollary 4.13 (or Corollary 5.17) we get,
$$\text{K-dim}I_{E}(R)=\text{K-dim}O_{v}$$ and since $I_{E}(R)$ is
integral over $R$, we get
$$\text{K-dim}I_{E}(R)=\text{K-dim}R.$$

A similar proof to the proof of Remark 5.13 shows that
$$\Cal U_{R}=\{ O_{u} \mid O_{u} \text{ is an expansion of } (R,K)
\text{ for some maximal ideal } K \text{ of } R \}.$$

Now assume that $R$ has an infinite number of maximal ideals; then
we get an infinite number of $O_{u}$'s all of which lie over
$O_{v}$, a contradiction.

Q.E.D.\enddemo

%Note: it is not difficult to see that $\Cal U_{R}$ as defined in
%the previous Theorem equals $\{ O_{u} \mid O_{u} $ is an expansion
%of $(R,K)$ for some maximal ideal $K$ of $R \}$; so, the notation
%$\Cal U_{R}$ is consistent with $\Cal U_{w}$ (actually, it is
%consistent with $\Cal U_{O_{w}}$, but to simplify notation we
%write $\Cal U_{w}$).

We note that since $I_{E} (O_{w})$ is integral over $O_{w}$, there
is a surjective mapping $\text{Spec} (I_{E}(O_{w})) \rightarrow
\text{Spec} (O_{w})$ given by $K \rightarrow K \cap O_{w}$ sending
$\{$maximal ideals of $I_{E}(O_{w})\} \rightarrow \{$maximal
ideals of $O_{w}\}$. See [Row, p. 184-189] for further discussion.

\proclaim {Lemma 5.20} $\Cal U_{w}=\Cal U'_{w}$ where $\Cal
U'_{w}=\{O_{u} \mid O_{u}$ is an expansion of $(I_{E}(O_{w}),K)$
for some maximal ideal $K$ of $I_{E}(O_{w})\}$.\endproclaim

\demo {Proof} Let $O_{u} \in \Cal U'_{w}$. Then $O_{u} \supseteq
I_{E}(O_{w}) \supseteq O_{w}$, and $I_{u} \supseteq K' \supseteq
K$ for maximal ideals $K'$ and $K$ in $I_{E}(O_{w})$ and $O_{w}$
respectively. So, $O_{u} \in \Cal U_{w}$. Conversely, if $O_{u}
\in \Cal U_{w}$, then $I_{u} \cap I_{E}(O_{w})$ is a maximal ideal
of $I_{E}(O_{w})$ by Lemma 5.14, so $(I_{E}(O_{w}),I_{u} \cap
I_{E}(O_{w}))$ has an expansion $(O_{u},I_{u})$.

Q.E.D\enddemo

%For every $P \in \text{Spec}(O_{v})$ denote $\Cal Q_{P}= \{ Q \in
%\text{Spec} (O_{w'})| Q \cap O_{v} = P \}$. Now,
Let $P  \vartriangleleft O_{v}$ be a prime ideal and let $S=O_{v}
\setminus P$; then $S^{-1}O_{v}$ is a valuation ring of $F$ with
valuation ideal $S^{-1}PO_{v}$. Note that $S^{-1}O_{w} \supseteq
S^{-1}O_{v}$ and every maximal ideal of $S^{-1}O_{w}$ contains
$S^{-1}PO_{v}$; thus by Theorem 5.19,
K-dim$S^{-1}O_{v}=$K-dim$S^{-1}O_{w} $. (Note that one can also
deduce this equation by Theorem 5.15 and Corollary 4.13 (or
Corollary 5.17)). Let $Q \vartriangleleft O_{w}$ be a prime ideal
lying over $P$. We expand $(S^{-1}O_{w},S^{-1}QO_{w})$ for every
prime ideal $Q$ lying over $P$, to get valuation rings lying over
$S^{-1}O_{v}$. Note that if $Q \neq Q'$ lie over $P$ then we get
different valuation rings. Indeed, assume to the contrary that we
get $(U,I)$ from both of the expansions; then $I \cap
S^{-1}O_{w}$, which is a proper ideal of $S^{-1}O_{w}$, contains
 the maximal ideals $S^{-1}QO_{w}$ and $S^{-1}Q'O_{w}$, a
contradiction. So, we have $$(*) \mid \{ Q \mid Q \text{ is a
prime ideal of } O_{w} \text{ lying over } P \} \mid$$  $$\leq
\mid \{ U \mid U \text{ is an expansion of } (S^{-1}O_{w},K)
\text{ for some maximal ideal } K \text{ of } S^{-1}O_{w} \}
\mid$$
$$
 \leq [E:F]_{sep}.$$ The last inequality is since, as before, every such $U$ lies
 over the valuation ring
$S^{-1}O_{v}$. So, we have the following theorem:

\proclaim{Theorem 5.21} $$\text{K-dim}O_{v} \leq \mid \text{Spec}
(O_{w}) \mid \leq \mid \text{Spec}(I_{E}(O_{w})) \mid \leq
[E:F]_{sep} \cdot\text{K-dim}O_{v}.$$\endproclaim

\demo {Proof} The first inequality is due to Lemma 3.12. The
second inequality is due to the fact that $I_{E}(O_{w})$ is
integral over $O_{w}$. For the last inequality, note that
$I_{E}(O_{w})=O_{w'}$ where $w'$ is the quasi-valuation extending
$v$ defined by $w'(x)= \min \{ u(x) \mid u \text{\ is a valuation
whose valuation ring is in\ } \Cal U_{w} \} $ for all $x \in E$.
($w'$ was defined before Theorem 5.19). For every $P \in
\text{Spec}(O_{v})$ denote $$\Cal Q_{P}= \{ Q \in \text{Spec}
(O_{w'})| Q \cap O_{v} = P \};$$ thus
$\text{Spec}(O_{w'})=\bigcup_{P \in \text{Spec}(O_{v})} \Cal
Q_{P}$. (In fact $\{ \Cal Q_{P}\}_{P \in \text{Spec}(O_{v})}$ is a
partition of $\text{Spec}(O_{w'})$.) Finally, by $(*)$ we get
$$\mid \bigcup_{P \in \text{Spec}(O_{v})} \Cal
Q_{P} \mid \leq  \text{K-dim}O_{v} \cdot [E:F]_{sep}.$$

%, and from the
%fact that $\text{K-dim}O_{v}=\text{K-dim}O_{w'}$ (Corollary 5.17)

Q.E.D.\enddemo

\heading \S 6 a Bound on the Quasi-Valuation and the Height
Formula\endheading

In this section we continue to assume that $w(E\setminus \{ 0 \})$
is torsion over $\Gamma_{v}$ (and thus the value monoid of the
quasi-valuation is a group). We show that any such quasi-valuation
on $E$ extending $v$ is dominated by any valuation $u$ on $E$
extending $v$ satisfying $O_{u} \supseteq O_{w}$. We also prove
that the quasi-valuation ring satisfies the height formula.

\proclaim{Lemma 6.1} Let $w$ be a quasi-valuation on $E$ extending
$v$ and let $u$ be a valuation on $E$ extending $v$ such that
$O_{u} \supseteq O_{w}$; then $I_{u} \supseteq I_{w}$.\endproclaim

\demo{Proof} $K = I_{u} \cap O_{w}$ is a prime ideal of $O_{w}$.
Assuming $K$ is not maximal in $O_{w}$, one has $K \cap O_{v}$ not
maximal in $O_{v}$ (by INC). This contradicts the fact that $I_{u}
\cap O_{v} = I_{v}$. So $I_{u} \cap O_{w}$ is a maximal ideal in
$O_{w}$ and thus (Lemma 5.3) contains $I_{w}$.

Q.E.D.\enddemo

The following lemma is valid without the assumption that
$\Gamma_{w}$ and $\Gamma_{w'}$ are torsion over $\Gamma_{v}$.

 \proclaim{Lemma 6.2}  Let $w$ and $w'$ be quasi-valuations
on $E$ extending $v$ such that $\Gamma_{w}$ and $\Gamma_{w'}$
embed in a group $G$ and assume that $w$ is an exponential
quasi-valuation. If there exists $x \in E$ such that $w(x)<w'(x)$,
then there exists $y \in E$ such that $0=w(y)<w'(y)$.\endproclaim

\demo{Proof}  Write $\sum_{i=0}^n \alpha_{i}x^{i} = 0$ for
$\alpha_{i} \in F$. Since the sum is zero we must have $k,l \in
\Bbb N$ such that
$$w(\alpha_{k}x^{k}) = w(\alpha_{l}x^{l}).$$ Assuming $k < l$ we
have, using the fact that $w$ is exponential,
$w(x^{l-k}\alpha_{l}\alpha_{k}^{-1}) = ~0$. Now, since $w(x)
<w'(x)$, we have %$w'(x^{l-k}) > w(x^{l-k})$. (Indeed
$$w'(x^{l-k}) \geq (l-k)w'(x) > (l-k)w(x) = w(x^{l-k}).$$ Therefore
$w'(x^{l-k}\alpha_{l}\alpha_{k}^{-1}) >
w(x^{l-k}\alpha_{l}\alpha_{k}^{-1}) = 0$. So take
$y=x^{l-k}\alpha_{l}\alpha_{k}^{-1}$.

Q.E.D.\enddemo

%Note that in Lemma 6.2 we do not have to assume that $\Gamma_{w}$
%and $\Gamma_{w'}$ are torsion over $\Gamma_{v}$.

We now show that the values of the quasi-valuation are bounded by
the values of a valuation.

\proclaim{Theorem 6.3} Let $w$ be a quasi-valuation on the field
$E$ extending a valuation $v$. Then there exists a valuation $u$
extending $v$ such that $u$ dominates $w$ i.e., $\forall x \in E$,
$w(x) \leq u(x)$. Moreover, for every valuation $u$ on $E$
extending $v$ satisfying $O_{u} \supseteq O_{w}$, $u$ dominates
$w$.\endproclaim

\demo{Proof} First note that $\Gamma_{u}$ and $\Gamma_{w}$ embed
in $\Gamma_{\text{div}}$, so we refer to the ordering in
$\Gamma_{\text{div}}$. By Corollary 5.9, there exists a valuation
$u$ that extends $v$, such that $O_{u} \supseteq O_{w}$. We now
prove that every such $u$ dominates $w$.

Let $x \in E$ and assume to the contrary that $u(x) < w(x)$. Note
that $u$ is a valuation and thus in particular an exponential
quasi-valuation. Therefore, by Lemma 6.2 there exists $y \in E$
such that $0=u(y)<w(y)$, i.e., $y \in I_{w} \setminus I_{u}$, in
contradiction to Lemma 6.1.

Q.E.D. \enddemo

%Let $R$ be a ring; $R$ is said to satisfy the height formula if
%for any ideal $P \in \text{Spec} (R)$, we have
%K-dim$R=h_{R}(P)+$K-dim$(R/P)$.

%The following lemma is a well known result, seen by matching
%chains of prime ideals:

%\proclaim{Lemma 4.4} If $C \subseteq R$ are commutative rings such
%that $C$ satisfies the height formula and $R$ satisfies GU, GD and
%INC over $C$, then $R$ satisfies the height formula.

%\endproclaim

The following theorem is a special case of Lemma 4.16 (by taking
$R=O_{w}$ and $C=O_{v}$) in view of the previous results. We shall
prove it here for the reader's convenience.

\proclaim{Theorem 6.4} $O_{w}$ satisfies the height formula.
\endproclaim

\demo{Proof} First note that $O_{v}$ satisfies the height formula
because $\text{Spec} (O_{v})$ is linearly ordered. Also note that
if K-dim$O_{v}=\infty$ then by Corollary 4.13 (or Corollary 5.17),
K-dim$O_{w}=\infty$. In this case, if $Q \in \text{Spec} (O_{w})$
has infinite height, then obviously
K-dim$O_{w}=h_{w}(Q)+$K-dim$(O_{w}/Q)$; if $Q \in \text{Spec}
(O_{w})$ has finite height $t$ then $P=Q \cap O_{v} \in
\text{Spec} (O_{v})$ has finite height $t$ by INC and GD and thus
by GU, K-dim$(O_{w}/Q)=\infty$. So, we may assume that
K-dim$O_{v}=m$ and by Corollary 4.13 (or Corollary 5.17), we have
K-dim$O_{w}=m$. Now, let $Q $ be a prime ideal of $O_{w}$ and $P=Q
\cap O_{v} \in \text{Spec} (O_{v})$; then, by Lemma 4.15,
$h_{w}(Q)=h_{v}(P)$. Write $h_{w}(Q)=t$, $Q=Q_{t}$, $P=P_{t}$.  We
have a chain $$P_{t} \subset P_{t+1} \subset ... \subset
P_{m}=I_{v}$$ of size $m-t+1$ and by GU, we have a chain
$$Q_{t} \subset Q_{t+1} \subset ... \subset Q_{m}.$$ Assume to the
contrary that there exists a longer chain. Let $Q'_{t} \subset
Q'_{t+1} \subset ... \subset Q'_{l}$ be a chain of longer size. We
may assume that $Q'_{t}=Q_{t}$. Therefore, by INC, there exists a
chain $$P'_{t} \subset P'_{t+1} \subset ... \subset P'_{l}$$ where
$P'_{t}=P_{t}$ and thus a chain
$$P'_{0} \subset P'_{1} \subset ... \subset P'_{l} \subseteq
I_{v}$$ of size $>m+1$, a contradiction.

Q.E.D.\enddemo

%We note that [Sa, Proposition 4.7] can be seen to be a consequence
%of Proposition 5.17 because under the hypothesis that $M_{w}$ is
%torsion over the group $\Gamma_{v}$, the positive part of
%$\sqrt{H}$ is the unique PIM lying over $H^{\geq0}$ and thus
%equals $\overline{hull_{M}(H^{\geq0})}$

\proclaim{Proposition 6.5} Let $P \in \text{Spec} (O_{v})$, $H$
its corresponding isolated subgroup in $\Gamma_{v}$, and
$hull_{\Gamma_{w}}(H)$ the corresponding isolated subgroup in
$\Gamma_{w}$. Denote $$I = \{ x \in O_{w}| w(x) \notin
hull_{\Gamma_{w}}(H) \}$$ and let $ \Cal Q    = \{ Q \in
\text{Spec}(O_{w})| Q \cap O_{v} = P \}$. Then
$$\sqrt{I} = \bigcap_{Q \in \Cal Q} Q.$$
\endproclaim

\demo{Proof} Under the hypothesis that $w(E\setminus \{ 0 \})$ is
torsion over $\Gamma_{v}$ (and thus the value monoid of the
quasi-valuation is a group), $hull_{\Gamma_{w}}(H^{\geq0})$ is the
unique PIM lying over $H^{\geq0}$ and thus equals
$\overline{hull_{\Gamma_{w}}(H^{\geq0})}$. Now, use Proposition
4.19.

Q.E.D.\enddemo

\proclaim{Corollary 6.6} $\sqrt{I_{w}} = J_{w}$. \endproclaim

\demo{Proof} First note that by Lemma 5.3 and INC (Remark 3.5 and
Theorem 3.7), the set of prime ideals lying over $I_{v}$ is
exactly the set of maximal ideals of $O_{w}$. Thus, the corollary
is a special case of the previous proposition, seen by taking $P =
I_{v}$.

Q.E.D.\enddemo

\heading \S 7 Exponential Quasi-Valuations \endheading

%We start with the following definition:
In this section we study exponential quasi-valuations extending a
valuation in a finite dimensional field extension. We show that
the associated rings are integrally closed; we also show that if
two exponential quasi-valuations are not equal then their rings
cannot be equal. Finally, we prove that exponential
quasi-valuations have a unique form and we obtain a bound on the
number of these quasi-valuations.

%\demo {Definition 5.1} A monoid $M$ is said to be {\it weakly
%cancellative} if for any $a,b \in M$, $a+b=a$ implies $b=0$.

%\enddemo

\demo {Definition 7.1} An additive monoid $M$ is called {\it
weakly cancellative} if for any $a,b \in M$, $a+b=a$ implies
$b=0$.\enddemo

Note that, by Lemma 5.1, every torsion monoid over a group is a
group and thus in particular weakly cancellative.

 In this section we do not assume that $w(E\setminus \{ 0 \})$
is torsion over $\Gamma_{v}$. Instead, we assume the weaker
hypothesis that $M_{w}$ (the value monoid) is weakly cancellative.

Recall from Definition 1.7 that a quasi-valuation $w$ is called
exponential if for every $x \in E$, $$w(x^{n})=nw(x), \ \ \
\forall n \in \Bbb N.$$

%The following remark shows that it implies that $M_{w}$ is torsion
%over $\Gamma_{v}$.

\proclaim {Lemma 7.2} If $w$ is an exponential quasi-valuation
then for each nonzero $ b \in E$ there exists $t=t(b)$, $1 \leq t
\leq n$, such that $w(b^{t}) \in \Gamma_{v}$, and $$ \{
w(b^{i})+\Gamma_{v} \mid i \in \Bbb N \}=\{ w(b^{i})+\Gamma_{v}
\mid 0 \leq i < t \}.$$\endproclaim

\demo {Proof} By Lemma 2.1, the set $\{ 1,b,...,b^{n}\}$ satisfies
$$w(b^{j}) + \Gamma_{v}= w(b^{i})+\Gamma_{v}$$ for some $0 \leq i
< j \leq n$. Therefore $jw(b) + \Gamma_{v}= iw(b)+\Gamma_{v}$ and
we get $jw(b) + \alpha= iw(b)$ for some $\alpha \in \Gamma_{v}$.
Now, add $(j-i)w(b)$ to both sides and get $$jw(b) + \alpha+
(j-i)w(b)= jw(b).$$ Thus (since $M_{w}$ is weakly cancellative) we
get $\alpha+ (j-i)w(b)=0$, i.e., $w(b^{j-i})=(j-i)w(b) \in
\Gamma_{v}$. As for the second part, let $t \leq k \in \Bbb N$ and
write $k=qt+r$ for $0 \leq r < t$. We have:
$$w(b^{k}) =w(b^{qt+r})= w(b^{qt})+w(b^{r}) \in w(b^{r}) +
\Gamma_{v}.$$

Q.E.D.\enddemo

%\proclaim{Corollary 5.2'} Let $w$ be an exponential
%quasi-valuation. The following are equivalent:

%(1). $w(E\setminus \{ 0 \})$ is torsion over $\Gamma_{v}$.

%(2). The value monoid, $M_{w}$, is a group.
%\endproclaim

%\demo {Proof} (1)$\Rightarrow$(2). True for any quasi-valuation
%due to Lemma 3.1. (2)$\Rightarrow$(1). Let $0 \neq x \in E$, by
%Lemma 7.2 there exists a $1 \leq t=t(b) \leq n$ such that
%$w(x^{t}) \in \Gamma_{v}$. Therefore $tw(x)=w(x^{t})\in
%\Gamma_{v}$ and $w(E \setminus \{0\})$ is torsion over
%$\Gamma_{v}$.
%\enddemo

\proclaim {Lemma 7.3} If $w$ is an exponential quasi-valuation and
$M_{w}$  is weakly cancellative, then $[M_{w}:\Gamma_{v}]$ is
finite, and consequently $M_{w}$ is a group. \endproclaim

\demo {Proof} First note that every element of $M_{w}$ is of the
form $~\sum_{i=1}^{k} w(x_{i}) \  $ for $x_{i} \in E$. By Lemma
2.1, there exists a set $\{ b_{1},...,b_{m} \}$ of $E$, where $ m
\leq n$, such that \break $\{ w(b_{i})+\Gamma_{v} \mid 1\leq i
\leq m \}$ comprises all the cosets of $\Gamma_{v}$ in
$w(E\setminus \{ 0 \})$. Thus, one can replace every coset
$\sum_{i=1}^{k} w(x_{i}) + \Gamma_{v}$ by $\sum_{i=1}^{m}
n_{i}w(b_{i}) + \Gamma_{v}$ for appropriate $n_{i}$'s. Now, since
$w$ is exponential,
 $$\sum_{i=1}^{m} n_{i}w(b_{i}) + \Gamma_{v} =
\sum_{i=1}^{m} w(b_{i}^{n_{i}}) + \Gamma_{v}.$$  Letting
$t_{i}=t(b_{i})$ in Lemma 7.2, we see that elements of the form
$\sum_{i=1}^{k} w(x_{i}) $ lie in $\leq \sum_{i=1}^{m} t_{i}$
cosets over $\Gamma_{v}$ and therefore
$$[M_{w}:\Gamma_{v}] \leq \sum_{i=1}^{m} t_{i}.$$

Now, using the fact that $[M_{w}:\Gamma_{v}] < \infty $ and the
weak cancellation property we show that $M_{w}$ is torsion over
$\Gamma_{v}$. Let $\delta \in M_{w}$. Since ~$[M_{w}:~\Gamma_{v}]
< ~\infty $, there exist $i<j \in \Bbb N$ such that
$i\delta+\Gamma_{v}=j\delta+\Gamma_{v}$; i.e., there exists
$\alpha \in \Gamma_{v}$ such that $i\delta=j\delta+\alpha$. Thus
$i\delta=i\delta+(j-i)\delta+\alpha$ and by the weak cancellation
property we have $(j-i)\delta+\alpha=0$, i.e.,
$(j-i)\delta=-\alpha \in \Gamma_{v}$. Therefore, $M_{w}$ is
torsion over $\Gamma_{v}$, and then by Lemma 5.1, $M_{w}$ is a
group.

%Now, the fact that $M_{w}$ is a group follows easily. Indeed,
%$M_{w}$ is torsion over $\Gamma_{v}$ by the weakly cancellative
%property, and then by Lemma 5.1, $M_{w}$ is a group.

% Indeed, let $0 \neq m \in M_{w}$; we look at the set $\{ km
%\}_{1 \leq k \leq (\sum_{i=1}^{m} t_{i})+1}$. There exist \break
%$1 \leq s<t \leq (\sum_{i=1}^{m} t_{i})+1$ such that
%$tm+\Gamma_{v}=sm+\Gamma_{v}$ and thus $(t-s)m \in \Gamma_{v}$.
%Therefore $m$ is invertible in $M_{w}$.

Q.E.D.\enddemo

\proclaim{Corollary 7.4} Let $w$ be an exponential quasi-valuation
with $M_{w}$ not necessarily weakly cancellative. The following
are equivalent:

(1) $w(E\setminus \{ 0 \})$ is torsion over $\Gamma_{v}$.

(2) $M_{w}$ is a torsion group over $\Gamma_{v}$.

(3) $M_{w}$ is weakly cancellative.
\endproclaim

\demo {Proof} (1)$\Rightarrow$(2). True for any quasi-valuation
since if $w(E\setminus \{ 0 \})$ is torsion

\ \ \ \ \ \ over $\Gamma_{v}$, then $M_{w}$ is torsion over
$\Gamma_{v}$. Now, use Lemma 5.1.

\ \ \ \ \ \ (2)$\Rightarrow$(3). Every group satisfies the weak
cancellation property.

\ \ \ \ \ \  (3)$\Rightarrow$(2). By Lemma 7.3.

\ \ \ \ \ \ (2)$\Rightarrow$(1). Obvious.

Q.E.D.\enddemo

%\proclaim{Corollary 7.4} If $w$ is an exponential quasi-valuation
%then $M_{w}$ is a group.\endproclaim

%\demo {Proof} Let $0 \neq x \in E$, by Lemma 5.2 there exists a $1
%\leq t=t(b) \leq n$ such that $w(x^{t}) \in \Gamma_{v}$. Therefore
%$tw(x)=w(x^{t})\in \Gamma_{v}$ and $w(E \setminus \{0\})$ is
%torsion over $\Gamma_{v}$; thus $M_{w}$ is torsion over
%$\Gamma_{v}$. Now, use Lemma 5.1.\enddemo

We note at this point that since $M_{w}$ is a torsion group over
$\Gamma_{v}$, one can use the results obtained in sections 5 and
6.

A natural question that arises now is the connection between the
quasi-valuation ring $O_{w}$ and the integral closures (of $O_{v}$
and of $O_{w}$) inside $E$. %We have the following Lemma:

\proclaim{Lemma 7.5} If $w$ is an exponential quasi-valuation,
then $O_w$ contains $I_E(O_{v})$. In fact, $O_w=I_E(O_{w})$.
\endproclaim

\demo{Proof} Let $x \in I_E(O_{w})$ and write $x^{n}+
\sum_{i=0}^{n-1} \alpha_{i} x^{i}=0 $ for $\alpha_{i} \in O_{w}$.
Then $$nw(x)=w (x^{n}) \geq \min_{0 \leq i \leq n-1} \{
w(\alpha_{i} x^{i})\}$$ $$ = w(\alpha_{i_{0}} x^{i_{0}}) \geq
w(\alpha_{i_{0}})+ i_{0}w(x)$$ for an appropriate $i_{0}$.
Therefore $$(n-i_{0})w(x) \geq w(\alpha_{i_{0}})  \geq 0,$$ i.e.,
$w(x) \geq 0$.

Q.E.D
\enddemo

 However, if $w$ is not exponential, $O_{w}$ does not
necessarily contain $I_E(O_{v})$. For example, take $F=\Bbb Q$,
$v=v_{p}$, $M=\Bbb Z$, $e=i$ and $\gamma = 1$ in Example ~2.9.
Note that $i \notin O_{w}$ whereas $x^{2}+1 \in O_v[x]$.

\proclaim{Lemma 7.6} Let $w$ and $w'$ be two exponential
quasi-valuations on $E$ extending $v$ on $F$ with value groups
$M_{w}$ and $M_{w'}$. If $O_{w'}=O_{w}$ then $w'=w$.
\endproclaim

\demo{Proof} First note that by Lemma 7.3, both $M_{w}$ and
$M_{w'}$ are actually torsion groups over $\Gamma_{v}$ and thus
can be embedded in $\Gamma_{\text{div}}$. Assume to the contrary
that there exists an $x \in E$ such that $w'(x)>w(x)$; then, By
Lemma 6.2 we get an element $y \in E$ with $w'(y)
> w(y) =0$. Thus, $y \in I_{w'} \setminus I_{w}$. Now, by
Corollary 6.6, $\sqrt{I_{w}} = J_{w}$ and $\sqrt{I_{w'}} =
J_{w'}$; also, since $w$ and $w'$ are exponential
quasi-valuations, we have $\sqrt{I_{w}} = I_{w}$ and
$\sqrt{I_{w'}} = I_{w'}$. However, since $O_{w'}=O_{w}$, we get
$$I_{w}=\sqrt{I_{w}} = J_{w}=J_{w'}=\sqrt{I_{w'}} = I_{w'}$$ which
contradicts $y \in I_{w'} \setminus I_{w}$.

Q.E.D
\enddemo

Note that in general $O_{w'}=O_{w}$ does not imply $w'=w$; in
fact, even for valuations $v_1$ and $v_2$ on a field $F$,
$O_{v_1}=O_{v_2}$ only implies the equivalence of $v_1$ and $v_2$.
However, in Lemma 7.6 we consider exponential quasi-valuations
extending $v$ where the value groups are embedded in
$\Gamma_{\text{div}}$ (and we get $O_{w'}=O_{w}$ implies $w'=w$).
In particular, if $u_1$ and $u_2$ are valuations on $E$ extending
$v$ on $F$ where $E$ is a finite field extension of $F$ and
$O_{u_1}=O_{u_2}$, then $u_1=u_2$.

We recall from section 5 that if $u_1,...,u_n$ are valuations on a
field $E$ which extend a given valuation $v$ on $F$, then $\min
\{u_1,...,u_n\}$ is an exponential quasi-valuation on $E$
extending $v$. The corresponding quasi-valuation ring is then the
intersection of the valuation rings of the $u_i$'s. Such a
quasi-valuation ring is integrally closed.

We shall now prove that every exponential quasi-valuation $w$
extending $v$ (with $M_{w}$ weakly cancellative) is of
the form above. %exponential quasi-valuations are unique in the
%sense that every exponential quasi-valuation is of the form above.

\proclaim{Theorem 7.7} Let $(F,v)$ be a valued field and let $w$
be an exponential quasi-valuation extending $v$ on $E$. Then $$w=
\min \{u_1,...,u_k\}$$ for valuations $u_i$ on $E$ extending $v$.
\endproclaim

\demo{Proof} Let $w$ be an exponential quasi-valuation. By Lemma
7.5, Remark 5.13 and Remark 5.10, $$O_{w} =
I_{E}(O_{w})=\cap_{i=1}^k O_{u_{i}}$$ for valuation rings
$O_{u_{i}} $ where $u_{i}$ extends
$v$. %By Theorem 6.1, $w \leq u_{i}$ for every $i$.
 Denote $$w'(x)=\min \{ u_{i}(x)\}.$$ Then $w'$ is an exponential
quasi-valuation on $E$ extending $v$ and $O_{w}=O_{w'}$.
Therefore, by Lemma 7.6, $w=w'$.

%$w \leq w'$. Assume the contrary, that there exists an $x \in
%E$ such that $w(x) < w'(x)$. By Lemma 6.2 we get an element
%$y \in E$ with $w'(y)
%> w(y) =0$. Thus, $y \in I_{w'} \setminus I_{w}$. However, by
%Corollary 6.6, $\sqrt{I_{w}} = J_{w}$ (and $\sqrt{I_{w'}} =
%J_{w'}$) and since $w$ (resp. $w'$) is exponential $\sqrt{I_{w}} =
%I_{w}$ (resp. $\sqrt{I_{w'}} = I_{w'}$). However, since
%$$I_{w}=\sqrt{I_{w}} = J_{w}=J_{w'}=\sqrt{I_{w'}} = I_{w'}$$ we have a
%contradiction.

Q.E.D
\enddemo

The following Corollary summarizes the tight connection between
exponential quasi-valuations and integrally closed quasi-valuation
rings:

\proclaim{Corollary 7.8} Every exponential quasi-valuation $w$
extending $v$ induces an integrally closed (quasi-valuation) ring,
namely $O_{w}$. Conversely, every integrally closed
quasi-valuation ring is of the form $O_{w} =
I_{E}(O_{w})=\cap_{i=1}^k O_{u_{i}}$ where each $u_{i}$ is a
valuation on $E$ extending $v$, and thus $O_{w}$ has a (unique)
exponential quasi-valuation ($w=\min \{ u_{i}\}$).\endproclaim

\demo{Proof} Let $$A= \{ \text {exponential quasi-valuations
extending\ } v\}$$ and $$B= \{ \text{integrally closed
quasi-valuation rings}\}.$$ We prove that there is a 1:1
correspondence between $A$ and
 $B$. We define $f:A \rightarrow B$ by $f(w)=O_{w}$. Assuming $w_{1}\neq
w_{2}$, we use Lemma 7.6 to get $O_{w_{1}} \neq O_{w_{2}}$. Now,
if $O_{w} \in B$ then, by Remark 5.13 and Remark 5.10,
$$O_{w} = I_{E}(O_{w})=\cap_{i=1}^k O_{u_{i}};$$ we take $w=\min
\{u_1,...,u_k\} $.

 %By 5.12 every $O_{w}\in B$ is of the form $O_{w} = I_{E}(O_{w})=\cap_{i=1}^k O_{u_{i}}$
%for valuation rings $O_{u_{i}} $ where $u_{i}$ extends $v$.

Q.E.D
\enddemo

\proclaim{Corollary 7.9} There are at most $\sum_{i=1}^{n} \binom
n i=2^{n}-1$ exponential quasi-valuations
extending~ $v$.\endproclaim %\vfill \eject

\demo{Proof} By Theorem 7.7 every exponential quasi-valuation $w$
extending $v$ is of the form $w= \min \{u_1,...,u_t\}$ for
valuations $u_i$ on $E$ extending $v$ where the number of these
valuations is bounded by $[E:F]=n$.

Q.E.D.\enddemo

%\proclaim{Corollary 7.10} Let $\{u_1,...,u_t\}$ denote the set of
%all valuation on $E$ extending $v$

%There are at most $\sum_{i=1}^{n} \binom n i=2^{n}-1$ exponential
%quasi-valuations
%extending~ $v$.\endproclaim %\vfill \eject

We note that although an exponential quasi-valuation induces an
integrally closed (quasi-valuation) ring and thus is Bezout (since
it is a finite intersection of valuation rings), it is not the
case for every quasi-valuation ring. Take for example, again,
$F=\Bbb Q$, $v=v_{p}$, $M=\Bbb Z$, $e=i$ and $n = 1$ in Example
2.9.

The quasi-valuation determines the quasi-valuation ring. Some very
important questions that one may consider are: "Is there any sort
of converse"? "Is there any sort of canonical quasi-valuation
associated to a quasi-valuation ring"? "Is there a ring-theoretic
characterization of a quasi-valuation ring"? The answers to these
questions, for $M_{w}$ a monoid, are affirmative. We shall answer
these questions in the next sections.

{\it From now on we study quasi-valuations extending a valuation
on a finite field extension in the general case where the values
of the quasi-valuation lie inside a monoid and we do not assume or
deduce (as in sections 5,6 and 7) that $M_{w}$ (the value monoid)
is a group. Namely, $F$ denotes a field with a valuation $v$,
$E/F$ is a finite field extension with $n=[E:F]$, and $w:E
\rightarrow M \cup \{ \infty \}$ is a quasi-valuation on $E$ such
that $w|_{F}=v$, where $M$ is a totally ordered abelian monoid.}

%We note that in sections 5, 6 and 7, an additional restriction was
%required of $M_{w}$, such as $M_{w}$ being a group. In the next
%chapters we consider a quasi-valuation $w$ extending a valuation
%on a finite field extension with value monoid, without any
%restrictions.

We shall see in section 9 that quasi-valuations can be defined and
actually exist in any ring. Moreover, we will show that any
subring of $E$ lying over $O_{v}$ has a quasi-valuation extending
$v$, so the results we get (for quasi-valuation rings with value
monoid) apply to all such rings.

\heading \S 8 Bounding the Quasi-Valuation - the Value Monoid case
\endheading

In this section we show that any quasi-valuation extending $v$ (in
a finite field extension $E/F$) is bounded by some valuation on
$E$ extending $v$. In order to compare the values of the
quasi-valuation with values of a valuation, we construct a total
ordering on a suitable amalgamation of $M_{w}$ and
$\Gamma_{\text{div}}$.

%and in section 3 we study quasi-valuations extending a valuation
%in a finite field extension in general; i.e., the monoid we deal
%with is a totally ordered abelian monoid containing the value
%group, and no other assumptions are made (as opposed to [Sa,
%sections 3,4 and 5] in which the value monoid was actually a
%group).

\demo{Remark 8.1} Let $M$ be a totally ordered abelian monoid, let
$m \in M$, $n \in \Bbb N$ and assume $nm=0$; then $m=0$. Indeed,
if $m>0$ then $nm \geq m > 0$, and thus $nm> 0$. If $m<0$ the
argument is the same.
\enddemo

\demo{Definition 8.2} A totally ordered monoid $M$ is called {\it
$\Bbb N$-strictly ordered} if $a < b$ implies $na < nb$ for every
$n \in \Bbb N$. A totally ordered monoid $M$ is called {\it
strictly ordered} if $a < b$ implies $a+c < b+c$ for every $c \in
M$.\enddemo

Note that a totally ordered abelian monoid is strictly ordered iff
it is cancellative. Also note that strictly ordered implies $\Bbb
N$-strictly ordered. On the other hand, if $M$ is a cancellative
monoid, then totally ordered, $\Bbb N$-strictly ordered and
strictly ordered are all equivalent. Thus, we need to distinguish
among those three types of total ordering only when considering a
monoid in general.

From now on until Proposition 8.12 we consider some general
properties of a totally ordered abelian monoid $M$ containing a
group $\Gamma$ (without considering a valuation nor a
quasi-valuation).

%The following 3 remarks are valid for any totally ordered abelian
%monoid $M$ containing a group $\Gamma$.

\demo{Remark 8.3} If $m \in M$ is torsion over $\Gamma$, then $m$
is invertible. Indeed if $nm=\gamma$ then $-m=(n-1)m-\gamma$. More
generally, if $nm$ is invertible, then $m$ is invertible.\enddemo

%\demo{Remark 2} Let $x,y,z \in M$ and assume $x \geq y > z$; then
%$x>z$. Indeed otherwise $z \geq x \geq y $, a contradiction. We
%have the same result if $x \leq y < z$.\enddemo

\demo{Remark 8.4} Let $m,m' \in M$, $n \in \Bbb N$ and assume
$nm=nm'$ is invertible; then $m=m'$. Indeed by Remark 8.3, $m$ and
$m'$ are invertible and thus $n(m-m')=0$. We conclude by Remark
8.1 that $m-m'=0$; i.e., $m=m'$.\enddemo

\demo{Remark 8.5} Let $m,m' \in M$, $n \in \Bbb N$ and assume $nm
\leq \gamma \leq nm'$ with $\gamma$ invertible; then $m \leq m'$.
(We are done unless $m \geq m'$, so $nm \geq nm'$, implying
$nm=\gamma=nm' $ is invertible, so $m=m'$ by Remark 8.4).
\enddemo

\demo{Remark 8.6} Let $m,m' \in M$ and assume $m$ is invertible;
then for every $n \in \Bbb N$, $m<m' \Rightarrow nm<nm'$ and $m>m'
\Rightarrow nm>nm'$. In other words, $nm \leq nm' \Rightarrow m
\leq m'$ and $ nm \geq  nm' \Rightarrow m \geq m' $. Indeed if
$m<m'$ then $nm \leq nm'$; assume to the contrary that $nm =nm'$,
then, by Remark 8.4, we have $m=m'$, a contradiction. The proof of
the second assertion is the same.
\enddemo

Recall that $\Gamma_{\text{div}}$ is the divisible hull of
$\Gamma$. Also recall that inside $\Gamma_{\text{div}}$,
$(\gamma_{1}, n_{1})=(\gamma_{2}, n_{2})$ iff
$n_{2}\gamma_{1}=n_{1}\gamma_{2}$. We shall now show that one can
construct a total ordering on a suitable amalgamation of $M$ and
$\Gamma_{\text{div}}$ that allows us to compare elements of
$\Gamma_{\text{div}}$ with elements of $M$.

Let $M_{\text{div}}$ denote the divisible hull of $M$; namely,
$M_{\text{div}}=(M \times \Bbb N  ) / \equiv $, where $\equiv$ is
the equivalence relation defined by $(m_{1}, n_{1}) \equiv (m_{2},
n_{2})$ iff there exists a $t \in \Bbb N$ such that
$tn_{2}m_{1}=tn_{1}m_{2}$. Note that the monoid $M$ is not
necessarily $\Bbb N$-strictly ordered and thus $M$ does not
necessarily embed in $M_{\text{div}}$; i.e., the function
$\varphi_{1}: M \rightarrow M_{\text{div}}$ defined by
$\varphi_{1}(m)= (m,1)$ need not be injective. Indeed, we may have
$m_{1},m_{2} \in M$, $m_{1}<m_{2}$ and $t \in \Bbb N$ such that
$tm_{1}=tm_{2}$ and thus $(m_{1},1) = (m_{2},1)$ in
$M_{\text{div}}$. However, $\Gamma_{\text{div}}$ does embed in
$M_{\text{div}}$ via the function $\varphi_{2}:
\Gamma_{\text{div}} \rightarrow M_{\text{div}}$ defined by
$\varphi_{2}((\gamma,n))=(\gamma,n)$, since every element of
$\Gamma$ is invertible. The total ordering in $M_{\text{div}}$ is
defined by: $(m_{1}, n_{1}) \leq (m_{2}, n_{2})$ iff $tn_{2}m_{1}
\leq tn_{1}m_{2}$ for some $t \in \Bbb N$. Note that $(m_{1},
n_{1}) < (m_{2}, n_{2})$ iff for every $t \in \Bbb N$ we have
$tn_{2}m_{1}<tn_{1}m_{2}$. Also note that $M_{\text{div}}$ is a
monoid with the usual addition.

Let $T$ denote the disjoint union of $M$ and
$\Gamma_{\text{div}}$. Let $\varphi : T \rightarrow
M_{\text{div}}$ denote the function defined by:

$$ \varphi(x)= \cases \varphi_{1}(x)  \text{ \ \ \ \ if     } x \in M; \\
\varphi_{2}(x) \text{\ \ \ \ \ if     } x \in \Gamma_{\text{div}}.
\endcases$$

We define an equivalence relation $\sim$ on $T$ in the following
way:
For every $x, y \in T$,  $$ x \sim y \text{ iff } \cases x=y  \text{ \ \ \ \ \ \ \ \ \ \ \ \ for     } x,y \in M; \\
\varphi(x)=\varphi(y) \text{\ \ \ \ \ otherwise.} \endcases$$

It is easy to check that $\sim $ is indeed an equivalence relation
on $T$. We shall only check transitivity. Suppose $x \sim y$ and
$y \sim z$. The assertion is clear unless $x,z \in M$ and $y \in
\Gamma_{\text{div}}$. Write $x=m,z=m' \in M$, $y=(\gamma, n) \in
\Gamma_{\text{div}}$; then $t_{1}nm=t_{1}\gamma$ and
$t_{2}\gamma=t_{2}nm'$ for some $t_{1},t_{2} \in \Bbb N$ and thus
$$t_{1}t_{2}nm=t_{1}t_{2}\gamma=t_{1}t_{2}nm'.$$ Hence, by Remark
8.4, $m=m'$; i.e., $x \sim z$.

\demo{Definition 8.7} Let $x \in M $. $x$ is called {\it singular}
if $x \sim y$ implies $x=y$.\enddemo

\demo{Remark 8.8} If $x,y \in T$ and $x \sim y$, then $x$ is
singular iff $y$ is singular. Indeed, if $x$ is singular, we get
$x=y$.\enddemo

\demo{Remark 8.9} If $x \in T$ is not singular then there exists
$(\gamma, n) \in \Gamma_{\text{div}}$ such that $x \sim (\gamma,
n)$.
\enddemo
\demo{Proof} If $x \in \Gamma_{\text{div}}$ the assertion is
clear. If $x \in M$, then there exists $y \in T$, $y \neq x$ such
that $x \sim y$; by the definition of $\sim$, $y \in
\Gamma_{\text{div}}$.

Q.E.D.
\enddemo

We denote $\overline{T}=T/ \sim$ and define an order on
$\overline{T}$: for every $[x],[y] \in \overline{T}$,

 $$ [x] \leq [y] \text{ iff } \cases x \leq y  \text{\ \ \ \ \ \ \ \ \ \ \ \ \ \ for } x,y \text{ singular;     }  \\
\varphi(x) \leq \varphi(y) \text{\ \ \ \ \   otherwise.}
\endcases$$

\proclaim {Lemma 8.10} $\leq$ on $\overline{T}$ is well defined.
\endproclaim

\demo{Proof} Let $x,y,a,b \in T$ and assume $[x] \leq [y]$, $x
\sim a$ and $y \sim b$; we prove $[a] \leq [b]$. If $x$ and $y$
are singular then $x \leq y$, $x = a$, $y = b$, and $a$ and $b$
are singular; thus $a \leq b$ and $[a] \leq [b]$. If $x$ and $y$
are not both singular then $ \varphi(x) \leq \varphi(y)$, and $a$
and $b$ are not both singular. Since $x \sim a$ and $y \sim b$ we
have $ \varphi(x) = \varphi(a)$ and $ \varphi(y) = \varphi(b)$;
thus $ \varphi(a) \leq \varphi(b)$ and $[a] \leq [b]$.

Q.E.D.
\enddemo

\demo{Remark 8.11} For $m \in M$ and $(\gamma_{1}, n_{1}),
(\gamma_{2}, n_{2}) \in \Gamma_{\text{div}}$ we have the
following:

$\varphi(m) \leq \varphi((\gamma_{1}, n_{1})) \Rightarrow n_{1}m
\leq \gamma_{1}$;

$\varphi((\gamma_{1}, n_{1})) \leq \varphi(m) \Rightarrow
\gamma_{1} \leq n_{1}m$;

$\varphi((\gamma_{1}, n_{1})) \leq \varphi((\gamma_{2}, n_{2}))
\Rightarrow n_{2}\gamma_{1} \leq n_{1}\gamma_{2}$.

\enddemo
\demo {Proof}  We shall only prove the first assertion. By the
definition of the ordering on $M_{\text{div}}$, $\varphi(m) \leq
\varphi((\gamma_{1}, n_{1})) \Rightarrow tn_{1}m \leq t\gamma_{1}$
for some $t \in \Bbb N$; by Remark 8.6, we get $n_{1}m \leq
\gamma_{1}$.

Q.E.D.\enddemo

Note that by Remark 8.11 and by the definition of $\leq$ on
$\overline{T}$, one can exclude, for elements of $\overline{T}$,
the natural number $t$ in the definition of the ordering on
$M_{\text{div}}$. For example, assuming $[m] \leq [(\gamma, n)]$
for $m \in M$ and $(\gamma, n) \in \Gamma_{\text{div}}$. By the
definition of $\leq$ on $\overline{T}$, $\varphi(m) \leq
\varphi((\gamma, n))$, i.e., $(m,1) \leq (\gamma, n)$ in
$M_{\text{div}}$. Thus, by the definition of $\leq$ on
$M_{\text{div}}$, $tnm \leq t\gamma$ for some $t \in \Bbb N$.
However, by Remark 8.11, $\varphi(m) \leq \varphi((\gamma, n))$
implies $nm \leq \gamma$; so, $t$ is not needed. We shall use this
fact in the proof of the next proposition.

\proclaim {Proposition 8.12} $\leq$ is a total ordering on
$\overline{T}$. \endproclaim

\demo {Proof} Reflexivity and antisymmetry follow from the
reflexivity and antisymmetry of the total ordering of $M$ and
$M_{\text{div}} $. As for transitivity, assume $[x] \leq [y]$ and
$[y] \leq [z]$.

%Note that if $a \in T$ is not singular then $a \sim (\gamma, n)$
%for some $(\gamma, n) \in \Gamma_{\text{div}}$; since $leq$ is well
%defined, Assuming $[a] \leq [b]$ or $[a] \leq [b]$, we may write
%$a = (\gamma, n)$.

If $x,y,z$ are all singular, the assertion is clear.

If $x$ and $y$ are singular and $z$ is not, write
$x=m_{1},y=m_{2}$ and $z = (\gamma, n)$. (Note that by Remark 8.9,
$z \sim (\gamma, n)$ for some $(\gamma, n) \in
\Gamma_{\text{div}}$ and since $\leq$ is well defined we may write
$z = (\gamma, n)$; we shall use this fact throughout the proof).
Then $m_{1} \leq m_{2}$ and $\varphi(m_{2}) \leq \varphi((\gamma,
n))$; i.e., $nm_{2} \leq \gamma$ and thus
$$nm_{1} \leq nm_{2} \leq \gamma.$$ Hence $\varphi(m_{1}) \leq
\varphi((\gamma, n))$, namely $\varphi(x) \leq \varphi(z)$; thus
$[x] \leq [z]$.

If $x$ and $z$ are singular and $y$ is not, write
$x=m_{1},z=m_{2}$ and $y=(\gamma, n)$. Then $\varphi(m_{1}) \leq
\varphi((\gamma, n)) \leq \varphi(m_{2})$, which implies that
$nm_{1} \leq \gamma$ and  $\gamma \leq nm_{2}$. Thus
$$nm_{1} \leq \gamma \leq
nm_{2}.$$ Therefore, by Remark 8.5, $m_{1} \leq m_{2}$.

%Thus $nm_{1} \leq nm_{2}$ and therefore $m_{1} \leq m_{2}$
%(indeed, otherwise $m_{1}>m_{2} \Rightarrow nm_{1} \geq nm_{2}$
%and thus $nm_{1}  = nm_{2}= \gamma$. Now, by Remark 8.4,
%$m_{1}=m_{2}$, a contradiction).

If $y$ and $z$ are singular and $x$ is not, write
$y=m_{1},z=m_{2}$ and $x=(\gamma, n)$. Then $\varphi((\gamma, n))
\leq \varphi(m_{1})$ and $m_{1} \leq m_{2}$; i.e., $$\gamma \leq
nm_{1} \leq nm_{2}.$$ Thus $\varphi((\gamma, n)) \leq
\varphi(m_{2})$.

The rest of the cases are not difficult, we shall only prove one
of them (the others are proved in the same manner):

If $x$ is singular and $y$ and $z$ are not, write $x=m$,
$y=(\gamma_{1}, n_{1})$ and $z=(\gamma_{2}, n_{2})$. Then
$\varphi(m) \leq \varphi((\gamma_{1}, n_{1})) \leq
\varphi((\gamma_{2}, n_{2}))$; i.e., $\varphi(m) \leq
\varphi((\gamma_{2}, n_{2}))$.

%$ n_{1}m \leq \gamma_{1}$ and $  n_{2}\gamma_{1} \leq
%n_{1}\gamma_{2}$. Thus $ n_{2}n_{1}m \leq n_{2}\gamma_{1} \leq
%n_{1}\gamma_{2}$; by Remark 8.6, $ n_{2}m \leq \gamma_{2}$ and
%therefore

%$\varphi((\gamma, n)) \leq \varphi(m_{2})$.

Q.E.D. \enddemo

 Note that $\overline{T}$ is not a monoid since no operation has been
defined on it. $\overline{T}$ is merely a set that "preserves" the
values of $M$ and $\Gamma_{\text{div}}$ and allows us to compare
them.

% We recall the following result from section 5 (Proposition 5.7).

%\proclaim{Proposition 8.13} Let $R \subseteq E$ be a ring such
%that $R \cap F=O_{v}$. Let $K$ be a maximal ideal of $R$ such that
%$K \supseteq I_{v}$. Then there exists a valuation $u$ of ~$E$
%that extends~ $v$, such that $O_{u} \supseteq R$ and $I_{u}
%\supseteq K$.\endproclaim

%Also recall that the pair $(O_{u},I_{u})$ is called an {\it
%expansion} of $(R,K)$ to $E$. Our application is for $R=O_{w}$.%, noting that every maximal ideal
%$K$ of ~$O_{w}$ contains $I_{v}$, by GU.

\proclaim{Lemma 8.13} Let $w$ be a quasi-valuation extending $v$;
then there exists a valuation $u$ of $E$ extending $v$ such that
$O_{u} \supseteq O_{w}$ and $I_{u} \supseteq I_{w}$.\endproclaim

\demo {Proof} Take $R=O_{w}$ in Proposition 5.7; take a maximal
ideal $K$ of $O_{w}$ such that $K \supseteq I_{w}$. Now, expand
the pair $(O_{w},K)$.\enddemo

Note that Lemma 6.1 is a stronger statement than Lemma 8.13; in
Lemma 6.1, $O_{u} \supseteq O_{w}$ implies $I_{u} \supseteq I_{w}$
(the existence of such $O_{u}$ is shown in Corollary 5.9) whereas
in Lemma 8.13 we only have the existence of a valuation ring such
that $O_{u} \supseteq O_{w}$ and $I_{u} \supseteq I_{w}$. The
reason is that in section 6 the value monoid of the
quasi-valuation is a group.

The proof of the following theorem uses ideas which are close to
the ideas of the proofs of Lemma 6.2 and Theorem 6.3. However, one
must be very careful since we are not comparing the values inside
a group but rather inside $\overline{T}$.

\proclaim{Theorem 8.14} Let $w$ be a quasi-valuation on $E$
extending a valuation $v$ on $F$ and let $u$ be a valuation on $E$
extending $v$ such that $O_{u} \supseteq O_{w}$ and $I_{u}
\supseteq I_{w}$; then $u$ dominates $w$; i.e., $\forall x \in E$,
$w(x) \leq u(x)$ in $\overline{T}$ (see above).
\endproclaim

\demo{Proof} Let $x \in E$ and assume to the contrary that $u(x) <
w(x)$; write $\sum_{i=0}^n \alpha_{i}x^{i} = 0$ for $\alpha_{i}
\in F$. Since the sum is zero we must have $k,l \in \Bbb N$ such
that $u(\alpha_{k}x^{k}) = u(\alpha_{l}x^{l})$. Assuming $k < l$
we have $$u(x^{l-k}\alpha_{l}\alpha_{k}^{-1}) = 0$$ i.e.,
$x^{l-k}\alpha_{l}\alpha_{k}^{-1} \notin I_{u}$. Now, since $w(x)
> u(x)$, we have $w(x^{l-k}) > u(x^{l-k})$. Indeed, $$w(x^{l-k}) \geq
(l-k)w(x)$$ and $$(l-k)w(x) > (l-k)u(x) = u(x^{l-k}).$$

Therefore, by the stability of $\alpha_{l}\alpha_{k}^{-1}$ (with
respect to $u$ and $w$), we get,
$$w(x^{l-k}\alpha_{l}\alpha_{k}^{-1})
> u(x^{l-k}\alpha_{l}\alpha_{k}^{-1}) = 0.$$

%(Indeed if equality held, then reducing
%$v(\alpha_{l}\alpha_{k}^{-1})$ we would get
%$$w(x^{l-k}) = u(x^{l-k})).$$

Thus $x^{l-k}\alpha_{l}\alpha_{k}^{-1} \in I_{w}\setminus I_{u}$,
a contradiction.

Q.E.D.\enddemo

\proclaim{Theorem 8.15} Let $w$ be a quasi-valuation of the field
$E$ extending a valuation $v$. Then there exists a valuation $u$
of $E$ extending $v$ such that $O_{u} \supseteq O_{w}$, $I_{u}
\supseteq I_{w}$ and $u$ dominates $w$.
\endproclaim

\demo{Proof} The existence of a valuation $u$ with $O_{u}
\supseteq O_{w}$ and $I_{u} \supseteq I_{w}$ is by Lemma 8.13.
Now, apply Theorem 8.14.

 Q.E.D.\enddemo

In section 10 (Theorem 10.5) we prove a stronger version of
Theorem 8.15 for the filter quasi-valuation. We show there that
for every valuation $u$ of $E$ extending $v$ such that $O_{u}
\supseteq O_{w}$, $u$ dominates the filter quasi-valuation.

\heading \S 9 Filter Quasi-Valuations - Extending Valuations to
Quasi-Valuations
\endheading

Our goal in this section is to construct a quasi-valuation
extending a given valuation. First we obtain a value monoid (we
call it the cut monoid), constructed from the value group. We then
show that one can extend a valuation to a quasi-valuation with
values inside this particular monoid.

 \heading Cuts of ordered sets \endheading

We start this section by reviewing some of the notions of cuts. We
shall review the parts needed to construct the cut monoid; for
more information about cuts see, for example, [FKK] or [Weh].

In this subsection $T$ denotes a totally ordered set.

\demo{Definition 9.1} A subset $S$ of $T$ is called initial (resp.
final) if for every $\gamma \in S$ and $\alpha \in T$, if $\alpha
\leq \gamma$ (resp. $\alpha \geq \gamma$), then $\alpha \in
S$.\enddemo

\demo{Definition 9.2} A cut $\Cal A=(\Cal A^{L}, \Cal A^{R})$ of
$T$ is a partition of $T$ into two subsets $\Cal A^{L}$ and $\Cal
A^{R}$, such that, for every $\alpha \in \Cal A^{L}$ and $\beta
\in \Cal A^{R}$, $\alpha<\beta$.\enddemo

 Note that the set of all cuts $\Cal A=(\Cal A^{L}, \Cal A^{R})$
 of the ordered set $T$ contains the two cuts $(\emptyset,T)$ and
 $(T,\emptyset)$; these are commonly denoted by $-\infty$ and
 $\infty$, respectively. However, we shall not use the symbols $-\infty$ and
 $\infty$ to denote the above cuts since we
 shall define a "different" $\infty$.

Given $\alpha \in T$, we denote
$$(-\infty,\alpha]=\{\gamma \in T \mid \gamma \leq \alpha \}$$ and
$$(\alpha,\infty)=\{\gamma \in T \mid \gamma > \alpha \}.$$

One defines similarly the sets $(-\infty,\alpha)$ and
$[\alpha,\infty)$.

To define a cut we shall often write $\Cal A^{L}=S$, meaning the
$\Cal A$ is defined as $(S, T \setminus S)$ when $S$ is an initial
subset of $T$.

\demo{Definition 9.3} The ordering on the set of all cuts of $T$
is defined by $\Cal A \leq \Cal B$ iff $\Cal A^{L} \subseteq \Cal
B^{L}$ (or equivalently $\Cal A^{R} \supseteq \Cal
B^{R}$).\enddemo

Given $S  \subseteq T$, $S^{+}$ is the smallest cut $\Cal A$ such
that $S  \subseteq \Cal A^{L}$. So, for $\alpha \in T$ we have
$\{\alpha\}^{+}=((-\infty,\alpha],(\alpha,\infty))$. We denote
$\{\alpha\}^{+}$ by $\alpha^{+}$.

\heading The cut monoid \endheading

In this subsection, $\Gamma$ will denote a totally ordered abelian
group and $\Cal M(\Gamma)$ will denote the set of all cuts of
$\Gamma$.

\demo{Definition 9.4} Let $S,S' \subseteq \Gamma$ and $n \in \Bbb
N$, we define

$$S+S'=\{ \alpha+\beta \mid \alpha \in S, \beta \in S' \};$$

$$nS=\{ s_{1}+s_{2}+...+s_{n} \mid s_{1},s_{2},...,s_{n} \in S
\}.$$

\enddemo
%$\gamma +S= \{ \gamma+\alpha \mid \alpha \in S \}$ (Do I need it?)

\demo{Definition 9.5} For $\Cal A , \Cal B \in \Cal M(\Gamma)$,
their (left) sum is the cut defined by $$(\Cal A + \Cal
B)^{L}=\Cal A^{L} + \Cal B^{L}.$$\enddemo

One can also define the right sum; however, we shall not use it.
Note that under the above definitions, the zero in $\Cal
M(\Gamma)$ is the cut $0^{+}=((-\infty,0],(0,\infty))$.

\demo{Definition 9.6} For $\Cal A \in \Cal M(\Gamma)$ and $n \in
\Bbb N$, we define the cut $n\Cal A$ by $$(n\Cal A)^{L}=n\Cal
A^{L}.$$\enddemo

The following lemma is well known (see for example [FKK] or
[Weh]).

\proclaim{Lemma 9.7} $(\Cal M(\Gamma),+,\leq)$ is a totally
ordered abelian monoid.

\endproclaim

\demo{Definition 9.8} We call $\Cal M(\Gamma)$ the \it {cut
monoid} of $\Gamma$. \enddemo

\demo{Remark 9.9} Note that there is a natural monomorphism of
monoids $\varphi : \Gamma \rightarrow \Cal M(\Gamma)$ defined in
the following way: for every $\alpha \in \Gamma$, $$\varphi
(\alpha) =\alpha^{+} $$ %((-\infty,\alpha],(\alpha,\infty)) $$ since
since $\alpha<\beta \in \Gamma$ implies $(-\infty,\alpha] \subset
(-\infty,\beta]$. Therefore, when we write $N \cap \Gamma$ for a
subset $N \subseteq \Cal M(\Gamma)$, we refer to the intersection
of $N$ with the copy $\{\alpha^{+} \mid {\alpha \in \Gamma}\}$ of
$\Gamma$, which could also be viewed as the set $\{ \alpha \in
\Gamma \mid \alpha^{+} \in N \}$. To simplify notation, for
$\alpha \in \Gamma$, when viewing $\alpha$ inside $ \Cal
M(\Gamma)$, we shall write $\alpha$ instead of $\alpha^{+}$. For
example, for $\Cal B \in \Cal M(\Gamma)$, $\Cal B+\alpha$ is the
cut defined by $(\Cal
B+\alpha)^{L} =  \Cal B^{L}+(-\infty,\alpha]$.% and $\Cal
%B+(-\alpha)$ is the cut defined by $(\Cal B+(-\alpha))^{L} = \Cal
%B^{L}+(-\infty,-\alpha]$. Finally, note that $\alpha \in \Gamma$
%is invertible in $\Cal M(\Gamma)$.

%Moreover, note that if $\Cal B,\Cal C \in \Cal M(\Gamma)$ and
%$\alpha \in \Gamma$, then $(-\infty,\alpha]+ \Cal
%B^{L}=(-\infty,\alpha]+\Cal C^{L}$ implies $\Cal B^{L}=\Cal C^{L}$
%(Indeed $c \in \Cal C^{L} \Rightarrow \alpha+c \in
%(-\infty,\alpha]+\Cal B^{L} \Rightarrow (\alpha - \gamma)+c \in
%\Cal B^{L}$ for some $\gamma \in (-\infty,\alpha]$, i.e., $\alpha
%- \gamma>0$; thus $c \in \Cal B^{L}$); in other words,
%$\alpha+\Cal B=\alpha+\Cal C$ implies $\Cal B=\Cal C$.

\enddemo

\demo{Definition 9.10} Let $\alpha \in \Gamma$ and $\Cal B \in
\Cal M(\Gamma)$. We write $\Cal B-\alpha$ for the cut $\Cal
B+(-\alpha)$ (viewing $-\alpha$ as an element of $\Cal
M(\Gamma)$).

\enddemo

\demo{Remark 9.11} $\Cal B+\alpha$ is actually the cut defined by
$(\Cal B+\alpha)^{L}=\{ \beta+\alpha \mid \beta \in \Cal B^{L}
\}$.\enddemo

\demo{Proof} Obviously $\{ \beta+\alpha \mid \beta \in \Cal B^{L}
\} \subseteq (\Cal B+\alpha)^{L}$. For the converse, if $\gamma
\in (\Cal B+\alpha)^{L}$, then for some $\varepsilon \in \Cal
B^{L}$ and $\delta \in (-\infty,\alpha]$,
$\gamma=\varepsilon+\delta \leq \varepsilon+\alpha$. So,
$\varepsilon \geq \gamma-\alpha$ and thus $\gamma-\alpha \in \Cal
B^{L}$, whence, $\gamma=(\gamma-\alpha)+\alpha \in \{ \beta+\alpha
\mid \beta \in
\Cal B^{L}\}$. %On the other hand, let
%$\gamma=\beta-\alpha$ for $\beta \in \Cal B^{L}$; then obviously
%$\gamma \in \Cal B^{L}+(-\infty,-\alpha]$.

\enddemo

%\demo{Remark 9.11} $\Cal B-\alpha$ is actually the cut defined by
%$(\Cal B-\alpha)^{L}=\{ \beta-\alpha \mid \beta \in \Cal B^{L}
%\}$.\enddemo

%\demo{Proof} Let $\gamma \in (\Cal B-\alpha)^{L}$ then
%$\gamma=\beta+\delta$ for $\beta \in \Cal B^{L}$ and $\delta \leq
%-\alpha$. So, $\beta \geq \beta+\delta+\alpha$ and thus
%$\beta+\delta+\alpha \in \Cal B^{L}$; hence
%$\gamma=\beta+\delta+\alpha-\alpha$. On the other hand, let
%$\gamma=\beta-\alpha$ for $\beta \in \Cal B^{L}$; then obviously
%$\gamma \in \Cal B^{L}+(-\infty,-\alpha]$.

%\enddemo

%\demo{Remark} Let $\alpha \in \Gamma$, and $\Cal B, \Cal C \in
%\Cal M(\Gamma)$. Viewing $\alpha$ as an element of $\Cal
%M(\Gamma)$, we obviously have $\alpha+(-\alpha)=-\alpha+\alpha=0$.
%So, $\alpha$ is invertible in $\Cal M(\Gamma)$. Thus, for example,
%$\Cal B+\alpha=\Cal C$ iff $\Cal B=\Cal C-\alpha$; $\alpha+\Cal
%B=\alpha+\Cal C$ implies $\Cal B=\Cal C$.
%\enddemo

\proclaim{Lemma 9.12} The cut monoid  $\Cal M(\Gamma)$ is $\Bbb
N$-strictly ordered.
\endproclaim

\demo{Proof} Let $\Cal A,\Cal B \in \Cal M(\Gamma)$ and assume
$\Cal A<\Cal B$, i.e., $\Cal A^{L} \subset \Cal B^{L}$. Then there
exists $\beta \in \Cal B^{L} \setminus \Cal A^{L}$; namely
$\beta>\alpha$ for all $\alpha \in \Cal A^{L}$. Thus $n\beta
> \alpha_{1}+\alpha_{2}+...+\alpha_{n}$ for all $\alpha_{1},\alpha_{2},...,\alpha_{n} \in \Cal
A^{L}$; therefore $n\beta \in n \Cal B^{L} \setminus n \Cal
A^{L}$, i.e., $n\Cal A<n\Cal B$.

Q.E.D.\enddemo

%\demo{Example 9.13} Let $\Gamma$ be a totally ordered abelian
%group, then $(\emptyset,\Gamma),(\Gamma,\emptyset) \in \Cal
%M(\Gamma) $; ignoring those two cuts, we have:

% If $\Gamma= \Bbb Z$ then every cut of $\Cal M(\Bbb Z)$ is of the
% form $((-\infty,\alpha],(\alpha,\infty))$, for $\alpha \in \Bbb Z$.

% If $\Gamma= \Bbb Q$ then there are three kinds of cuts in $\Cal M(\Bbb Q)$:

%\ \ 1. Rational cuts of the form
%$((-\infty,\alpha],(\alpha,\infty))$, for $\alpha \in \Bbb Q$;

% \ \ 2. Rational cuts of the form $((-\infty,\alpha),[\alpha,\infty))$, for $\alpha \in \Bbb Q$;

%\ \  3. Irrational cuts of the form $(\{ \alpha \in \Bbb Q \mid
%\alpha<r\}, \{ \alpha \in \Bbb Q \mid \alpha>r\})$
 % for $r \in \Bbb R \setminus \Bbb Q$.

 % If $\Gamma= \Bbb R$ then there are two kinds of cuts in $\Cal M(\Bbb R)$:

%\ \ 1. Cuts of the form $((-\infty,\alpha],(\alpha,\infty))$, for
%$\alpha \in \Bbb R$;

% \ \ 2. Cuts of the form $((-\infty,\alpha),[\alpha,\infty))$, for $\alpha \in \Bbb
% R$.
%\enddemo

 \heading Constructing the filter quasi-valuation \endheading

\demo{Definition 9.13} Let $v$ be a valuation on a field $F$ with
value group $\Gamma_{v}$. Let $O_{v}$ be the valuation ring of $
v$ and let $R$ be an algebra over $O_{v}$. For every $x \in R$,
the $O_{v}$-{\it support} of $x$ in $R$ is the set
$$S^{R/O_{v}}_{x}=\{ a \in O_{v} | xR \subseteq aR\}.$$ We suppress
$R/O_{v}$ when it is understood. Note that $xR \subseteq aR$ iff
$x \in aR$.
\enddemo

\demo{Definition 9.14} Let $A$ be a collection of sets and let $B$
be a subset of $A$; we call $B$ a {\it filter} of $A$ if $B$
satisfies the following:

1. If $C \in B$ and $D$ is an element of $A$ containing $C$, then
$D \in B$;

2. If $C,D \in B$ then $C \cap D \in B$;

3. $\emptyset \notin B$.

\enddemo

\demo{Remark 9.15} For every $x \in R$, the set $\{ aO_{v} \mid a
\in S_{x}^{R/O_{v}}\}$ is a filter of $\{ aO_{v} \mid a \in
O_{v}\}$.
\enddemo

For every $A \subseteq F$ we denote $v(A)=\{ v(a) \mid a \in A
\}$. Note that for $A \subseteq O_{v}$, we have $v(A)=(v(A))^{\geq
0}$ (viewing $v(A)$ as a subset of $\Gamma_{v} \cup \{ \infty
\}$). In particular,

$$v(S_{x})=(v(S_{x}))^{\geq 0}=\{ v(a) | a \in S_{x} \};$$ the reason for
this notion is the following observation.

 \proclaim{Lemma 9.16} Let $0 \neq x \in R$. Then $v(S_{x})$ is an initial subset of $\ (\Gamma_{v})
^{\geq 0}$.
\endproclaim

\demo{Proof} First note that $1 \in S_{x}$; so $S_{x} \neq
\emptyset$. Now, let $\alpha \in v(S_{x})$; then there exists an
$a \in O_{v}$ with $v(a)= \alpha$ such that $x \in aR$. Let $0
\leq \beta \leq \alpha$; then there exists a $b \in O_{v}$ with
$v(b) = \beta$. Thus, $a \in bO_{v}$ and $x \in aR \subseteq bR$.
Therefore $b \in S_{x}$ and $\beta \in v(S_{x})$.

Q.E.D.\enddemo

%To be able to work with cuts, we denote, for every $A \subseteq
%O_{v}$, $v(A)=(v(A))^{\geq 0} \cup (\Gamma_{v})^{<0}$; in
%particular,
%$$v(S_{x})=(v(S_{x}))^{\geq 0} \cup (\Gamma_{v})^{<0}.$$ With this new notation and by Lemma 9.16, $v(S_{x})$
%is an initial subset of $\Gamma_{v}$.

Note that by Lemma 9.16, we have, for $0 \neq x \in R$,
$$(v(S_{x})^{+})^{L}=v(S_{x}) \cup (\Gamma_{v})^{<0}.$$ Also
note that if $A$ and $ B$ are subsets of $O_{v}$ such that $A
\subseteq B$ then $v(A) \subseteq v(B)$.

Recall that we do not denote the cut $(\Gamma_{v},\emptyset) \in
\Cal M(\Gamma_{v}) $ as $\infty$. So, as usual, we adjoin to $\Cal
M(\Gamma_{v}) $ an element $\infty$ greater than all elements of
$\Cal M(\Gamma_{v})$; for every $\Cal A \in \Cal M(\Gamma_{v})$
and $\alpha \in \Gamma_{v}$ we define $\infty+\Cal A=\Cal
A+\infty=\infty$ and $\infty-\alpha=\infty$.

 \demo{Definition 9.17} Let $v$ be a valuation on a field $F$
with value group $\Gamma_{v}$. Let $O_{v}$ be the valuation ring
of $v$ and let $R$ be an algebra over $O_{v}$. Let $\Cal
M(\Gamma_{v})$ denote the cut monoid of $\Gamma_{v}$. We say that
a quasi-valuation $w:R \rightarrow \Cal M(\Gamma_{v}) \cup \{
\infty \}$ is {\it induced by} $(R,v)$ if $w$ satisfies the
following:

 1. $w(x)=v(S_{x})^{+}$ for every $0 \neq x \in
 R$. I.e.,  $w(x)^{L}=v(S_{x}) \cup (\Gamma_{v})^{<0}$;

 2. $w(0)=\infty$.

 %3. $w(x) \geq 0$ for every $x \in R$.
\enddemo

\demo{Remark 9.18} Notation as in Definition 9.17 and let $0 \neq
r \in R$. We note that it is possible to have $v(S_{r})=\
(\Gamma_{v}) ^{\geq 0}$ and thus $w(r)=(\Gamma_{v},\emptyset)$;
for example, take $R=F$. However, $ 1 \in S_{r} $ and thus $0 \in
v(S_{r})$; therefore $(v(S_{r})^{+})^{L} \supseteq (-\infty,0] $;
i.e., $w(r) \geq 0$. Thus, for $w$ a quasi-valuation induced by
$(R,v)$ and $r \in R$, we have $w(r) \geq 0$ (recall that by
definition $w(0)=\infty$). In particular, $w$ cannot satisfy
$w(r)=(\emptyset, \Gamma_{v})$, i.e., $(\emptyset, \Gamma_{v})
\notin im(w)$.

%note that it is possible to have $v(S_{x})=\Gamma_{v}$; however,
%by definition $v(S_{x}) \supseteq (\Gamma_{v})^{<0} $. Thus, a
%quasi-valuation induced by $(R,v)$ cannot satisfy
%$w(x)^{L}=\emptyset$, i.e., $(\emptyset, \Gamma_{v}) \notin
%im(w)$.

%For $R$ an algebra over $O_{v}$ and $0 \neq r \in R$, note that it
%is possible to have $v(S_{r})=\Gamma_{v}$ and thus
%$w(r)=(\Gamma_{v},\emptyset)$. However, by definition, $v(S_{r})
%\supseteq \{ \alpha \in \Gamma_{v} \mid \alpha < 0\} $ and $0 \in
%v(S_{r})$; therefore $v(S_{r}) \supseteq \{ \alpha \in \Gamma_{v}
%\mid \alpha \leq 0\} $. Thus, for every $r \otimes \frac{1}{\beta}
%\in R \otimes_{O_{v}}F$, $$(W(r \otimes
%\frac{1}{\beta}))^{L}=v(S_{r})+(-\infty,-v(\beta)] \supseteq
%(-\infty,-v(\beta)],$$ i.e., $(W(r \otimes \frac{1}{\beta}))^{L}
%\neq \emptyset$; hence $(\emptyset, \Gamma_{v}) \notin im(W)$.

\enddemo

 The following theorem holds for arbitrary algebras $R$
(not necessarily integral domains).

\proclaim{Theorem 9.19} Let $v$ be a valuation on a field $F$ with
value group $\Gamma_{v}$. Let $O_{v}$ be the valuation ring of $\
v$ and let $R$ be an algebra over $O_{v}$. Let $\Cal
M(\Gamma_{v})$ denote the cut monoid of $\Gamma_{v}$. Then there
exists a quasi-valuation $w:R \rightarrow \Cal M(\Gamma_{v}) \cup
\{ \infty \}$ induced by $(R,v)$.
\endproclaim

\demo{Proof} We define for every $0 \neq x \in R$ a function $w: R
\rightarrow \Cal M(\Gamma_{v}) \cup \{ \infty \}$ by
$$w(x)=v(S_{x})^{+}; $$
and $w(0)= \infty$. Note that by definition $w $ satisfies
conditions 1 and 2 of Definition 9.17. We prove that $w$
is indeed a quasi-valuation.% and satisfies condition 3 of
%Definition 9.17.

%We prove that $w$ is indeed a quasi-valuation on $R$.

First note that if $x,y \in R$ such that at least one of them is
zero then it is easily seen that $ w(xy) \geq w(x)+w(y)$ and $
w(x+y) \geq \min \{ w(x),w(y) \}$; so, we may assume that $x,y \in
R$ are both non-zero.

Let $x,y \in R$ and note that $S_{xy} \supseteq S_{x} \cdot S_{y}$
where $S_{x} \cdot S_{y} = \{ a \cdot b | a \in S_{x}, b \in S_{y}
\}$. Indeed, let $a \in S_{x}$, $b \in S_{y}$; then $x=ar$, $y=bs$
for some $r,s \in R$ and thus $xy=(ab)(rs)$. Therefore, $$v(
S_{xy}) \supseteq v(S_{x} \cdot S_{y})=v(S_{x})+v(S_{y}).$$
%$$xy= (ar)(bs)=a(r(bs))=a(b(rs))=(ab)(rs).$$

Hence,
$$w(xy) = v( S_{xy})^{+} \geq (v(S_{x})+v(S_{y}))^{+}$$$$=v(S_{x})^{+}+v(S_{y})^{+}
 =w(x)+w(y);$$ i.e.,
$$ w(xy) \geq w(x)+w(y).$$

Now, assume $w(x) \leq w(y)$ i.e., $v(S_{x}) \subseteq v(S_{y})$.
Thus $S_{x} \subseteq S_{y}$. (Indeed, assume to the contrary that
there exists an $a \in S_{x} \setminus S_{y}$; then for every $a'
\in O_{v}$ with $v(a')=v(a)=\alpha$, one has $a' \notin S_{y}$ and
thus $\alpha \notin v(S_{y})$, a contradiction). Consequently, if
$a \in S_{x}$ we have
$$(x+y)R \subseteq xR+yR \subseteq aR+aR = aR ;$$ i.e., $a \in
S_{x+y}$. Therefore $$v( S_{x+y})  \supseteq v(S_{x});$$ i.e., $v(
S_{x+y})^{+} \geq   v(S_{x})^{+}.$ Thus
$$ w(x+y) \geq w(x)= \min \{ w(x),w(y) \}.$$

%Finally, note that if $x \in R$ then $w(x) \geq
%0=((-\infty,0],(0,\infty))$. Indeed, for $x \neq 0$, we have $ 1
%\in S_{x} \subseteq O_{v}$ and thus $0 \in (v(S_{x}))^{\geq 0} $
%i.e., $v(S_{x}) \supseteq (-\infty,0]$; if $x=0$ then, by
%definition, $w(x)= \infty$.

%Finally, let $x \in R$ and $0 \neq c \in O_{v}$; then $$a \in
%S_{x} \Rightarrow x \in aR \Rightarrow xc \in ac R \Rightarrow ac
%\in S_{xc}.$$ On the other hand, $$b \in S_{xc} \Rightarrow xc \in
%bR \Rightarrow x \in bc^{-1} R. $$ Now, If $bc^{-1} \notin O_{v}$
%then $bc^{-1} \notin S_{x}$; otherwise $bc^{-1} \in O_{v}$ and
%thus $bc^{-1} \in S_{x}$. Consequently, $v(b) \geq v(c)$, i.e.,
%$b=ac$ for some $a \in O_{v}$. Also note that since $a=bc^{-1}$,
%we have $a \in S_{x}$. Thus $w(xc)=v(S_{xc})=\{ v(ac) \mid ac \in
%S_{xc} \}= \{ v(a)+v(c) \mid a \in S_{x} \}=\{ v(a) \mid a \in
%S_{x} \}+v(c)=w(x)+v(c)$. If $c=0$ the assertion is clear.

Q.E.D.\enddemo

\demo {Definition 9.20} In view of Remark 9.15, the
quasi-valuation constructed in Theorem 9.19 is called the {\it
filter quasi-valuation} induced by $(R,v)$.
\enddemo

%\demo {Remark 4.14} Note that by definition $w(r)=\infty$ iff
%$r=0$.
%\enddemo

\demo {Example 9.21} Let $v$ be a valuation on a field $F$. Let
$O_{v}$ denote the valuation ring and let $I \vartriangleleft
O_{v}$. Then $O_{v} / I $ is an algebra over $O_{v}$ whose filter
quasi-valuation induced by $(O_{v} / I,v)$
is:    $$ w(x+I)=\cases v(x)  \text{ \ \  if     } x \notin I \\
\infty \text{\ \ \ \ \ otherwise.} \endcases$$ Here $w: O_{v} / I
\rightarrow \Cal M(\Gamma_{v}) \cup \{ \infty \} $. Note that in
the definition of $w$ we write $v(x)$ as an element of $\Cal
M(\Gamma_{v})$, i.e., $v(x)=((-\infty,v(x)],(v(x),\infty))$. Also
note that $w$ is well defined. We show now that $w$ is indeed the
filter quasi-valuation induced by $(O_{v} / I,v)$. Let $w_{v}$
denote the filter quasi-valuation induced by $(O_{v} / I,v)$; if
$x \in I$ then $x \equiv 0$ and $w_{v}(x+I)=w_{v}(0)=\infty$. If
$x \notin I$ then obviously $w_{v}(x+I) \geq v(x)$; assume that
$w_{v}(x+I)
> v(x)$. Then there exists $y \in O_{v}$ with $$w_{v}(x+I)\geq v(y) >
v(x)$$ such that $x+I \in y O_{v}/ I$; i.e., $x-ya \in I$ for some
$a \in O_{v}$. Note that $v(x)= v(x-ya)$ and thus $x \in I$, a
contradiction. So, $w_{v}=w$. %i.e., $w$ as defined above is indeed
%the filter quasi-valuation induced by $(O_{v} / I,v)$.

%$$ x \equiv y \text{ \ \  iff     } \cases x=y  \text{ \ \  if     } x,y \in M \\
%\varphi(x)\equiv \varphi(y)   \text{\ \ \ \ \
%otherwise.}\endcases$$

%$w$ is called the {\it lopped quasi-valuation}.
\enddemo

Note that if $I=P$ is a prime ideal of $O_{v}$, then $w$ above is
actually a valuation on $O_{v} / P $.

%Indeed $O_{v}/P$ is a valuation ring and one can define a
%corresponding valuation by:
%$$\overline{v}(x+P)=v(x) \text{\ if\ } x \notin P \text{\ and\ }
%\overline{v}(x+P)=\infty \text{\ otherwise}.$$

\demo {Remark 9.22} A different way to define the filter
quasi-valuation described in Example 9.21 is by using the notion
of final subsets. %Recall that a subset $B$ of an ordered set $A$
%is called major if $$\{ y \in A \mid x \leq y \} \subseteq B$$ for
%every $x \in B$ (see [Bo, p. 391]).
Assuming $v$ is a valuation on a field $F$ with valuation ring
$O_{v}$ and value group $\Gamma_{v}$, it is known that there is a
one to one correspondence between the set of final subsets of
$\Gamma_{v}$ and the set of $O_{v}$-submodules of $F$. (See [Bo,
p. 391-392]; final subsets are called there major subsets). Thus,
one can equivalently define

$$ w(x+I)=\cases v(x)  \text{ \ \  if     } v(x) \notin B \\
\infty \text{\ \ \ \ \ otherwise.} \endcases$$

where $B$ denotes the final subset corresponding to $I$ in Example
9.21. \enddemo

\demo {Remark 9.23} Notation as in Theorem 9.19. For every $c \in
O_{v}$, we have $$v(S_{c \cdot 1_{R}}) \supseteq [0,v(c)],$$ i.e.,
$w(c \cdot 1_{R}) \geq v(c)$. Indeed $c \cdot 1_{R} \in cR$ and
thus $c \in S_{c \cdot 1_{R}}$. Note that one can easily have $w(c
\cdot 1_{R}) > v(c)$. For example, take $0 \neq c \in O_{v}$ such
that $c \cdot 1_{R}=0$ (in particular, $R$ is not faithful); then
$\infty = w(0)=w(c \cdot 1_{R})
> v(c)$. In Example 9.28 we shall see a situation in which $\infty
>w(c)>v(c)$.
\enddemo

\demo {Remark 9.24} Let $R$ be a torsion free algebra over an
integral domain ~$C$. Let $0 \neq c \in C$, $b \in C$ satisfying
$c^{-1}b \in C$; let $x,r \in R$ and assume $cx=br$. Since
$$br=c(c^{-1}b)r,$$ we may cancel $c$ and conclude that
$x=(c^{-1}b)r$.
%Then $x=(c^{-1}b)r$. Indeed, write $y=(c^{-1}b)r \in (c^{-1}b)R
%\subseteq R$; then $cy=br=cx$ i.e., $c(y-x)=0$ which implies that
%$y=x$ since $R$ is torsion free over $C$.
Note that if $R$ is not torsion free, this fact is not true, as we
shall see in Corollary 9.26.
%indeed if $0 \neq c \in C$ and $0 \neq x \in R$ such that $cx=0$,
%then writing $cx=c \cdot 0$, we have $x \neq (c^{-1}c) \cdot 0=0$.

\enddemo

 \proclaim{Lemma 9.25} Notation as in Theorem
9.19, assume in addition that $R$ is torsion free over $O_{v}$;
then
%$v(S_{cx})=[0,v(c)]+v(S_{x})$. More generally,
$$w(cx)=v(c)+w(x)$$ for every $c \in O_{v}$, $x \in R$.
\endproclaim

\demo{Proof} First note that if $c=0$ or $x=0$ then
$w(cx)=v(c)+w(x)$ is clear. Now, by Remark 9.23 and by the fact
that $w$ is a quasi-valuation, we have
$$w(cx)=w(c \cdot 1_{R}\cdot x) \geq w(c \cdot 1_{R})+w(x) \geq
v(c)+w(x).$$

For the other direction we show that if $v(b) \in v(S_{cx})$ for
$b \in O_{v}$, then $$v(b) \in [0,v(c)]+v(S_{x}).$$ Note that if
$v(b) < v(c)$ then clearly $v(b) \in [0,v(c)]+v(S_{x})$. (Indeed,
$0 \in v(S_{x})$ and $[0,v(c)]$ is an initial subset of
$(\Gamma_{v})^{\geq 0}$ and thus $v(b) \in [0,v(c)] $). Thus, we
may assume that $v(b) \geq v(c)$, i.e., $c^{-1}b \in O_{v}$.
Therefore, by the definition of $S_{cx}$ and Remark 9.24, we have
$$b \in S_{cx} \Rightarrow cx \in bR \Rightarrow x \in c^{-1}b
R.$$ So we have $c^{-1}b \in S_{x}$, and writing $b= c(c^{-1}b) $,
we conclude that $$v(b)=v(c)+v(c^{-1}b) \in [0,v(c)]+v(S_{x}).$$

% Now, if $bc^{-1} \notin O_{v}$ i.e., $bc^{-1} \notin S_{x}$,
%then $v(b)<v(c)$ and thus $v(b)=0+v(b) \in v(S_{x})+[0,v(c)]$
%(Indeed, $0 \in S_{x}$ and $[0,v(c)]$ is a positive isolated
%subset). If $bc^{-1} \in O_{v}$ i.e., $bc^{-1} \in S_{x}$, then
%writing $b= (bc^{-1}) c$, we have $v(b)=v(bc^{-1})+ v(c) \in
%v(S_{x})+[0,v(c)]$. On the other hand, let $v(a)\in v(S_{x})$ and
%$v(t) \in [0,v(c)]$, then $v(t) \leq v(c)$ and $$a \in S_{x}
%\Rightarrow x \in aR \Rightarrow xc \in ac R \Rightarrow ac \in
%S_{xc}.$$ Thus $v(ac) \in v(S_{xc})$ and obviously $v(ac) \geq
%v(a)+v(t)$. Therefore, by Lemma 4.10, $v(a)+v(t) \in v(S_{xc})$.

%$$a \in S_{x} \Rightarrow x \in aR \Rightarrow
%xc \in ac R \Rightarrow ac \in S_{xc}.$$ On the other hand, $$b
%\in S_{xc} \Rightarrow xc \in bR \Rightarrow x \in bc^{-1} R. $$
%Now, If $bc^{-1} \notin O_{v}$ then $bc^{-1} \notin S_{x}$;
%otherwise $bc^{-1} \in O_{v}$ and thus $bc^{-1} \in
%S_{x}$.

Q.E.D.\enddemo

We deduce the following corollary:

\proclaim{Corollary 9.26} Let $R$ be an algebra over $O_{v}$ and
let $w$ denote the filter quasi-valuation induced by $(R,v)$; then
$R$ is torsion free iff $$w(cx)=v(c)+w(x)$$ for every $c \in
O_{v}$ and $x \in R$.\endproclaim

\demo{Proof} $(\Rightarrow)$ By Lemma 9.25. $(\Leftarrow)$ If $R$
is not torsion free over $O_{v}$, then there exist $0 \neq c \in
O_{v}$ and $0 \neq x \in R$ such that $cx=0$. Hence $$\infty =
w(0)=w(c x)
> v(c)+w(x)$$ (since $v(c), w(x) < \infty$).

Q.E.D.\enddemo

\demo {Remark 9.27} Note that even in the case where $R$ is a
torsion free algebra over $O_{v}$, one does not necessarily have
$w(c \cdot 1_{R})=v(c)$ for $c \in O_{v}$, despite the fact that
$$w(c \cdot 1_{R}) = v(c)+w( 1_{R})$$ by Lemma 9.25. The reason is that $w(
1_{R})$ is not necessarily $0$, as shown in the following example.

\enddemo

\demo {Example 9.28} Let $v$ be a valuation on a field $F$,
$O_{v}$ its valuation ring and let $ O_{u} \subseteq F$ be any
ring strictly containing $O_{v}$; then $O_{u}$ is a valuation ring
of $F$ and $O_{u}= S^{-1}O_{v}$ where $S=O_{v} \setminus P$ for
some non-maximal prime ideal $P \vartriangleleft O_{v}$; see [Bo,
section 4] for further discussion. Taking $R=O_{u}$ in Theorem
9.19 and noting that $O_{u}$ is obviously torsion free over
$O_{v}$, we see that for every $c \in O_{v}$, one can write $c=cs
\cdot s^{-1} \in csR$ for every $s \in S$. Also, if $t \notin S$
then $c \notin ctR$ (since otherwise $t^{-1} \in R$ which is
impossible). Thus $$v(S_{c})=\{ \alpha \in (\Gamma_{v})^{\geq 0}
\mid \alpha \leq v(cs) \text{ for some } s \in S\}.$$ Writing $H$
for the isolated subgroup corresponding to $P$, we deduce that
$$w(c)^{L}= (-\infty,v(c)]+H^{\geq 0}.$$ %where $H^{\geq 0}=\{ h \in H \mid h \geq 0\}$.
In particular, $w(1)^{L}=(-\infty,0]+H^{\geq 0}$; thus $w(1)>0$
(note that $P \neq I_{v}$ and thus $H \neq \{  0\}$).

We note that $(-\infty,0]+H^{\geq 0}=(-\infty,0] \cup H^{\geq 0}$.

%Note that $[0,v(c)]+H$ refers to sums in $\Cal M$ and not to
%modulo $H$.

\enddemo

Note that in Example 9.28 $w(1) \in \overline{hull_{\Cal
M(\Gamma_{v})}(H^{\geq0})}$ yet $w(1) > \alpha$ for every $\alpha
\in H^{\geq0}$ (though $w$ does not extend $v$); see also Lemma
10.6 for further discussion on such elements.

%Thus $S_{c}=\{ cs \mid s \in S\}$.

The next observation is well known.

 \demo{Remark 9.29} Let $C$ be an integral domain, $S$ a
multiplicative closed subset of $C$ with $0 \notin S$, and $R$ an
algebra over $C$. We claim that every $x \in R \otimes_{C}CS^{-1}$
is of the form $r \otimes \frac{1}{\beta}$ for $r \in R$ and
$\beta \in S$. Indeed, write $x=\sum_{i=1}^{m}(r_{i} \otimes
\frac{\alpha_{i}}{\beta_{i}})$ where $r_{i} \in R$, $\alpha_{i}
\in C$ and $\beta_{i} \in S$. Let $\beta=\Pi_{i=1}^{m}\beta_{i}$
and $\alpha_{i}'=\alpha_{i}\beta \beta_{i}^{-1} \in C$. Thus,
$$\sum_{i=1}^{m}(r_{i} \otimes \frac{\alpha_{i}}{\beta_{i}})=
\sum_{i=1}^{m}(r_{i} \otimes \frac{\alpha_{i}'}{\beta})=
\sum_{i=1}^{m} (\alpha_{i}' r_{i} \otimes \frac{1}{\beta})=r
\otimes \frac{1}{\beta}.$$

Here $r=\sum_{i=1}^{m} \alpha_{i}' r_{i}$.\enddemo

We now consider the tensor product $R \otimes_{O_{v}}F$ where $R$
is a torsion free algebra over $O_{v}$. Our goal is to construct a
quasi-valuation on $R \otimes_{O_{v}}F$ using the filter
quasi-valuation induced by $(R,v)$ that was constructed in Theorem
9.19.

\demo{Remark 9.30} Note that if $R$ is a torsion free algebra over
$O_{v}$, then there is an embedding $R \hookrightarrow R
\otimes_{O_{v}}F$; we shall see that in this case the
quasi-valuation on $R \otimes_{O_{v}}F$ extends the
quasi-valuation on $R$.\enddemo

%Therefore, in case $R$ is a torsion free algebra over $O_{v}$,
%every element $\sum_{i=1}^{m}r_{i} \otimes
%\frac{\alpha_{i}}{\beta_{i}} \in R \otimes_{O_{v}}F$ can be
%written {\bf uniquely} as $\sum_{i=1}^{m} r_{i}' \otimes
%\frac{1}{\beta}$.\enddemo

\proclaim{Lemma 9.31} Let $v, F, \Gamma_{v}$ and $O_{v}$ be as in
Theorem 9.19. Let $R$ be a torsion free algebra over $O_{v}$, $S$
a multiplicative closed subset of $O_{v}$, $0 \notin S$, and let
$w: R \rightarrow M \cup \{ \infty \}$ be any quasi-valuation
where $M$ is any totally ordered abelian monoid containing
$\Gamma_{v}$ and $w(cx)=v(c)+w(x)$ for every $c \in O_{v}$, $x \in
R$. Then there exists a quasi-valuation $W$ on $R
\otimes_{O_{v}}O_{v}S^{-1}$, extending $w$ on $R$ (under the
identification of $R$ with $R \otimes_{O_{v}} 1$), with value
monoid $ M \cup \{ \infty \}$.\endproclaim

\demo{Proof} In view of Remark 9.29, let $r \otimes
\frac{1}{\beta} \in R \otimes_{O_{v}}O_{v}S^{-1}$ and define
$$W( r\otimes \frac{1}{\beta})=w(r)-v(\beta)\ (\ =w(r)+(-v(\beta))\ ).$$
Note that $W$ is well defined since if $r \otimes \frac{1}{\beta}=
s \otimes \frac{1}{\delta}$ then there exists an $\alpha \in
O_{v}$ such that $\alpha (\delta r- \beta s)=0$ and thus, since
$R$ is torsion free, $\delta r=\beta  s$. Therefore, by our
assumption that $w(cx)=v(c)+w(x)$ for every $c \in O_{v}$ and $x
\in R$, we have
$$v(\delta)+w(  r)=v(\beta)+w(  s);$$ i.e., $W(r \otimes
\frac{1}{\beta})=W( s \otimes \frac{1}{\delta})$.

%$\sum_{i=1}^{m} r_{i} \otimes \frac{1}{\beta} \in R
%\otimes_{O_{v}}F$ and define
%$$W(\sum_{i=1}^{m} r_{i} \otimes \frac{1}{\beta})=w(\sum_{i=1}^{m}
%r_{i})-v(\beta).$$
%Note that $W$ is well defined since if
%$\sum_{i=1}^{m} r_{i} \otimes \frac{1}{\beta}=\sum_{j=1}^{k} s_{j}
%\otimes \frac{1}{\delta}$ then there exist an $\alpha \in O_{v}$
%such that $\alpha (\delta \sum_{i=1}^{m} r_{i}- \beta
%\sum_{j=1}^{k} s_{j})=0$ and thus, since $R$ is torsion free,
%$\delta \sum_{i=1}^{m} r_{i}=\beta \sum_{j=1}^{k} s_{j}$.
%Therefore, by Lemma 4.19, $v(\delta)+w( \sum_{i=1}^{m}
%r_{i})=v(\beta)+w( \sum_{j=1}^{k} s_{j})$; i.e., $W(\sum_{i=1}^{m}
%r_{i} \otimes \frac{1}{\beta})=W(\sum_{j=1}^{k} s_{j} \otimes
%\frac{1}{\delta})$.

We prove now that $W$ satisfies the axioms of a quasi-valuation.
First note that $W(0 \otimes 1)=w(0)-v(1)=\infty$. Next, note that
for every two elements $r \otimes \frac{1}{\beta} , s \otimes
\frac{1}{\delta} \in R \otimes_{O_{v}}O_{v}S^{-1}$, assuming that
$v(\beta) \leq v(\delta)$, we have $\delta= \alpha \beta$ for some
$\alpha \in O_{v}$ and thus $$r \otimes \frac{1}{\beta} = r
\otimes \frac{\alpha}{\alpha\beta}= \alpha r \otimes
\frac{1}{\delta}.$$ Therefore, we may assume that we have elements
$r \otimes \frac{1}{\delta} , s \otimes \frac{1}{\delta} \in R
\otimes_{O_{v}}O_{v}S^{-1}$; then $$W(r \otimes \frac{1}{\delta}+s
\otimes \frac{1}{\delta})=W((r+s)\otimes \frac{1}{\delta})$$
$$=w(r+s)-v(\delta) \geq \min \{ w(r), w(s) \} -v(\delta)$$ $$= \min \{
W(r \otimes \frac{1}{\delta}),W(s \otimes \frac{1}{\delta}) \}$$

Now, let $r \otimes \frac{1}{\beta} , s \otimes \frac{1}{\delta}
\in R \otimes_{O_{v}}O_{v}S^{-1}$; then $$W(r \otimes
\frac{1}{\beta} \cdot s \otimes \frac{1}{\delta})=W(rs \otimes
\frac{1}{\beta \delta})$$
$$=w(rs)-v(\beta \delta) \geq w(r)+w(s)-v(\beta)-v(\delta)$$ $$ =W(r
\otimes \frac{1}{\beta})+W(s \otimes \frac{1}{\delta}).$$

Finally note that $R$ embeds in $R \otimes_{O_{v}}O_{v}S^{-1}$ and
$$W(r \otimes 1)=w(r)-v(1)=w(r).$$

Q.E.D.\enddemo

\proclaim{Theorem 9.32} Let $v, F, \Gamma_{v}, O_{v}$ and $\Cal M
(\Gamma_{v})$ be as in Theorem 9.19. Let $R$ be a torsion free
algebra over $O_{v}$ and let $w$ denote the filter quasi-valuation
induced by $(R,v)$; then there exists a quasi-valuation $W$ on $R
\otimes_{O_{v}}F$, extending $w$ on $R$, with value monoid $\Cal
M(\Gamma_{v}) \cup \{ \infty \}$ and $O_{W}=R \otimes_{O_{v}}
1$.\endproclaim

\demo{Proof} Note that by Lemma 9.25, $w(cx)=v(c)+w(x)$ for every
$c \in O_{v}$, $x \in R$. Thus we can use Lemma 9.31 by taking
$S=O_{v} \setminus \{ 0 \}$ (and get $F=O_{v}S^{-1}$). Also note
that by Remark 9.29, every $x \in R \otimes_{O_{v}}F$ is of the
form $r \otimes \frac{1}{\beta}$ for $r \in R$ and $\beta \in
O_{v}$ (by taking $C=O_{v}$ and $S=O_{v} \setminus \{ 0 \}$).

So there exists a quasi-valuation $W$ on $R \otimes_{O_{v}}F$,
extending $w$ on $R$, with value monoid $ \Cal M(\Gamma_{v}) \cup
\{ \infty \}$. $W$ is given by $$W( r\otimes
\frac{1}{\beta})=w(r)-v(\beta)$$ for every $r \otimes
\frac{1}{\beta} \in R \otimes_{O_{v}}F$.

%$\sum_{i=1}^{m} r_{i} \otimes \frac{1}{\beta} \in R
%\otimes_{O_{v}}F$ and define
%$$W(\sum_{i=1}^{m} r_{i} \otimes \frac{1}{\beta})=w(\sum_{i=1}^{m}
%r_{i})-v(\beta).$$
%Note that $W$ is well defined since if
%$\sum_{i=1}^{m} r_{i} \otimes \frac{1}{\beta}=\sum_{j=1}^{k} s_{j}
%\otimes \frac{1}{\delta}$ then there exist an $\alpha \in O_{v}$
%such that $\alpha (\delta \sum_{i=1}^{m} r_{i}- \beta
%\sum_{j=1}^{k} s_{j})=0$ and thus, since $R$ is torsion free,
%$\delta \sum_{i=1}^{m} r_{i}=\beta \sum_{j=1}^{k} s_{j}$.
%Therefore, by Lemma 4.19, $v(\delta)+w( \sum_{i=1}^{m}
%r_{i})=v(\beta)+w( \sum_{j=1}^{k} s_{j})$; i.e., $W(\sum_{i=1}^{m}
%r_{i} \otimes \frac{1}{\beta})=W(\sum_{j=1}^{k} s_{j} \otimes
%\frac{1}{\delta})$.

%Now, let $r \otimes \frac{1}{\beta} , s \otimes \frac{1}{\delta}
%\in R \otimes_{O_{v}}F$; then $$W(r \otimes \frac{1}{\beta} \cdot
%s \otimes \frac{1}{\delta})=W(rs \otimes \frac{1}{\beta \delta})$$
%$$=w(rs)-v(\beta \delta) \geq w(r)+w(s)-v(\beta)-v(\delta)$$ $$ =W(r
%\otimes \frac{1}{\beta})+W(s \otimes \frac{1}{\delta}).$$

Finally note that for every element $r \in R$, we have, by Remark
9.18, that $w(r) \geq 0$ and thus $W(r \otimes 1)=w(r) \geq 0$. On
the other hand, let $r \otimes \frac{1}{\beta}  \in R
\otimes_{O_{v}}F$ with $W(r \otimes \frac{1}{\beta}) \geq 0$; then
$w(r) \geq v( \beta)$ i.e., $\beta \in S_{r}$ and thus one can
write $r=\beta r'$ for some $r' \in R$. Hence, $$r \otimes
\frac{1}{\beta}=\beta r' \otimes \frac{1}{\beta}=r' \otimes 1.$$
Consequently, $O_{W}=R \otimes 1$.

%Now, assume that $W(r \otimes \frac{1}{\beta}) \leq W( s \otimes
%\frac{1}{\delta})$ and also assume that $v(\beta) \leq v(\delta)$.
%Then $\delta= \alpha \beta$ for some $\alpha \in O_{v}$ and we
%have $r \otimes \frac{1}{\beta} = r \otimes
%\frac{\alpha}{\alpha\beta}= \alpha r \otimes \frac{1}{\delta}$.
%Thus $w(\alpha r)-v(\delta) \leq w(s)-v(\delta)$. Now, $W(r
%\otimes \frac{1}{\beta}+s \otimes \frac{1}{\delta})=W(\alpha r
%\otimes \frac{1}{\delta}+s \otimes \frac{1}{\delta})=W((\alpha r
%+s)\otimes \frac{1}{\delta})=w(\alpha r +s)-v(\delta) \geq \min \{
%w(\alpha r),w(s) \}-v(\delta)= w(\alpha r)-v(\delta)=W(r \otimes
%\frac{1}{\beta})=\min \{W(r \otimes \frac{1}{\beta}) , W( s
%\otimes \frac{1}{\delta}) \}$

Q.E.D.\enddemo

%$R=O_{w}$ and $w$ extends $v$ (on $F$).

\demo{Remark 9.33} Let $R$ be an algebra over $O_{v}$ and let $0
\neq r \in R$. By Remark 9.18, $w(r) \geq 0 $. Thus, for every $r
\otimes \frac{1}{\beta} \in R \otimes_{O_{v}}F$ where $r \neq 0$,
$$W(r \otimes \frac{1}{\beta})=w(r)-v(\beta) \geq-v(\beta)$$ i.e., $(W(r \otimes
\frac{1}{\beta}))^{L} \neq \emptyset$. Note that $(W(0 \otimes
\frac{1}{\beta}))=\infty-v(\beta)=\infty$.

Hence, $$(\emptyset, \Gamma_{v}) \notin \text{im}(W).$$

\enddemo

%$$(W(r \otimes
%\frac{1}{\beta}))^{L}=v(S_{r})+(-\infty,-v(\beta)] \supseteq
%(-\infty,-v(\beta)],$$

\proclaim{Theorem 9.34} Let $v, F, \Gamma_{v}, O_{v}$ and $\Cal
M(\Gamma_{v})$ be as in Theorem 9.19 and let $A$ be an
$F-$algebra. Let $R$ be a subring of $A$ such that $R \cap F
=O_{v}$. Then there exists a quasi-valuation $W$ on $RO_{v}^{-1} =
\{ rs^{-1} | r \in R, s \in O_{v} \setminus \{ 0 \} \}$ with value
monoid $\Cal M (\Gamma_{v})\cup \{ \infty \}$ such that $R=O_{W}$
and {\bf $ W $ extends $v$ (on $F$)}.
\endproclaim

\demo{Proof} Viewing $R$ as an algebra over $O_{v}$, $R$ has the
filter quasi-valuation defined in Theorem 9.19. Note that
$RO_{v}^{-1} \cong R \otimes_{O_{v}}F$ and thus, by Theorem 9.32,
there exists a quasi-valuation $W$ on $RO_{v}^{-1}$ such that
$R=~O_{W}$.

We shall now prove that $W$ extends $v$. Note that if $ 0 \neq x
\in O_{v}$ then $x \in S_{x}$ and thus $v(S_{x}) \supseteq
[0,v(x)]$ i.e., $w(x) \geq v(x)$. Moreover, for every $a \in
S_{x}$ one has $v(a) \leq v(x)$. (Indeed, if $x=ar$ for some $r
\in R$ then $xa^{-1}=r \in O_{v}$; i.e., $x \in aO_{v}$).
Therefore $ v(S_{x}) \subseteq [0,v(x)]$ and $w(x) = v(x)$. Now,
if $x \in F \setminus O_{v}$ then $x=\frac{\alpha}{\beta} $ where
$\alpha,\beta \in O_{v}$ and by the definition of $W$, we have
$$W(\frac{\alpha}{\beta})=W(\alpha \otimes
\frac{1}{\beta}) =w(\alpha)-v(\beta)$$
$$=v(\alpha)-v(\beta)=v(\frac{\alpha}{\beta})$$

The third equality is since $\alpha \in O_{v}$ and as proven
before $w(\alpha)=v(\alpha)$.

Q.E.D.\enddemo

%which is equal to $v(\alpha)-v(\beta)$ since
%$w(\alpha)=v(\alpha)$.

%We define $w:RO_{v}^{-1} \rightarrow \Cal M \cup \{ \infty \}$ by
%$$w(x) = w(rs^{-1}) = w(r)-w(s)$$ for every $ x \in RO_{v}^{-1}$ .
%Note that every $0 \neq s \in F$ is stable and thus $w$ is
%well-defined. Also note that if $x \in F$ then $x = rs^{-1}$ for
%$r,s \in O_{v}$, so $w(x) = w(rs^{-1}) = w(r)-w(s)$. However $w(r)
%= v(r)$ and $w(s) = v(s)$ since $r,s \in O_{v}$. So, $w(x)=v(x)$.
%Now, letting $r_{1}s_{1}^{-1}, r_{2}s_{2}^{-1} \in RO_{v}^{-1}$,
%we have
%$$w(r_{1}s_{1}^{-1}r_{2}s_{2}^{-1}) =
%w(r_{1}r_{2}s_{1}^{-1}s_{2}^{-1}) =
%w(r_{1}r_{2}(s_{2}s_{1})^{-1})$$  $$= w(r_{1}r_{2})-w(s_{2}s_{1})
%\geq w(r_{1})+w(r_{2})-(w(s_{2})+w(s_{1}))$$ $$=
%w(r_{1}s_{1}^{-1})+w(r_{2}s_{2}^{-1}).$$ Now assume
%$w(r_{1}s_{1}^{-1}) \leq w(r_{2}s_{2}^{-1})$; then
%$w(r_{1})+w(s_{2}) \leq w(r_{2})+w(r_{1})$. Therefore,

 %$$w(r_{1})+w(s_{2}) =
%w(r_{1}s_{2}) \leq w(r_{2}s_{1}) = w(r_{2})+w(s_{1}).$$ So, we
%have $$w(r_{1}s_{1}^{-1}+r_{2}s_{2}^{-1}) =
%w((r_{1}s_{2}+r_{2}s_{1})(s_{2}s_{1})^{-1})$$ $$=
%w((r_{1}s_{2}+r_{2}s_{1})-w((s_{2}s_{1}))) \geq
%w(r_{1}s_{2})-w(s_{2}s_{1})$$ $$=
%w(r_{1})+w(s_{2})-w(s_{2})-w(s_{1}) = w(r_{1}s_{1}^{-1})$$ $$= \min
%\{ w(r_{1}s_{1}^{-1}), w(r_{2}s_{2}^{-1}) \}.$$

%We now show $O_{w}=R$. Note that $R \subseteq O_{w}$ by Theorem
%4.26. As for the other direction, If $rs^{-1} \in O_{w}$ then
%$w(rs^{-1}) \geq 0$ i.e., $w(r) \geq w(s) =v(s)$. Therefore
%$v(S_{r}) \supseteq \Cal I_{v(s)}$ i.e., $v(s) \in v(S_{r})$ and
%thus $s \in S_{r}$. Namely, $r \in sR$, i.e., $rs^{-1} \in R$.

We continue our discussion in case $R$ is an integral domain. We
extend $w$ to $E$, the field of fractions of $R$. We assume that
$E$ is finite dimensional over $F$. Note that $RO_{v}^{-1} = \{
rs^{-1} | r \in R, s \in O_{v} \setminus \{ 0 \} \}$ is an
integral domain contained in $E$; i.e., $RO_{v}^{-1}$ is an
integral domain finite dimensional over $F$ and is thus a field.
Also note that $RO_{v}^{-1}$ is a field containing $R$ and thus
contains its field of fractions. Therefore $RO_{v}^{-1} = E$. So
every $x \in E$ can be written as $rs^{-1}$ where $r \in R, s \in
O_{v}$. We have the following important theorem:

\proclaim{Theorem 9.35} In view of Theorem 9.34, the
quasi-valuation $W$, as defined in Theorem 9.32, is a
quasi-valuation on $E$ extending $v$, with \linebreak $R=O_{W}$.
In other words, if $R \subseteq E$ satisfies $R \cap F=O_{v}$ and
$E$ is the field of fractions of $R$, then $R$ and $v$ induce a
quasi-valuation $W$ on $E$ such that $R=O_{W}$ and $W$ extends
$v$.
\endproclaim

Note: $W$ as described in Theorem 9.35 will also be called the
filter quasi-valuation induced by $(R,v)$.

Note that the filter quasi-valuation induced by $(R,v)$ is not
necessarily an exponential quasi-valuation. However, we cannot
hope for an exponential quasi-valuation in light of the following
example.

\demo{Example 9.36} Let $v$ denote a $p$-adic valuation on $ F=
\Bbb Q$; let $E=\Bbb Q[i]$ and $R=O_{v}[pi]$. Let $w$ denote a
quasi-valuation extending $v$ with $O_{w}=R$; then we must have
$w(i)<0$ (since $i \notin R$) whereas $w(i^{2})=w(-1)=0$.
% we have $w(i)=-1$ while $w(-1)=0$.
\enddemo

 \demo{Example 9.37} Recall from Lemma 2.8 that if $w$ is
a quasi-valuation on a field $E$ extending $v$ on $F$ and
$[E:F]<\infty$ then for any nonzero $x \in E$, $w(x)<\alpha$ for
some $\alpha \in \Gamma_{v}$. If $\Gamma_{v}=\Bbb Z$ and $R
\subseteq E$ is a
ring such that $R \cap F=O_{v}$, %$v$ is a discrete rank one valuation (i.e.,
 then by Example 9.13 and Remark 9.33, %so is the filter
%quasi-valuation induced by $v$ (i.e.,
$M_{w}$ can be identified as $\Bbb Z$; where $M_{w}$ is the value
monoid of the filter quasi-valuation induced by $(R,v)$. Namely,
the filter quasi-valuation induced by $(R,v)$ is a quasi-valuation
extending $v$ with $M_{w}$ a group; therefore $R=O_{w}$ has also
the properties of a quasi-valuation ring studied in
sections 5 and 6.% one can apply on $O_{w}$ the methods developed in [Sa].
\enddemo

We summarize the main results we have obtained using the theory of
quasi-valuation extending a valuation on a finite field extension.

\proclaim{Theorem 9.38} Let $v$ be a valuation on a field $F$ with
a valuation ring $O_{v}$ and value group $\Gamma_{v}$. Let $R$ be
an integral domain with field of fractions $E$ finite dimensional
over $F$ and $R \cap F=O_{v}$. Then:

1. Every f.g. ideal of $R$ is generated by $m \leqslant [E:F]$
generators.

(See Corollary 2.4).

2. There exists a quasi-valuation $W$ on $E$ extending $v$ such
that $R=O_{W}$.

(See Theorem 9.35)

3. $R$ satisfies LO, INC and GD over $O_{v}$.

(See Lemma 3.12, Remark 3.5 and Theorem 3.7, and Lemma 4.12

respectively).

4. K-dim$R =$K-dim$O_{v}$. (See Corollary 4.13]).

5. If $R$ satisfies GU over $O_{v}$ then $R$ satisfies the height
formula.

\ \ \ \ (See Theorem 4.17).

6. If there exists a quasi-valuation $w$ on $E$ extending $v$ with
$R=O_{w}$ such that

\ \ \ \ $w(E \setminus \{0\})$ is torsion over $\Gamma_{v}$, then

  \ \ \ \ (a). $R$ satisfies GU over $O_{v}$. (See Theorem
  5.16).

 % \ \ \ \ (b).

  \ \ \ \ (b).  K-dim$O_{v} \leq \mid \text{Spec}(R) \mid \leq \mid \text{Spec}(I_{E}(R)) \mid
\leq [E:F]_{sep} \cdot$K-dim$O_{v}$.

 \ \ \ \ (See Theorem
  5.21).

   \ \ \ \ (c). $R$ has finitely many maximal ideals, the number of which is

 \ \ \ \ less or equal to $[E:F]$. In fact, for each $P \in \text{Spec}(O_{v})$
there are at most

 \ \ \ \ $[E:F]$ prime ideals $Q \in \text{Spec}(R)$ lying over $P$. (See
Theorem 5.19 and

 \ \ \ \  Theorem 5.21 and the discussion before
Theorem 5.21).
\endproclaim

Note, for example, that in view of Example 9.37, if $F=\Bbb Q$ %is a number field
and $R$ is as above, then $R$ satisfies properties 1-6.

%8. The spectrum of $R$ is nicely ordered.

\heading \S 10 Properties of the filter quasi-valuation
\endheading

In this section we prove some properties of the cut monoid (Lemma
10.1 and Lemma 10.6). These properties are valid in general (for
the cut monoid induced by an arbitrary totally ordered abelian
group). In addition, we prove some facts regarding filter
quasi-valuations extending a valuation on a finite field
extension.

\proclaim{Lemma 10.1} Let $\Gamma$ be a totally ordered abelian
group; then the only PIM in the cut monoid $\Cal M(\Gamma)$ lying
over the isolated subgroup $\{0\}$ is the set
$\{((-\infty,0],(0,\infty))\}$; i.e., the set containing the $0$
of $\Cal M(\Gamma)$.
%$\{\{0\}\}$
\endproclaim

\demo{Proof} First note that $hull_{\Cal
M(\Gamma)}({\{0\}})=\{((-\infty,0],(0,\infty))\}$ lies over
$\{0\}$. Now, if there exists another PIM $N \neq hull_{\Cal
M(\Gamma)}({\{0\}})$ lying over $\{0\}$ then take an element
$((-\infty,0],(0,\infty)) \neq \Cal A \in N$. Next, take $ 0 <
\alpha \in \Cal A^{L}$ and get $N \cap \Gamma \neq \{ 0 \}$.

Q.E.D.\enddemo

Let $R$ be a ring; we denote by $J(R)$ its jacobson radical. Note
that if $C \subseteq R$ are commutative rings and $R$ satisfies GU
over $C$ then every maximal ideal of $R$ lies over a maximal ideal
of $C$. In particular, if $R$ satisfies GU over $C$ and $C$ is
local then $J(R) \supseteq J(C)$.

\proclaim{Proposition 10.2} Let $v$ be a valuation on a field $F$
and let $E/F $ be a finite field extension. Let $R $ be a subring
of $ E$ satisfying $R \cap F=O_{v}$ and $J(R) \supseteq I_{v}$.
Let $w$ be the filter quasi-valuation induced by $(R,v)$ (and thus
$R=O_{w}$). Then
$$\sqrt{I_{w}} = J_{w}(=J(R)).$$
\endproclaim

\demo{Proof} Apply proposition 4.19 by taking $P = I_{v}$; then
use Lemma 10.1, the assumption that $J(R) \supseteq I_{v}$ (to
deduce $\sqrt{I_{w}} \subseteq J_{w}$) and the fact that $O_{w}$
satisfies INC over $O_{v}$ (to deduce $\sqrt{I_{w}} \supseteq
J_{w}$).

Q.E.D.\enddemo

For example, if $\Gamma_{v}= \Bbb Z$ then, by Example 9.37,
$M_{w}$ can be identified as $ \Bbb Z$ (where $w$ denotes the
filter quasi-valuation induced by $(R,v)$) and thus, by Theorem
5.16, $R$ satisfies GU over $O_{v}$ and therefore $J(R) \supseteq
I_{v}$. Hence, by Proposition 10.2, $\sqrt{I_{w}} =J(R)$.

 \proclaim{Lemma 10.3} Let $v$ be a valuation on a field $F$ and
let $E/F $ be a finite field extension. Let $R $ be a subring of $
E$ satisfying $R \cap F=O_{v}$ and let $w$ be the filter
quasi-valuation induced by $(R,v)$ (and thus $R=O_{w}$). Then
$I_{w}=I_{v}R$. \endproclaim

\demo{Proof} Let $x \in I_{w}$; then there exists $a \in O_{v}$
with $v(a)>0$ such that $x \in aR$. Thus, $x \in aR \subseteq
I_{v}R$. On the other hand, let $x \in I_{v}R$; then one can write
$x=ar$ for $a \in I_{v}$ and $r \in R$. Thus $w(x) \geq
w(a)=v(a)>0$ and $x \in I_{w}$.
\enddemo

%\proclaim{Corollary 10.3} Notation as in Proposition 10.2. Then

%${I_{w}} \subseteq M$ for every maximal ideal $M$ of $O_{w}$.
%\endproclaim

%Compare Corollary 10.3 with [Sa, Lemma 5.3] where we prove this
%property under the stronger hypothesis that $M_{w}$ is a group.

We note now that since $\Cal M(\Gamma)$ is $\Bbb N$-strictly
ordered, by Lemma 9.12, then one can embed $\Cal M(\Gamma)$ in its
divisible hull $\Cal M_{\text{div}}=(\Cal M(\Gamma) \times \Bbb N
) / \sim $ where $\sim$ is the equivalence relation defined by
$$(m_{1}, n_{1}) \sim (m_{2}, n_{2}) \text{\ iff\ } n_{2}m_{1}=n_{1}m_{2}.$$
Note that $\Gamma_{\text{div}} $ and $\Cal M(\Gamma)$ embed in
this divisible hull.

Now, we aim for a stronger version of Theorem 8.15 when dealing
with filter quasi-valuations. This time we may consider the
ordering inside the divisible hull of $\Cal M(\Gamma)$. We start
with the following lemma:

\proclaim{Lemma 10.4} Notation as in Lemma 10.3. Let $u$ be a
valuation on $E$ extending $v$ such that $O_{u} \supseteq O_{w}$;
then $I_{u} \supseteq I_{w}$.\endproclaim

\demo{Proof} By Lemma 10.3, $I_{w}=I_{v}R$. Thus,
$$I_{u}=I_{u}O_{u} \supseteq I_{v}O_{u} \supseteq I_{v}R =I_{w}.$$

%As in [Sa, Lemma 3.23]. $K = I_{u} \cap O_{w}$ is a prime ideal of
%$O_{w}$. Assuming $K$ is not maximal in $O_{w}$, one has $K \cap
%O_{v}$ not maximal in $O_{v}$ (by INC). This contradicts the fact
%that $I_{u} \cap O_{v} = I_{v}$. So $I_{u} \cap O_{w}$ is a
%maximal ideal in $O_{w}$ and thus (Corollary 10.3) contains
%$I_{w}$.

Q.E.D.\enddemo

\proclaim{Theorem 10.5} Notation as in Lemma 10.3. Then there
exists a valuation $u$ of $E$ extending $v$ on $F$ such that
$O_{u} \supseteq O_{w}$. Moreover, for every such $u$, $u$
dominates $w$.\endproclaim

\demo{Proof} The first part is true for any quasi-valuation, as
proved in Lemma 8.13. As for the second part, note that by Lemma
10.4, $O_{u} \supseteq O_{w}$ implies $I_{u} \supseteq I_{w}$ and
apply Theorem 8.14.

Q.E.D.\enddemo

\proclaim{Lemma 10.6} Let $H$ be an isolated subgroup of a totally
ordered abelian group $\Gamma$. Then there exist at most two PIMs
in $\Cal M(\Gamma)$ lying over $H^{\geq0}$; namely, $$hull_{\Cal
M(\Gamma)}(H^{\geq0}) \text{ and } \overline{hull_{\Cal
M(\Gamma)}(H^{\geq0})}.$$\endproclaim

\demo{Proof} First note that $H^{+} \in \Cal M(\Gamma)$ and
$(H^{+})^{L}=(-\infty,0] \cup H^{\geq 0}$. Now, if $H=\{0\}$ then
by Lemma 10.1 there is only one PIM lying over $H=H^{\geq0}$. So,
we assume that $H \neq \{ 0 \}$. We prove that
$$X= \{ \Cal A \in \Cal M(\Gamma) \mid 0 \leq \Cal A < H ^{+}\}
$$ is equal to $hull_{\Cal M(\Gamma)}(H^{\geq0})$. We denote %by $\Cal C$ the cut $(H^{\geq0}
%\cup \Gamma^{<0}, \Gamma \setminus (H^{\geq0} \cup \Gamma^{<0}))$
%and
 $$Y= hull_{\Cal M(\Gamma)}(H^{\geq0}) \cup \{ H ^{+} \}= \{ \Cal A \in \Cal M(\Gamma) \mid 0 \leq \Cal A \leq  H ^{+}\}.$$ We shall prove
that $Y$ is equal to $\overline{hull_{\Cal
M(\Gamma)}(H^{\geq0})}$. Indeed if $\Cal A \in X$ then there
exists $\alpha \in H^{\geq0}$ such that $\alpha \notin \Cal A^{L}
$ and thus $$(-\infty,0] \subseteq \Cal A^{L} \subset
(-\infty,\alpha];$$ i.e., $\Cal A \in hull_{\Cal
M(\Gamma)}(H^{\geq0})$. On the other hand, if $\Cal A \in
hull_{\Cal M(\Gamma)}(H^{\geq0})$ then $$(-\infty,0] \subseteq
\Cal A^{L} \subseteq (-\infty,\alpha]$$ for $0 \neq \alpha \in
H^{\geq0}$; thus for example, $\Cal A^{L} \subset
(-\infty,2\alpha]$ i.e., $(\Cal A^{L})^{\geq 0} \subset H ^{\geq
0} $. Thus $\Cal A \in X$.

Next, if $\Cal A \in Y$ then obviously $\Cal A \in
\overline{hull_{\Cal M(\Gamma)}(H^{\geq0})}$. Now, assuming there
exists $\Cal B \in \overline{hull_{\Cal M(\Gamma)}(H^{\geq0})}
\setminus Y$, then $(\Cal B^{L})^{\geq 0} \supset H ^{\geq 0}$;
take $\beta \in (\Cal B^{L})^{\geq 0} \setminus H^{\geq0}$ then
$(-\infty,\beta] \supset H^{\geq0}$ and thus
$$\overline{hull_{\Cal M(\Gamma)}(H^{\geq0})} \cap \Gamma \supset H^{\geq0},$$ a
contradiction. Therefore $hull_{\Cal M(\Gamma)}(H^{\geq0})$ and
$\overline{hull_{\Cal M(\Gamma)}(H^{\geq0})}$ differ only by one
element and thus there are no other PIMs lying over $H^{\geq0}$.

Q.E.D.\enddemo

%Note: one can visualize the sets $X$ and $Y$ in the proof above as
%$X= \{ \Cal A \in \Cal M(\Gamma) \mid 0 \leq \Cal A < \Cal C\} $
%and $Y= \{ \Cal A \in \Cal M(\Gamma) \mid 0 \leq \Cal A \leq \Cal
%C\} $.

Here is an example of a group contained in a monoid for which
there are more than two PIMs lying over the positive part of an
isolated subgroup:

 \demo{Example 10.7} We consider the set $\Bbb N \cup \{ 0\}$
with the maximum operation, i.e., $i+j=\max\{ i,j\}$ for all $i,j
\in \Bbb N \cup \{ 0\}$. Let $\Bbb Z \times (\Bbb N \cup \{ 0\})$
denote the totally ordered abelian monoid with addition defined
componentwise and the left to right lexicographic order, i.e.,

$$(z,i) \leq (z',j)\text{\ iff\ } z < z' \text{\ or\ }
(z=z' \text{\ and\ } i \leq j).$$ Viewing $\Bbb Z$ inside $\Bbb Z
\times (\Bbb N \cup \{ 0\})$ via the natural monomorphism $z
\rightarrow (z,0)$, we see that there is an infinite number of
PIMs lying over $\{0\}$. Namely, $$\{ (0,i ) \}_{i \leq j}$$ for
every $j \in \Bbb N \cup \{ 0\}$.

%containing $\Gamma$ and there are an infinite number of PIMs lying
%over $\{0\}$. Namely, $\{ 0\} , \{ m \in M \mid m \leq n \alpha
%_{i}$ for some $n \in \Bbb N$ and $\alpha _{i}\in M \}=\{ 0\} , \{
%0\} \cup \{ \alpha _{k}\}_{k \leq i}$ for every $i \in \Bbb
%N$.
\enddemo

We present now an example of a quasi-valuation such that inside
$M_{w}$ there are two PIMs lying over $\{ 0 \}$ ; in particular,
this is not a filter quasi-valuation.

\demo{Example 10.8} Let $v$ denote a p-adic valuation on $\Bbb Q$
and let $M=\{ \alpha _{0}, \alpha _{1}\}$ be the totally ordered
abelian monoid with the maximum operation where $\alpha _{1}>
\alpha _{0}$. Let $\Bbb Z \times M$ denote the totally ordered
abelian monoid defined in a similar way as in Example 10.7 and
adjoin a largest element $\infty$ (as usual). Define $w$ on $\Bbb
Q[\sqrt{p}]$ by

$$w(a+b \sqrt{p})= \cases \infty  \text{ \ \ \ \ \ \ \ \ \ \ \ \ \ \ \ \ \ \ \ \ \ \ \ \ \ \ \ \ \ \ \ \ \ \ \ if     } a=b=0
\\ (v(b),\alpha _{1}) \text{\ \ \ \ \ \ \ \ \ \ \ \ \ \ \ \ \ \ \ \ \ \ \ \ \ \ \ if     } a=0, b \neq 0
\\ (v(a),\alpha _{0}) \text{\ \ \ \ \ \ \ \ \ \ \ \ \ \ \ \ \ \ \ \ \ \ \ \ \ \ \  if     } b=0, a \neq 0
\\ \min \{ (v(a),\alpha _{0}), (v(b),\alpha _{1})\} \text{\ \ \ \ \ if }  a \neq 0, b \neq 0.
\endcases    $$

We note that the elements $(\infty,\alpha _{0})$ and
$(\infty,\alpha _{1})$ are not defined and thus the case
distinction above is needed.

It is not difficult to check that $w$ is a quasi-valuation on
$\Bbb Q[\sqrt{p}]$ extending $v$ with $M_{w}=\Bbb Z \times M$ and
2 PIMs over $\{ 0 \}$.

Note, for example, that $w(\sqrt{p})=(0,\alpha _{1}) \in
\overline{hull_{M_{w}} (\{ 0 \})} $ while $w(\sqrt{p}) > 0$ which
is the only element in $H=\{  0 \}$ (where obviously $0 \in \Bbb
Z$ is identified as $(0,\alpha _{0})$ in $M_{w}$).
\enddemo

%Here is an example of a monoid for which there are more then two
%PIMs over the positive part of an isolated subgroup:

%\demo{Example 10.7} Let $\Gamma= \Bbb Z$ and let $M$ be the monoid
%generated by $\Gamma \cup \{ \alpha _{i}\}_{i \in \Bbb N}$ where
%$M$ is commutative, each $0< \alpha _{i} <1$, $\alpha _{i}<\alpha
%_{j}$ for $i < j$ and $\alpha _{i}+\alpha _{j}=\alpha _{j}$ for $i
%\leq j$. Then $M$ is a totally ordered abelian monoid containing
%$\Gamma$ and there are an infinite number of PIMs lying over
%$\{0\}$. Namely, $\{ 0\} , \{ m \in M \mid m \leq n \alpha _{i}$
%for some $n \in \Bbb N$ and $\alpha _{i}\in M \}=\{ 0\} , \{ 0\}
%\cup \{ \alpha _{k}\}_{k \leq i}$ for every $i \in \Bbb
%N$.\enddemo

%We present now an example of a quasi-valuation such that inside
%$M_{w}$ there are two PIMs lying over $\{ 0 \}$ ; in particular,
%this is not a filter quasi-valuation.

%\demo{Example 10.8} Let $v$ denote a p-adic valuation on $\Bbb Q$
%and let $M$ be the monoid generated by $ \Bbb Z \cup \{ m\}$ where
%$M$ is abelian, $0< m< 1$ and $nm=m$ for every $n \in \Bbb N$.
%Define $w$ on $\Bbb Q[\sqrt{p}]$ by $w(a+b \sqrt{p})=\min \{ v(a),
%v(b)+m\}$. Then $w$ is a quasi-valuation on $\Bbb Q[\sqrt{p}]$
%extending $v$ with $M_{w}=M$ and 2 PIMs over $\{ 0 \}$.\enddemo

Now, we show that the filter quasi-valuation construction respects
localization at prime ideals of $O_{v}$.

%can be constructed in a somewhat different way when it extends a
%valuation whose valuation ring contains another valuation ring.

\demo {Remark 10.9} Let $v$ be a valuation on a field $F$, let
$E/F $ be a finite field extension and let $R \subseteq E$ be a
ring such that $E$ is the field of fractions of $R$ and $R \cap
F=O_{v}$. Let $w_{v}$ be the filter quasi-valuation induced by
$(R,v)$. Thus $w_{v}$ extends $v$ and $O_{w_{v}}=R$. Let $\Cal
M(\Gamma_{v})$ denote its cut monoid. Let $P \vartriangleleft
O_{v}$ be a prime ideal of $O_{v}$, and $H$ its corresponding
isolated subgroup. Then $(O_{v})_{P}=S^{-1}O_{v}$, where $S=O_{v}
\setminus P$. $S^{-1}O_{v}$ is a valuation ring containing
$O_{v}$; we denote it by $O_{u}$ and its valuation by $u$. Recall
that for every $x \in F$, $u(x)=v(x)+H$ and
$\Gamma_{u}=\Gamma_{v}/H$. Note that $S^{-1}R$ is a subring of $E$
such that $S^{-1}R \cap F= O_{u}$ and thus there exists a filter
quasi-valuation, denoted $w_{u}$, induced by $(S^{-1}R,u)$; i.e.,
$w_{u}$ extends $u$ and $O_{w_{u}}=S^{-1}R$. We denote its cut
monoid by $M(\Gamma_{u})=\Cal M(\Gamma_{v} / H)$. Recall that
every element $z \in E$ can be written as $xy^{-1}$ where $x \in
R$ and $y \in O_{v}$. Also recall that $\Gamma_{v}$ embeds in
$\Cal M(\Gamma_{v})$ and when we write, for $y \in O_{v}$, $v(y)$
as an element of $\Cal M(\Gamma_{v})$ we refer to the cut
$((-\infty,v(y)],(v(y),\infty))$; the same notation holds for the
valuation $u$.\enddemo

\proclaim{Lemma 10.10} Notation as in Remark 10.9. If $x \in R$
then
$$\{ v(a)+H \mid a \in S_{x}^{R/O_{v}}\}= \{ v(b)+H \mid b \in
S_{x}^{S^{-1}R/S^{-1}O_{v}}\}.$$\endproclaim

\demo{Proof} $(\subseteq)$ If $a \in S_{x}^{R/O_{v}}$ then $x \in
a R$ for $a \in O_{v}$ and obviously $x \in a S^{-1}R$ and $a \in
S^{-1}O_{v}$, thus $a\in S_{x}^{S^{-1}R/S^{-1}O_{v}}$.
$(\supseteq)$ If $b \in S_{x}^{S^{-1}R/S^{-1}O_{v}}$ then $x \in b
S^{-1}R$ for $b \in S^{-1}O_{v}$. We have two possibilities:

Case I. $v(b )< 0$, then $v(b)\in H$ (since otherwise $b \notin
S^{-1}O_{v}$) and we are

 done.

Case II. $v(b) \geq 0$, assume to the contrary that $v(b)>v(a)+h$
for every

$a \in S_{x}^{R/O_{v}}$ and $h \in H$; we have $x=bs^{-1}r$ for $s
\in S$, $r \in R$ and thus

$v(bs^{-1})>0$ (since $v(b)>v(a)+h \geq h$ for every $h \in H$)
i.e., $bs^{-1}\in O_{v}$.

Therefore, writing $a=bs^{-1} \in S_{x}^{R/O_{v}}$, we have a
contradiction

 (since $v(b)=v(a)+v(s)$).

Q.E.D.\enddemo

%\vfill \eject

%$$w(z)^{L}=w(xy^{-1})^{L}=\{ v(a)+H \mid a \in S_{x}^{R/O_{v}}
%\} \cup (\Gamma_{v} / H)^{<0}+ (-\infty,-v(y)+H],$$

 \proclaim{Theorem 10.11}  Let $w : E
\rightarrow \Cal M(\Gamma_{v} / H)$ be the function defined by

$$w(z)=w(xy^{-1})=\{ v(a)+H \mid a \in S_{x}^{R/O_{v}}
\}^{+} + (-v(y)+H),$$ $\forall z \in E$ where $ z=xy^{-1} \text{
for } x \in R \text{ and }  0 \neq y \in O_{v}$. Then $w$ is the
filter quasi-valuation induced by $(S^{-1}R,u)$; i.e., $w=w_{u}$.
\endproclaim

\demo{Proof} Let $z \in E$ and write $z=xy^{-1}$ for $x \in R$ and
$0 \neq y \in O_{v}$. Note that $y \in O_{v} \subseteq
S^{-1}O_{v}$ is stable with respect to $w_{u}$ and $x \in R
\subseteq S^{-1}R$; thus
$$w_{u}(z)=w_{u}(xy^{-1})=w_{u}(x)-w_{u}(y)=w_{u}(x)-u(y).$$ Now,

$$w_{u}(x)-u(y)=u(S_{x}^{S^{-1}R/S^{-1}O_{v}})^{+}-u(y)$$ $$
=\{ u(b) \mid b \in S_{x}^{S^{-1}R/S^{-1}O_{v}} \}^{+}-u(y).$$
Note that $b \in O_{u}$ and $ u(b)=v(b)+H$. Moreover, by Lemma
10.10, $$\{ v(b)+H \mid b \in S_{x}^{S^{-1}R/S^{-1}O_{v}} \}=\{
v(a)+H \mid a \in S_{x}^{R/O_{v}} \}.$$ Also, $y \in O_{v}$ and
$u(y)=v(y)+H$. Consequently, $w=w_{u}$.

Q.E.D.\enddemo

%So, we have the following commutative diagram:

%u(S_{x}^{S^{-1}R/S^{-1}O_{v}})

%$$\CD
%R @> >> S^{-1}R\\
%@A  AA  @AA  A\\
%O_{v} @> >> O_{u}=S^{-1}O_{v}
%\endCD
%$$

%We have the following commutative diagram:

%$$\CD
%R @>\alpha>> S^{-1}R\\
%@A\alpha AA  @AA\alpha A\\
%O_{v} @>\alpha>> O_{u}=S^{-1}O_{v}
%\endCD
%$$

%$$\CD
%\Cal M(\Gamma_{v}) @>\alpha>> \Cal M(\Gamma_{u})\\
%@A\alpha AA  @AA\alpha A\\
%\Gamma_{v} @>f(\gamma)=\gamma+H>> \Gamma_{u}=\Gamma_{v}/H
%\endCD
%$$

%$$\CD
%w_{v} @>\alpha>> w_{u}\\
%@A\alpha AA  @AA\alpha A\\
%v @>f(v(x))=v(x)+H>> u
%\endCD
%$$

% $\{ u(b) \mid \in S_{x}^{S^{-1}R/S^{-1}O_{v}}
%\}=\{ v(b)+H \mid \in S_{x}^{S^{-1}R/S^{-1}O_{v}} \}$ and by Lemma
%5.9, $\{ v(b)+H \mid \in S_{x}^{S^{-1}R/S^{-1}O_{v}} \}=\{ v(a)+H
%\mid a \in S_{x}^{R/O_{v}} \}$\enddemo

%\vfill \eject

We shall now present the minimality of the filter quasi-valuation
with respect to a natural partial order.

For every ring $R \subseteq E$ satisfying $R \cap F=O_{v}$, we
denote by $$\Cal W_{R}=\{ w \mid~ w \text{\ is a quasi-valuation
on\ } E \text{\ extending\ } v \text{\ with\ } O_{w}=R\}.$$

%and by $w_{f}$ its filter quasi-valuation.

Note that $\Cal W_{R}$ is not empty by Theorem 9.35.

%\heading \S 4.3  The minimality of the filter quasi-valuation
%\endheading

%For every ring $R \subseteq E$ satisfying $R \cap F=O_{v}$,

We define an equivalence relation on $\Cal W_{R}$ in the following
way: $w_{1} \sim w_{2}$ iff for every $x,y \in E$, $$w_{1}(x)<
w_{1}(y) \Longleftrightarrow w_{2}(x)<w_{2}(y).$$

 We define a partial order on $\Cal
W_{R} / \sim$ in the following way: $[w_{1}] \leq [w_{2}]$ iff for
every $x,y \in E$, $w_{1}(x)< w_{1}(y)$ implies
$w_{2}(x)<w_{2}(y)$ (we say that $w_{1}$ is coarser than $w_{2}$
and that $w_{2}$ is finer than $w_{1}$). Note that $\leq$ is well
defined and is indeed a partial order on the set of equivalent
quasi-valuations corresponding to $R$.

We shall prove now that the equivalence class of the filter
quasi-valuation is the minimal one with respect to the partial
order defined above.

%filter quasi-valuation is unique in the sense that its equivalence
%class is the minimal one with respect to the partial order
%defined above.

\proclaim{Proposition 10.12} The (equivalence class of \ the)
filter quasi-valuation is the coarsest of all (equivalence classes
of) quasi-valuations in $\Cal W_{R} / \sim$.
\endproclaim

\demo{Proof} We first prove that if $x,y \in R$ and $w \in \Cal
W_{R}$ then $ w_{v}(x)< w_{v}(y)$ implies $ w(x)< w(y)$ (where
$w_{v}$ is the filter quasi-valuation induced by $(R,v)$ and thus
$O_{w_{v}}=R$). Let $x,y \in R$, $w \in \Cal W_{R}$ and assume $0
\leq w_{v}(x)< w_{v}(y)$; then there exists an $\alpha \in
\Gamma_{v}$ such that
$$w_{v}(x)< \alpha \leq w_{v}(y).$$ Let $a \in F$ with
$v(a)=\alpha$; then $w_{v}(ya^{-1}) \geq 0$ and thus $ya^{-1} \in
R$ and $w(ya^{-1}) \geq 0$, i.e., $w(y) \geq v(a)$. Note that
$w(x)<v(a)$ since otherwise $xa^{-1} \in R$ and $w_{v}(x) \geq
v(a)$, a contradiction.

Now, in the general case, let $x,y \in E$, $w \in \Cal W_{R}$ and
assume that $ w_{v}(x)< w_{v}(y)$. Write
$$w_{v}(x)=\Cal A- \alpha,\ w_{v}(y)=\Cal A'- \alpha'$$  for $\Cal A,\Cal A' \in
\Cal M(\Gamma_{v})$, $\Cal A,\Cal A' \geq 0$ and $\alpha,\alpha'
\in (\Gamma_{v})^{\geq 0};$ then $0 \leq \Cal A < \Cal
A'-\alpha'+\alpha$. Let $a \in F$ with $v(a)=\alpha$ and define
$x'=xa$, $y'=ya$; thus
$$0 \leq w_{v}(x')=\Cal A< w_{v}(y')=\Cal A'-\alpha'+\alpha$$ and, by the
proof of the first part, $ w(x')< w(y')$. Therefore
$$w(x)=w(x')-v(a)<w(y')-v(a)=w(y).$$

Q.E.D.\enddemo

I would like to thank  S. Margolis, P. J. Morandi, J.-P. Tignol,
U. Vishne, and A.R. Wadsworth for their helpful comments and
queries. I would also like to thank the referee for his helpful
comments and suggestions.

Special thanks are, of course, to my advisor, L. H. Rowen, who
spent hours, with and without me, going through the ideas
presented in this paper. Rowen suggested many crucial corrections,
improvements and questions. I owe a very special debt to him.

%I would like to express special gratitude to my advisor, L. H.
%Rowen, for his valuable comments, suggestions, and corrections.

%\vfill \eject

 \Refs
\def\hangbox to #1 #2{\vskip1pt\hangindent #1\noindent \hbox to #1{#2}$\!\!$}
\def\refn#1{\hangbox to 40pt {#1\hfill}}

\refn{[Bo]} N. Bourbaki, {\it Commutative Algebra}, Chapter 6,
Valuations, Hermann, Paris, 1961.

\refn{[Co]} P. M. Cohn, {\it An Invariant Characterization of
Pseudo-Valuations}, Proc. Camp. Phil. Soc. 50 (1954), 159-177.

\refn{[End]} O. Endler, {\it Valuation Theory}, Springer-Verlag,
New York, 1972.

\refn{[FKK]} A. Fornasiero, F.V. Kuhlmann and S. Kuhlmann, {\it
Towers of complements to valuation rings and truncation closed
embeddings of valued fields}, J. Algebra 323 (2010), no. 3,
574-600.

\refn{[Hu]} J. A. Huckaba {\it Extensions of Pseudo-Valuations},
Pacific J. Math. Volume 29, Number 2 (1969), 295-302.

\refn{[Kap]} I. Kaplansky, {\it Commutative Rings}, Allyn and
Bacon, Boston, 1970.

\refn{[KZ]} M. Knebusch and D. Zhang, {\it Manis Valuations and
Pr\"ufer Extensions}, Springer-Verlag, Berlin, 2002.

\refn{[MH]} M. Mahadavi-Hezavehi, {\it Matrix Pseudo-Valuations on
Rings and their Associated Skew Fields}, Int. Math. J. 2 (2002),
no. 1, 7--30.

\refn{[Mor]} P. J. Morandi, {\it Value functions on central simple
algebras}, Trans. Amer. Math. Soc., 315 (1989), 605-622.

\refn{[Row]} L.H. Rowen, {\it Graduate Algebra: Commutative View},
American Mathematical Society, Rhode Island, 2006.

\refn{[Ste]} W. Stein, {\it A Brief Introduction To Classical And
Adelic Algebraic Number Theory}, 2004. Link:
http://modular.math.washington.edu/papers/ant/.

\refn{[TW]} J.-P. Tignol and A.R. Wadsworth, {\it Value functions
and associated graded rimgs for semisimple algebras}, to appear.

\refn{[Wad]} A.R. Wadsworth, {\it Valuation theory on finite
dimensional division algebras}, pp. 385-449 in {\it Valuation
Theory and its Applications}, Vol. 1, eds. F.-V. Kuhlmann et al.,
Fields Inst. Commun., 32, American Mathematical Society,
Providence, RI, 2002.

\refn{[Weh]} F. Wehrung, {\it Monoids of intervals of ordered
abelian groups}, J. Algebra 182 (1996), no. 1, 287-328.

\endRefs

\end